\documentclass{article}
\usepackage{amsmath,amssymb,amsfonts}
\usepackage[T2A]{fontenc}
\newtheorem{lemma}{Lemma}
\newtheorem{theorem}{Theorem}
\newtheorem{proposition}{Proposition}
\newtheorem{corollary}{Corollary}
\newtheorem{conjecture}{Conjecture}

\def\D{{\cal D}}
\def\R{{\mathbb R}}
\def\C{{\mathbb C}}
\def\Z{{\mathbb Z}}

\def\H{{\mathbb H}}
\def\W{{\cal W}}
\def\Re{{\mathrm{Re}\, }}
\def\Im{{\mathrm{Im}\, }}
\def\Spec{{\mathrm{Spec}\, }}
\def\vol{{\mathrm{vol}\,}}
\def\ord{{\mathrm{ord}\,}}
\def\nm{{\mathrm{nm}}}
\def\j{{\bf j}}
\def\const{{\mathrm{const}}}
\def\Tr{{\mathrm{Tr}\, }}
\def\nil{{\mathrm{Nil}\, }}
\def\sol{{\mathrm{Sol}\, }}
\def\sll{\widetilde{SL}_2}
\def\U{{\bf U}}
\def\V{{\bf V}}
\def\Ker{\mathrm{Ker}\,}

\def\beq{\begin{equation}}
\def\eeq{\end{equation}}
\def\bet{\begin{theorem}}
\def\eet{\end{theorem}}
\def\bi{\begin{itemize}}
\def\ei{\end{itemize}}

\begin{document}

\title{Two-dimensional Dirac operator and surface theory
\thanks{The work is supported by Russian Foundation for Basic
Researches.}}
\author{Iskander A. TAIMANOV
\thanks{Institute of Mathematics, 630090 Novosibirsk, Russia;
e-mail: taimanov@math.nsc.ru}}
\date{}
\maketitle

\abstract{We give a survey on the Weierstrass representations of surfaces 
in three- and four-dimensional spaces, their applications to the theory of the 
Willmore functional and on related problems of spectral theory of the 
two-dimensional Dirac operator with periodic coefficients.}

\tableofcontents

\section{Introduction}

In this paper we survey some results and problems related to
global representations of surfaces in three- and four-spaces by
means of solutions to the Dirac equation and application of these constructions
to a study of the Willmore
functional and its generalizations.

This activity started ten years ago \cite{T1}. 
In this approach
the Gauss map of a surface is represented in terms of solutions
$\psi$ to the equation
$$
\D \psi = 0
$$
where $\D$ is the Dirac operator with potentials
$$
\D = \left(
\begin{array}{cc}
0 & \partial \\
-\bar{\partial} & 0
\end{array}
\right) + \left(
\begin{array}{cc}
U & 0 \\
0 & V
\end{array}
\right).
$$
Such a representation for surfaces in $\R^3$ has different forms
and some prehistory however in this explicit form involving the
Dirac equation first it was written down for inducing surfaces
admitting soliton deformations in \cite{K1} where these
deformations were also introduced.

The appearance of an operator with a well-developed spectral
theory makes it possible to use this theory for study  problems of
global surface theory. Moreover this approach explains an
importance of the Willmore functional since for surfaces in $\R^3$
it is up to a multiple is the squared $L_2$-norm of the potential
$U=V=\bar{U}$ of the operator $\D$ \cite{T1}.

The approach to proving the Willmore conjecture for tori proposed
by us in \cite{T1,T2} and based on the theory of spectral curves
(on one energy level) \cite{DKN} led to a very interesting paper
by Schmidt \cite{Schmidt} where a substantial progress was
achieved however the conjecture stayed unproved.

Therewith the spectral curve of $\D$ with double-periodic
potentials gives rise a notion of the spectral curve of a torus in
$\R^3$ \cite{T2}, in which it is encoded a lot of geometrical information on 
the surface.

Another approach to obtaining lower bounds for the Willmore
functional involved methods of the inverse spectral problem and
algebraic geometry of curves and led to obtaining such estimates
which are quadratic in the dimension of the kernel of $\D$. They
were first obtained for spheres of revolution and some their
generalizations and conjectured for all spheres in \cite{T21} by
using the inverse spectral problem and proved in the full
generality of surfaces of all genera in \cite{FLPP} where the
theory of algebraic curves was applied in a fabulous and unusual
way to surface theory.

Later this representation was generalized for surfaces in $\R^4$
\cite{PP,K2} and three-dimensional Lie groups \cite{T3,Berdinsky}.
In \cite{BFLPP,FLPP} it was proposed to consider the
representations of surfaces in $\R^3$ and $\R^4$ in the conformal
setting from the beginning. However for non-commutative noncompact
three-dimensional groups the analogs of the Willmore functional
appear to be of the form
$$
\int \left(\alpha H^2 + \beta \widehat{K} + \gamma \right)d\mu
$$
where $H$ is the mean curvature, $\widehat{K}$ is the sectional
curvature of the ambient space along the tangent pane to a surface
and $d \mu$ is the induced measure on a surface. 
We note that the functionals of the similar form
$$
\int \left(\alpha H^2 + \beta K + \gamma \right) d\mu
$$
for surfaces in $\R^3$  are well-known in physics as the Helfrich functionals
\cite{Helfrich} (see also, for instance, \cite{Helf1,Helf2}) and
for generic values of $\alpha,\beta,\gamma$ are not conformally
invariant even for surfaces in $\R^3$. For 
surfaces with boundary which are interesting for physical applications
the term containing 
the Gauss curvature $K$ is not reduced to a topological term.

Although until recently these representations were applied mostly
to the problems related to the Willmore functional and its
generalizations we think that they can be effectively used for
study other problems of global surface theory.

\section{Representations of surfaces in three- and four-spaces}
\label{sec3}

\subsection{The generalized Weierstrass formulas for surfaces in $\R^3$}
\label{subsec2.1}

The Grassmannians of oriented
$2$-planes in $\R^n$ are diffeomorphic to quadrics in $\C P^{n-1}$.

Indeed, take a $2$-plane and choose a positively oriented
orthonormal basis $u = (u_1,\dots,u_n), v =(v_1,\dots,v_n)$, i.e.
$|u|=|v|, (u,v)=0$, for the plane. It is defined by a vector $y =
u+iv \in \C^n$ such that
$$
y_1^2 + \dots + y_n^2 = [(u,u) - (v,v)] + 2i(u,v) = 0.
$$
The plane determines such a basis up to rotations of the plane by
an angle $\varphi, 0 \leq \varphi \leq 2\pi$, which result in
transformations $y \to re^{i\varphi}y$. Therefore the Grassmannian
$\widetilde{G}_{n,2}$ of oriented $2$-planes in $\R^n$ is
diffeomorphic to the quadric
$$
y_1^2 + \dots + y_n^2 = 0, \ \ \ (y_1: \ldots :y_n) \in \C P^{n-1},
$$
where $(y_1:\dots:y_n)$ are homogeneous coordinates in $\C P^{n-1}$.
The Grassmannian $G_{n,2}$ of unoriented $2$-planes in $\R^n$ is the
quotient of $\widetilde{G}_{n,2}$ with respect to a fixed-point free
antiholomorphic involution $y \to \bar{y}$.

Given an immersed surface
$$
f: \Sigma \to \R^n
$$
with a (local) conformal parameter $z$, the Gauss map of this
surface is
$$
\Sigma \to \widetilde{G}_{n,2} \ \ : \ \ P \to (x^1_z(P) : \ldots :
x^n_z(P))
$$
where $x^1,\dots,x^n$ are the Euclidean coordinates in $\R^n$ and $P
\in \Sigma$.

There are only two cases when the Grassmannian admits a rational
parameterization:
$$
\widetilde{G}_{3,2} = \C P^1, \ \ \ \widetilde{G}_{4,2} = \C P^1
\times \C P^1
$$
and only in these cases we have the Weierstrass representations of
surfaces.

First  we consider surfaces in $\R^3$.

The Grassmannian  $\widetilde{G}_{3,2}$ is the quadric
$$
y_1^2 + y_2^2 + y_3^2 = 0
$$
admitting the following rational parameterization
\footnote{It is well-known in the number theory as the Lagrange representation of all integer solutions
to the equation $x^2 + y^2 = z^2$.}
$$
y_1 = \frac{i}{2}(b^2 + a^2), \ \ y_2 = \frac{1}{2}(b^2 - a^2), \ \
y_3 = ab, \ \ (a:b) \in \C P^1.
$$
We put
$$
\psi_1 = a, \ \ \ \psi_2 = \bar{b}
$$
and substitute these expressions into the formulas for $x^k_z =
y_k, k=1,2,3$. Since $x^k\in \R$ for all $k$, we have
$$
\Im x^k_{z\bar{z}} = 0, \ \ \ k=1,2,3.
$$
In terms of $\psi$ this condition takes the form of the Dirac
equation
\begin{equation}
\label{dirac3} \D\psi = \left[
\left(\begin{array}{cc} 0 & \partial \\
-\bar{\partial} & 0 \end{array}\right) +
\left(\begin{array}{cc} U & 0 \\
0 & U \end{array}\right) \right] \left(\begin{array}{c} \psi_1 \\
\psi_2 \end{array}\right) = 0, \ \ \ U = \bar{U}.
\end{equation}
Moreover if for a complex-valued function $f$ we have $\Im
f_{\bar{z}} = 0$ then locally we have $f = g_z$ where $g$ is a
real-valued function of the form
$$
g = \int\left[\Re f dx - \Im f dy\right].
$$
We have the following theorem.

\begin{theorem}
1) (\cite{K1}) If $\psi$ meets the Dirac equation (\ref{dirac3})
then the formulas
\begin{equation}
\label{int3} x^k = x^k(0) + \int \left( x^k_z dz + \bar{x}^k_z
d\bar{z}\right), \ \ k=1,2,3,
\end{equation}
with
\begin{equation}
\label{int30}
x^1_z = \frac{i}{2}(\bar{\psi}^2_2 + \psi^2_1), \ \ \
\ x^2_z = \frac{1}{2}(\bar{\psi}^2_2 - \psi^2_1), \ \ \ \ x^3_z =
\psi_1\bar{\psi}_2
\end{equation}
give us a surface in $\R^3$.

2) (\cite{T1})
Every smooth surface in $\R^3$ is locally defined
by the formulas (\ref{int3}) and (\ref{int30}).
\end{theorem}

The proof of the second statement is given above and the proof of
the first statement is as follows: by the Dirac equation, the
integrands in (\ref{int3}) are closed forms and, by the Stokes
theorem, the values of integrals are independent of the choice of
a path in a simply connected domain in $\C$.

This representation of a surface is called {\it a Weierstrass
representation}. In the case $U=0$ it reduces to the classical
Weierstrass (or Weierstrass--Enneper) representation of minimal
surfaces.
\footnote[1]{
$^{^\ast}$ The spinor representation of minimal surfaces
originates in lectures of Sullivan in the late 1980s.
It was successively applied to some problems on
minimal surfaces by Kusner, Schmitt, Bobenko et al.
(see [85] and references therein) and this approach
deserves a complimentary survey.}$^{^\ast}$

The following proposition is derived by straightforward
computations.

\begin{proposition}
\label{prop1}
Given a surface $\Sigma$ defined by the formulas (\ref{int3}) and
(\ref{int30}),

1) $z$ is a conformal parameter on the surface and the induced
metric takes the form
$$
ds^2 = e^{2\alpha}dz d\bar{z}, \ \ \ \ e^\alpha = |\psi_1|^2 +
|\psi_2|^2,
$$

2) the potential $U$ of the Dirac operator equals to
$$
U = \frac{He^\alpha}{2},
$$
where $H$ is the mean curvature,
\footnote {We recall that the
normal vector $N$ meets the condition
$$
\Delta f = 2HN,
$$
where $\Delta = 4 e^{-2\alpha}\partial \bar{\partial}$ is the
Laplace--Beltrami operator corresponding on the surface.} i.e. $H
= \frac{\varkappa_1+\varkappa_2}{2}$ with $\varkappa_1,
\varkappa_2$ the principal curvatures of the surface,

3) the Hopf differential equals $A dz^2 = (f_{zz}, N) dz^2$ and
$$
|A|^2 = \frac{(\varkappa_1 - \varkappa_2)^2 e^{4\alpha}}{16},
$$
$$
A =  \bar{\psi}_2 \partial \psi_1 - \psi_1 \partial \bar{\psi}_2,
$$

4) the Gauss--Weingarten equations take the form
$$
\left[\frac{\partial}{\partial z} - \left(\begin{array}{cc}
 \alpha_z & A e^{-\alpha} \\
-U & 0
\end{array}
\right)\right]\psi = \left[\frac{\partial}{\partial \bar{z}} -
\left(\begin{array}{cc}
0 & U \\
-\bar{A}e^{-\alpha} & \alpha_{\bar{z}}
\end{array}
\right)\right]\psi = 0,
$$

5) the compatibility conditions for the Gauss--Weingarten
equations are the Gauss--Codazzi equations which are
$$
A_{\bar{z}} = (U_z - \alpha_z U)e^{\alpha}, \ \ \
\alpha_{z\bar{z}} + U^2 - A\bar{A}e^{-2\alpha} = 0
$$
and the Gaussian curvature equals $K =
-4e^{-2\alpha}\alpha_{z\bar{z}}$.
\end{proposition}

It is easy to notice that if $\varphi$ meets the Dirac equation
(\ref{dirac3}) then the vector function $\varphi^\ast$ defined by
the formula \beq \label{ast} \varphi = \left(
\begin{array}{c} \varphi_1 \\ \varphi_2
\end{array} \right) \to \varphi^\ast = \left(\begin{array}{c} -\bar{\varphi}_2
\\ \bar{\varphi}_1 \end{array} \right)
\eeq
also meets the Dirac equation.

Let us identify $\R^3$ with the linear space of $2\times 2$
matrices spanned over $\R$ by
$$
e_1 = \left(
\begin{array}{cc}
0 & -i \\
i & 0
\end{array}
\right), \ \ \
e_2 = \left(
\begin{array}{cc}
0 & -1 \\
1 & 0
\end{array}
\right), \ \ \
e_3 = \left(
\begin{array}{cc}
-1 & 0 \\
0 & 1
\end{array}
\right).
$$
We have the orthogonal representation of $SU(2)$ on $\R^3$ which is as follows
$$
e_k \to \rho(S)(e_k) =  \overline{S}^\top e_i S = S^\ast e_k S,  \ \ \ \ k=1,2,3,
$$
$$
S =
\left(
\begin{array}{cc}
\lambda & \mu \\
-\bar{\mu} & \bar{\lambda}
\end{array}
\right) \in SU(2), \ \ \ \mbox{i.e. $|\lambda|^2 + |\mu|^2 =1$},
$$
which descends through $SO(3) = SU(2)/\{\pm 1\}$.
The following lemma is proved by straightforward computations:

\begin{lemma}
\label{lemma-dilation}
If a surface $\Sigma$ is defined by $\psi$ via the Weierstrass
representation, then

1) $\lambda \psi + \mu \psi^\ast$ defines the surface obtained from $\Sigma$
under the transformation $\rho(S)$ of the ambient space $\R^3$.

2) $\lambda \psi$ with $\lambda \in \R$ defines an image of
$\Sigma$ under the dilation $x \to \lambda x$.
\end{lemma}

{\sc Remark.} The formulas (\ref{dirac3}) and (\ref{int30}) were
introduced in \cite{K1} for inducing surfaces which admit soliton
deformations described by the modified Novikov--Veselov equation.
They originate in some complex-valued formulas derived for other
reasons  by Eisenhart \cite{Eisenhart}. Similar representation for
CMC surfaces in terms of the Dirac operator was proposed in 1989
by Abresch (talk in Luminy). It was very soon understood that
these formulas give a local representation of a general surface
(see \cite{T1}; in the proof given above we follow \cite{T3},
later another proof was given in \cite{F1} and from the physical
point of view the representation was described in \cite{Mats}).
Moreover this representation appeared to be equivalent to the
Kenmotsu representation \cite{Kenmotsu} which does not involve
explicitly the Dirac operator.

\subsection{The global Weierstrass representation}
\label{subsec2.2}

In \cite{T1} the global Weierstrass representation was introduced.
For that is necessary to use special $\psi_1$-bundles over
surfaces and to consider the Dirac operator defined on sections of
bundles. In this event

\begin{itemize}
\item
the Willmore functional appears as the integral squared norm of the
potential $U$ and the conformal geometry of a surface is related to
the spectral properties of the corresponding Dirac operator;

\item
since it was proved in \cite{T1} that tori are deformed into tori by
the modified Novikov--Veselov flow and this flow preserves the
Willmore functional, the moduli space of immersed tori is embedded
into the phase space of an integrable system which has the Willmore
functional as an integral of motion.
\end{itemize}

By the uniformization theorem, any closed oriented surface $\Sigma$
is conformally equivalent to a constant curvature surface $\Sigma_0$
and a choice of a conformal parameter $z$ on $\Sigma$ fixes such an
equivalence $\Sigma_0 \to \Sigma$.

Since the quantities
$$
\bar{\psi}_2^2 dz, \ \ \psi_1^2dz, \ \  \psi_1\bar{\psi}_2 dz, \ \
e^{2\alpha}dz d\bar{z}, \ \ H = 2 U e^{-\alpha}
$$
are globally defined on a surface $\Sigma_0$, this leads to the
following description:

\begin{theorem}
[\cite{T1,T21}] Every oriented closed surface $\Sigma$, immersed
in $\R^3$, admits a Weierstrass representation of the form
(\ref{int3})--(\ref{int30}) where $\psi$ is a section of some
bundle $E$ over the surface $\Sigma_0$ which is conformally
equivalent to $\Sigma$ and has constant sectional curvature and
$\D \psi = 0$. Moreover

a) if $\Sigma = \C \cup \{\infty\}$ is a sphere then $\psi$ and
$U$ defined on $\C$ are expanded onto the neighborhood of the
infinity by the formulas:
\begin{equation}
\label{bundle0} (\psi_1,\bar{\psi}_2) \to (z\psi_1,z\bar{\psi}_2),
\ \ U \to |z|^2 U \ \ \ \mbox{as $z \to -z^{-1}$},
\end{equation}
and there is the following asymptotic of $U$:
$$
U = \frac{\const}{|z|^2} + O\left(\frac{1}{|z|^3}\right) \ \ \
\mbox{as $z \to \infty$}.
$$

b) if $\Sigma$ is conformally equivalent a torus $\Sigma_0 =
\R^2/\Lambda$, then
$$
U(z+\gamma,\bar{z}+\overline{\gamma}) = U(z,\bar{z}), \ \
\psi(z+\gamma,\bar{z}+\overline{\gamma}) = \mu(\gamma)
\psi(z,\bar{z}) \ \ \ \mbox{for all $\gamma \in \Lambda$}
$$
where $\mu$ is the character of $\Lambda \to \{\pm 1\}$ which
takes values in $\{\pm 1\}$ and determines the bundle
$$
E \stackrel{\C^2}{\longrightarrow} \Sigma_0
$$
such that $(\psi_1,\psi_2)^\perp$ is a section of $E$.

c) if $\Sigma$ is a surface of genus $g \geq 2$, then $\Sigma_0 =
{\cal H}/\Lambda$, where ${\cal H}$ is the Lobachevsky upper-half
plane and $\Lambda$ is a discrete subgroup of $PSL(2,\R)$ which
acts on ${\cal H} = \{ \Im z > 0\} \subset \C$ as
$$
z \to \gamma(z) = \frac{az + b}{cz + d} , \ \ \ \left(
\begin{array}{cc}
a & b \\ c & d
\end{array}
\right) \in SL(2,\R).
$$
The $\psi$-bundle
$$
E \stackrel{\C^2}{\longrightarrow} \Sigma_0
$$
is defined by the monodromy rules
\begin{equation}
\label{bundle} \gamma:(\psi_1, \bar{\psi}_2) \to
(cz+d)(\psi_1,\bar{\psi}_2).
\end{equation}
and
$$
U(\gamma(z),\overline{\gamma(z)}) = |cz+d|^2 U(z,\bar{z}).
$$

The bundle $E$ splits into the sum of two conjugate bundles $E =
E_0 \oplus \bar{E}_0$ which sections are $\psi_1$ and $\psi_2$
respectively.
\end{theorem}

Since $PSL(2,\R) = SL(2,\R)/\{\pm 1\}$, an element $\gamma \in
PSL(2,\R)$ defines a monodromy up to a sign. The same situation
holds for the torus. Therefore the bundles $E$ are called spin
bundles.

Given a conformal parameter on $\Sigma$, the potential $U$ is fixed
and is called the potential of the representation. Moreover we have
$$
{\cal W}(\Sigma) = 4 \int_\Sigma U^2 dx \wedge dy.
$$
Any section $\psi \in \Gamma(E)$ such that $\D\psi =0$ defines a
surface which is generically not closed but only have a periodic
Gauss map. Therewith the  Weierstrass formulas define an immersion
of the universal covering $\widetilde{\Sigma}$ of $\Sigma$. The
following proposition shows when such an immersion converts into
an immersion a a compact surface.

\begin{proposition}
The Weierstrass representation
defines an immersion of a compact surface $\Sigma$ if and only if
\begin{equation}
\label{period3} \int_{\Sigma_0} \bar{\psi}_1^2 \, d\bar{z} \wedge
\omega = \int_{\Sigma_0} \psi_2^2 \, d\bar{z} \wedge \omega =
\int_{\Sigma_0} \bar{\psi}_1 \psi_2 \, d\bar{z} \wedge \omega = 0
\end{equation}
for any holomorphic differential $\omega$ on $\Sigma_0$.
\end{proposition}

We see that to any immersed torus $\Sigma \subset \R^3$ with a
fixed conformal parameter $z$ it corresponds the Dirac operator
$\D$ with the double-periodic potential
$$
U \ = V  \ = \frac{He^\alpha}{2}
$$
where $H$ is the mean curvature and $e^{2\alpha} dz d\bar{z}$ is the
induced metric.

\subsection{Surfaces in three-dimensional Lie groups}
\label{subsec2.3}

For surfaces in three-dimensional Lie groups the Weierstrass
representation is generalized as follows.

Let $G$ be a three-dimensional Lie group with a left-invariant
metric and let
$$
f: \Sigma \to G
$$
be an immersion of a surface $\Sigma$ into $G$. We denote by
${\cal G}$ the Lie algebra of $G$. Let $z=x+iy$ be a conformal
parameter on the surface.

We take the pullback of $TG$ to a  ${\cal G}$-bundle over
$\Sigma$: ${\cal G} \to E = f^{-1}(TG) \stackrel{\pi}{\to}
\Sigma$, and consider the differential
$$
d_{\cal A}:\Omega^1(\Sigma;E) \to \Omega^2(\Sigma;E),
$$
which acts on $E$-valued $1$-forms:
$$
d_{\cal A} \omega = d'_{\cal A} \omega + d''_{\cal A} \omega
$$
where $\omega = u dz + u^\ast d\bar{z}$ and
$$
d'_{\cal A} \omega = -\nabla_{\bar{\partial}f} u dz \wedge
d\bar{z}, \ \ \ d''_{\cal A} \omega = \nabla_{\partial f}u^\ast dz
\wedge d\bar{z}.
$$
By straightforward computations we obtain the first derivational
equation
\begin{equation}
\label{d1} d_{\cal A} (df) = 0.
\end{equation}
The tension vector $\tau(f)$ is defined via the equation
$$
d_{\cal A} (\ast df) = f\cdot (e^{2\alpha} \tau(f)) dx \wedge dy =
\frac{i}{2} f\cdot(e^{2\alpha} \tau(f)) dz \wedge d\bar{z}
$$
where $f\cdot\tau(f) = 2HN$, $N$ is the normal vector and $H$ is
the mean curvature. This gives the second derivational equation:
\begin{equation}
\label{d2} d_{\cal A} (\ast df) = i e^{2\alpha} H N dz \wedge
d\bar{z}.
\end{equation}

Since the metric is left invariant we rewrite the derivational
equations in terms of
$$
\Psi = f^{-1}\partial f, \ \ \ \Psi^\ast = f^{-1} \bar{\partial} f
$$
as follows:
\begin{equation}
\label{h1}
\partial\Psi^\ast - \bar{\partial}\Psi + \nabla_{\Psi}\Psi^\ast -
\nabla_{\Psi^\ast}\Psi = 0,
\end{equation}
\begin{equation}
\label{h2}
\partial\Psi^\ast + \bar{\partial}\Psi + \nabla_{\Psi}\Psi^\ast +
\nabla_{\Psi^\ast}\Psi = e^{2\alpha} H f^{-1}(N).
\end{equation}
The equation \eqref{h1} is equivalent to \eqref{d1} and the
equation \eqref{h2} is equivalent to \eqref{d2}.

We take an orthonormal basis $e_1,e_2,e_3$ in the Lie algebra
${\cal G}$ of the group $G$ and decompose $\Psi$ and $\Psi^\ast$
in this basis:
$$
\Psi = \sum_{k=1}^3 Z_k e_k, \ \ \ \Psi^\ast = \sum_{k=1}^3
\bar{Z}_k e_k.
$$
Then the equations \eqref{h1} and \eqref{h2} take the form
\begin{equation}
\label{z1} \sum_j (\partial \bar{Z}_j - \bar{\partial} Z_j) e_j +
\sum_{j,k} ( Z_j \bar{Z}_k - \bar{Z}_j Z_k) \nabla_{e_j}e_k = 0,
\end{equation}
\begin{equation}
\label{z2}
\begin{split}
\sum_j (\partial \bar{Z}_j + \bar{\partial} Z_j) e_j +
\sum_{j,k} (Z_j \bar{Z}_k + \bar{Z}_j Z_k) \nabla_{e_j}e_k = \\
2iH \left[ (\bar{Z}_2 Z_3 - Z_2 \bar{Z}_3) e_1 + (\bar{Z}_3 Z_1 -
Z_3 \bar{Z}_1) e_2 + (\bar{Z}_1 Z_2 - Z_1 \bar{Z}_2) e_3 \right].
\end{split}
\end{equation}
Here we assumed that the basis $\{e_1,e_2,e_3\}$ is positively
oriented and therefore
$$
f^{-1}(N) = 
$$
$$
2 i  e^{-2\alpha} \left[ (\bar{Z}_2 Z_3 - Z_2
\bar{Z}_3) e_1 + (\bar{Z}_3 Z_1 - Z_3 \bar{Z}_1) e_2 + (\bar{Z}_1
Z_2 - Z_1 \bar{Z}_2) e_3 \right]
$$
(for $G=SU(2)$ with the Killing metric this formula takes the form
$f^{-1}(N)= 2ie^{-2\alpha} [\Psi^\ast,\Psi]$). Since the parameter
$z$ is conformal we have
$$
\langle \Psi,\Psi \rangle = \langle \Psi^\ast,\Psi^\ast\rangle =
0, \ \ \ \langle \Psi, \Psi^\ast\rangle = \frac{1}{2} e^{2\alpha}
$$
which is rewritten as
$$
Z_1^2 + Z_2^2 + Z_3^2 = 0, \ \ \ |Z_1|^2+|Z_2|^2+|Z_3|^2 =
\frac{1}{2} e^{2\alpha}.
$$
Therefore, as in the case of surfaces in $\R^3$, the vector $Z$ is
parameterized in terms of $\psi$ as follows:
\begin{equation}
\label{spinor} Z_1 = \frac{i}{2} ( \bar{\psi}_2^2 + \psi_1^2), \ \
\ Z_2 = \frac{1}{2} ( \bar{\psi}_2^2 - \psi_1^2), \ \ \ Z_3 =
\psi_1 \bar{\psi}_2.
\end{equation}

Let us show how to reconstruct a surface from $\psi$ meeting the
derivational equations (\ref{h1}) and (\ref{h2}). In the case of
non-commutative Lie groups that can not be done by the integral
Weierstrass formulas.

Let $\psi$ be defined on a surface $\Sigma$ with a complex
parameter $z$ and $\Psi$ constructed from $\psi$ meet (\ref{h1})
and (\ref{h2}). Let us pick up a point $P \in \Sigma$. We substitute
$\psi$ into the formula \eqref{spinor} for the components
$Z_1,Z_2,Z_3$ of $\Psi = \sum_{k=1}^3 Z_k e_k = f^{-1}
\partial f$ and solve the linear equation in the
Lie group $G$:
$$
f_z = f \Psi,
$$
with the initial data $f(P) = g \in G$. Thus we obtain the desired
surface as the mapping
$$
f: \Sigma \to G.
$$
For the group $\R^3$ a solution to such an equation is given by
the Weierstrass formulas (\ref{int3}) and (\ref{int30}).

From the derivation of (\ref{h1}) and (\ref{h2}) it is clear that
any surface $\Sigma$ in $G$ is constructed by this procedure which
is just the generalized Weierstrass representation for surfaces in
Lie groups. In this event we say that $\psi$ generates the surface
$\Sigma$.

Let us write down the derivational equations (\ref{h1}) and
(\ref{h2}) in terms of $\psi$. They are written as the Dirac
equation
\beq
\label{dirac-lie} \D \psi = \left[ \left(
\begin{array}{cc}
0 & \partial \\
-\bar{\partial} & 0
\end{array}
\right)
+
\left(
\begin{array}{cc}
U & 0 \\
0 & V
\end{array}
\right) \right] \psi = 0,
\eeq
the induced metric is given by the same formula
$$
ds^2 = e^{2\alpha} dz d\bar{z}, \ \ \ e^\alpha = |\psi_1|^2 + |\psi_2|^2,
$$
and the Hopf quadratic differential $Adz^2$ equals to
$$
A =  (\bar{\psi}_2 \partial \psi_1 - \psi_1 \partial \bar{\psi}_2)
+ \left(\sum_{j,k} Z_j Z_k \nabla_{e_j} e_k, N \right).
$$
For a compact Lie group with the Killing metric, in particular for
$G = SU(2)$, we have $\nabla_{e_j}e_k = - \nabla_{e_k}e_j$ and the
Hopf differential takes the same form as for surfaces in $\R^3$:
$A = \bar{\psi}_2 \partial \psi_1 - \psi_1 \partial \bar{\psi}_2$.

We consider three-dimensional Lie groups with Thurston's
geometries. Let us recall that by Thurston's theorem
\cite{Scott,Thurston} all three-dimensional maximal simply
connected geometries $(X,\mathrm{Isom}\, X)$ admitting compact
quotients are given by the following list:

1) the geometries with constant sectional curvature: $X = \R^3,
S^3$, and $H^3$;

2) two product geometries: $X = S^2 \times \R$ and $H^2 \times
\R$;

3) three geometries modelled on Lie groups $\nil, \sol$, and
$\sll$ with certain left invariant metrics.

The group $\R^3$ with the Euclidean metric was already considered
above. Hence we are left with four groups:
$$
SU(2) = S^3, \ \ \ \nil, \ \ \ \sol, \ \ \ \sll
$$
where $\nil$ is a nilpotent group, $\sol$ is a solvable group, and
$\sll$ is the universal cover of the group $SL_2(\R)$:
$$
\nil = \left\{\left(
\begin{array}{ccc}
1 & x & z \\
0 & 1 & y \\
0 & 0 & 1
\end{array}
\right)\right\}, \ \ \ \sol = \left\{\left(
\begin{array}{ccc}
e^{-z} & 0 & x \\
0 & e^z & y \\
0 & 0 & 1
\end{array}
\right)\right\},
$$
with $x,y,z \in \R$.

The case $G = SU(2)$ was studied in \cite{T3} and surfaces in the other groups were
considered in \cite{Berdinsky}:

\bi
\item
$G = SU(2)$:
$$
U = \bar{V} = \frac{1}{2}(H - i)(|\psi_1|^2 + |\psi_2|^2),
$$
the Gauss--Weingarten equations are
$$
\left[\frac{\partial}{\partial z} - \left(\begin{array}{cc}
 \alpha_z & A e^{-\alpha} \\
-U & 0
\end{array}
\right)\right]\psi = \left[\frac{\partial}{\partial \bar{z}} -
\left(\begin{array}{cc}
0 & \bar{U} \\
-\bar{A}e^{-\alpha} & \alpha_{\bar{z}}
\end{array}
\right)\right]\psi = 0,
$$
their compatibility conditions --- the Gauss--Codazzi equations
---
take the form
$$
\alpha_{z\bar{z}} + |U|^2 - |A|^2e^{-2\alpha} = 0, \ \ \
A_{\bar{z}} = (\bar{U}_z - \alpha_z\bar{U})e^{\alpha}
$$
with $A dz^2$ the Hopf differential:
$$
A = \bar{\psi}_2 \partial \psi_1 - \psi_1 \partial
\bar{\psi}_2,
$$

\item
$G = \nil$:
$$
U = V = \frac{H}{2}(|\psi_1|^2+|\psi_2|^2) +
\frac{i}{4}(|\psi_2|^2-|\psi_1|^2),
$$
the Gauss--Weingarten equations are
$$
\left[
\frac{\partial}{\partial
z} - \left(\begin{array}{cc}
 \alpha_z - \frac{i}{2}\psi_1 \bar{\psi}_2 & A e^{-\alpha} \\
-U & 0
\end{array}
\right)\right]\psi = 0, $$
$$
\left[\frac{\partial}{\partial
\bar{z}} - \left(\begin{array}{cc}
0 & U \\
-\bar{A}e^{-\alpha} & \alpha_{\bar{z}} - \frac{i}{2}\bar{\psi}_1\psi_2
\end{array}
\right)\right]\psi = 0,
$$
and the Gauss--Codazzi equations have the following shape
$$
\alpha_{z\bar{z}} - |A|^2 e^{-2\alpha} + \frac{H^2}{4} e^{2\alpha}
= \frac{1}{16}(3|\psi_1|^4 + 3 |\psi_2|^4 - 10 |\psi_1|^2
|\psi_2|^2),
$$
$$
A_{\bar{z}} - \frac{H_z}{2}e^{2\alpha} +
\frac{1}{2}(|\psi_2|^4-|\psi_1|^4)\psi_1 \bar{\psi}_2 = 0
$$
where the Hopf differential equals to
$$
A =  (\bar{\psi}_2 \partial \psi_1 - \psi_1 \partial \bar{\psi}_2)
+ i \psi_1^2 \bar{\psi}_2^2,
$$

\item
$G = \sll$:
$$
U = \frac{H}{2}(|\psi_1|^2+|\psi_2|^2) +
i\left(\frac{1}{2}|\psi_1|^2-\frac{3}{4}|\psi_2|^2\right),
$$
$$
V = \frac{H}{2}(|\psi_1|^2+|\psi_2|^2) +
i\left(\frac{3}{4}|\psi_1|^2-\frac{1}{2}|\psi_2|^2\right),
$$
the Gauss--Weingarten equations are written as
$$
\left[
\frac{\partial}{\partial
z} - \left(\begin{array}{cc}
 \alpha_z + \frac{5i}{4}\psi_1 \bar{\psi}_2 & A e^{-\alpha} \\
-U & 0
\end{array}
\right)\right]\psi = 0, $$
$$
\left[\frac{\partial}{\partial
\bar{z}} - \left(\begin{array}{cc}
0 & V \\
-\bar{A}e^{-\alpha} & \alpha_{\bar{z}} + \frac{5i}{4}\bar{\psi}_1\psi_2
\end{array}
\right)\right]\psi = 0,
$$
and the Gauss--Codazzi equations take the form
$$
\alpha_{z \bar{z}} -  e^{-2 \alpha} |A|^{2} + \frac{1}{4}
e^{2\alpha} H^{2} =
e^{2\alpha} - 5|Z_3|^2,
$$
$$
\bar{\partial}\left (A + \frac{5 Z_3^2}{2(H-i)} \right) =
\frac{1}{2} H_z e^{2\alpha} + \bar{\partial}
\left(\frac{5}{2(H-i)} \right) {Z_3}^2,
$$
where
$$
A = (\bar{\psi}_2 \partial \psi_1 - \psi_1 \partial \bar{\psi}_2)
-\frac{5i}{2}  \psi_1^2 \bar{\psi}_2^2,
$$

\item
$G = \sol$: we consider only domains where $Z_3 = \psi_1
\bar{\psi}_2$ in which
$$
U = \frac{H}{2}(|\psi_1|^2+|\psi_2|^2) + \frac{1}{2}\bar{\psi}_2^2
\frac{\bar{\psi}_1}{\psi_1},
$$
$$
V = \frac{H}{2}(|\psi_1|^2+|\psi_2|^2) + \frac{1}{2}\bar{\psi}_1^2
\frac{\bar{\psi}_2}{\psi_2},
$$
the Gauss--Weingarten equations consist in the Dirac equation and
the following system
$$
\partial \psi_1 = \alpha_z \psi_1 +
Ae^{-\alpha} \psi_2 -
\frac{1}{2} \bar{\psi}_2^3, \ \
\bar{\partial} \psi_2 = -\bar{A} e^{-\alpha} \psi_1  +
\alpha_{\bar{z}}\psi_2- \frac{1}{2} \bar{\psi}_1^3,
$$
the Gauss--Codazzi equations are
$$
\alpha_{z \bar{z}} -  e^{-2 \alpha} |A|^{2} + \frac{1}{4}
e^{2\alpha} H^{2} = \frac{1}{4} (6 |\psi_1|^2 |\psi_2|^2
-(|\psi_1|^4 +  |\psi_2|^4)),
$$
$$
A_{\bar{z}} - \frac{1}{2} H_z
e^{2\alpha} = (|\psi_2|^4 - |\psi_1|^4) \psi_1 \bar{\psi}_2
$$
with
$$
A =  (\bar{\psi}_2 \partial \psi_1 - \psi_1 \partial \bar{\psi}_2) +
\frac{1}{2} (\bar{\psi}_2^4 - \psi_1 ^4).
$$
\ei

It needs to make several explaining remarks:

1) for the last three groups the term $Z_3$ appears in the
formulas. The direction of the vector $e_3$ has different sense
for these groups:

1a) the groups $\nil$ and $\sll$ admit $S^1$-symmetry which is the
rotation around the geodesic drawn in the direction of $e_3$. This
rotation together with left shifts generate $\mathrm{Isom} G$;

1b) for $\sol$ the vectors $e_1$ and $e_2$ commute. Therefore
the equation $Z_3 = \psi_1 \bar{\psi}_2
=0$ can be valid in an open subset $B$ of a surface and
therewith the
Dirac equation is not extended by continuity onto the whole
surface.
Since $H=0$ in $B$ we put
$$
U = V = 0  \ \ \ \ \mbox{for $\psi_1 \bar{\psi}_2 = 0$ and $G =
\sol$}.
$$
However on the boundary $\partial B$ of the set $\{Z_3 \neq 0\}$
the potentials $U_\sol$ and $V_\sol$ are not always correctly
defined due to the indeterminacy of $\frac{\bar{\psi}_1}{\psi_1}$
for $\psi_1=0$ and the Dirac equation with given potentials is valid
outside $\partial B$;

2) for $G = \R^3$ or $SU(2)$
the Gauss--Codazzi equations are derived as follows.
We have
$$
R \psi = (\partial - {\cal A})(\bar{\partial} - {\cal B}) \psi -
(\bar{\partial} - {\cal B})(\partial - {\cal A}) \psi =
({\cal A}_{\bar{z}} - {\cal B}_z + [{\cal A},{\cal B}]) \psi = 0,
$$
where $(\partial - {\cal A})\psi = (\bar{\partial} - {\cal B})\psi
= 0$ are the Gauss--Weingarten equations and the vector function
$\psi^\ast$ (see (\ref{ast})) meets the same equation $R\psi^\ast =
0$ which together with $R\psi=0$ implies that $R = {\cal
A}_{\bar{z}} - {\cal B}_z + [{\cal A},{\cal B}] = 0$. For other
groups the equations $\D \psi^\ast=0$ and $R\psi^\ast=0$ does not
hold and, in particular, the kernel of the Dirac operator can not
be treated as a vector space over quaternions. Therefore the
Gauss--Codazzi equations are derived in \cite{Berdinsky} by other
methods;

3) in fact the Dirac equations in the case of non-commutative
groups are nonlinear in $\psi$ due to the constraints on the
potentials. Therefore if $\psi$ defines a surface then $\lambda
\psi$ does not define another surface for $|\lambda|\neq 1$ since
these groups do not admit dilations. For $SU(2)$ the mapping
(\ref{ast}) maps a solution of the Dirac equation to another
solution of it and the analog of part 1 of Lemma
\ref{lemma-dilation} holds: $\lambda \psi + \mu \psi^\ast,
|\lambda|^2 + |\mu|^2 =1$, defines an image of the initial surface
under some inner automorphism of $SU(2)$ corresponding to a
rotation of the Lie algebra.

Let us expose some corollaries. Since the case $G=SU(2)$ was
well-studied,~\footnote{For $SU(2)$ the minimal surface equations
are $\bar{\partial} \psi_1 =
-\frac{i}{2}(|\psi_1|^2+|\psi_2|^2)\psi_2$, $\partial \psi_2 =
\frac{i}{2}(|\psi_1|^2+|\psi_2|^2)\psi_1$, and CMC surfaces are
distinguished by the condition $A_{\bar{z}} = 0$.} we consider
only other groups.

\bet 1) Given $\psi$ generating a minimal surface in a Lie group,
the following equations hold:
$$
\bar{\partial} \psi_1 = \frac{i}{4}(|\psi_2|^2-|\psi_1|^2)\psi_2, \ \ \
\partial \psi_2 = - \frac{i}{4}(|\psi_2|^2-|\psi_1|^2) \psi_1 \ \
\ \ \mbox{for $G=\nil$},
$$
$$
\bar{\partial} \psi_1 =
i\left(\frac{3}{4}|\psi_1|^2-\frac{1}{2}|\psi_2|^2\right) \psi_2,
\
\partial \psi_2 =
-i\left(\frac{1}{2}|\psi_1|^2-\frac{3}{4}|\psi_2|^2\right)\psi_1 \
$$
$$
\ \mbox{for $G = \sll$},
$$
$$
\bar{\partial}\psi_1 = \frac{1}{2} \bar{\psi}_1^2 \bar{\psi}_2, \
\ \
\partial \psi_2 = - \frac{1}{2} \bar{\psi}_1 \bar{\psi}_2^2 \ \ \
\ \mbox{for $G = \sol$}.
$$

2) (Abresch \cite{Abresch}) If a surface has constant mean
curvature then the following quadratic differential
$\widetilde{A}dz^2$ is holomorphic:
$$
\widetilde{A} dz^2 = \left (A + \frac{{Z_3}^2}{2H+i} \right)dz^2 \
\ \ \ \mbox{for $G = \nil$},
$$
$$
\widetilde{A} dz^2 = \left(A + \frac{5}{2(H-i)}Z^2_3\right) dz^2 \
\ \ \ \mbox{for $G = \sll$}.
$$

3) If for a surface in $G = \nil$ the differential
$\widetilde{A}dz^2$ is holomorphic then the surface has constant
mean curvature.
\eet

It would be interesting to understand relations of formulas for
constant mean curvature surfaces in these groups to soliton
equations. Such relations are well-known for such surfaces in
$\R^3$ and $SU(2)$.

The analogs of the statement 2) are known also for surfaces in the
product geometries $S^2 \times \R$ and $H^2 \times \R$
\cite{AbreschR}. However only for surfaces in $\nil$ the converse
--- the statement 3) --- is also proved.
\footnote{After this paper was submitted for a publication it was shown that
for surfaces in $G = \sll$ and in $H^2 \times \R$ 
the analogous statement does not hold: there are
surfaces for which the quadratic differential $\widetilde{A}dz^2$
is holomorphic and which are not CMC surfaces \cite{Fernandez}.}

We remark that for minimal surfaces in $\nil$ and $\sol$ the
analog of the Weierstrass representation was obtained by other
methods in \cite{I1,I2}. Other approaches to study surfaces in Lie
groups were used in \cite{Daniel,FG}.

\subsection{The quaternion language and quaternionic function theory}
\label{subsec2.4}

Pedit and Pinkall wrote the Weierstrass representation
representation of surfaces in $\R^3$ in the quaternion language
and then extended it to surfaces in $\R^4$ \cite{PP} (see some
preliminary results in \cite{Kamberov,KPP,Richter}).

Indeed, the idea of using quaternions comes from the symmetry of
the kernel of the Dirac operator under the transformation
(\ref{ast}) (remark that that holds for surfaces in $\R^3$ and
$SU(2)$ when $U = \bar{V}$ and is not valid for surfaces in other
three-dimensional Lie groups).

We identify $\C^2$ with the space of quaternions $\H$ as follows
$$
(z_1,z_2) \to z_1 + \j z_2 = \left(
\begin{array}{cc}
z_1 & -\bar{z}_2 \\
z_2 & \bar{z}_1
\end{array}
\right)
$$
and consider two matrix operators
$$
\bar{\partial} = \left(
\begin{array}{cc}
\bar{\partial} & 0 \\
0 & \partial
\end{array}
\right), \ \ \ \j U = \j \left(\begin{array}{cc} U & 0 \\ 0 &
\bar{U} \end{array}\right) = \left(\begin{array}{cc} 0 & -\bar{U} \\
U & 0 \end{array}\right).
$$
Here $\j$ is one the standard generators of quaternions and we
have
$$
\j^2 = -1, \ \ \ z\j = j \bar{z}, \ \ \ \bar{\partial} \j = \j
\partial.
$$
Then the Dirac equation takes the form
$$
(\bar{\partial} + \j U)(\psi_1 + \j \psi_2) =
(\bar{\partial}\psi_1 - \bar{U}\psi_2) + \j (\partial \psi_2 + U
\psi_1) = 0.
$$
Since, by (\ref{bundle}), $\psi_1$ and $\bar{\psi}_2$ are sections
of the same bundle $E$ it is worth working in terms of quaternions
to rewrite the Dirac equation as
$$
(\bar{\partial} + \j U)(\psi_1 + \bar{\psi}_2 \j) = 0.
$$

One may treat $L = E_0 \oplus E_0$ as a quaternionic line bundle
whose sections take the form $\psi_1 + \bar{\psi}_2 \j$ and which
is endowed by some quaternion linear endomorphism $J$ such that
$J^2=-1$. In our case $J$ simply acts as the right-side
multiplication by $\j$:
$$
J:(\psi_1,\bar{\psi}_2) \to (-\bar{\psi}_2,\bar{\psi}_1) \ \
\mbox{or} \ \ \ \psi_1 + \bar{\psi}_2 \j \to (\psi_1 +
\bar{\psi}_2 \j)\j = -\bar{\psi_2} + \psi_1 \j.
$$
This mapping $J$ defines for any quaternion fiber a canonical
splitting into $\C \oplus \C$ (in our case this is a splitting
into $\psi_1$ and $\bar{\psi}_2$). In \cite{PP,BFLPP}
such a bundle is called a ``complex quaternionic line bundle''.

The Dirac operator in these terms is just
$$
\D \psi = (\bar{\partial} + \j U) (\psi_1 + \bar{\psi}_2 \j) =
(\bar{\partial} \psi_1 -\bar{U}\psi_2) + (\bar{\partial}
\bar{\psi}_2 + \bar{U}\bar{\psi}_1)\j
$$
and we see that its kernel is invariant under the right-side
multiplications by constant quaternions (see Lemma \ref{lemma-dilation})
and thus the kernel can be considered as a linear space over $\H$.
By (\ref{bundle}), we have an operator
$$
\D: \Gamma(L) \to \Gamma(\bar{K}L)
$$
where, given a bundle $V$, we denote by $\Gamma(V)$ the space of
sections of a bundle $V$ and $\bar{K}$ is the bundle of $1$-forms
of type $(0,1)$, i.e., of type $fd\bar{z}$, over a surface
$\Sigma_0$.

This operator, of course, is not linear with respect to right-side
multiplications on quaternion-valued functions and the following
evident formula holds:
$$
\D (\psi \lambda) = (\D \psi)\lambda + \psi_1 (\bar{\mu} +
\j \partial \eta) + \bar{\psi}_2
(-\bar{\partial}\eta + \j \partial \mu),
$$
where $\lambda = \mu + \j \eta = \mu +
\bar{\eta} \j$. In \cite{PP} this formula is written in a
coordinate-free form as
$$
\D (\psi \lambda) = (\D \psi) \lambda + \frac{1}{2} (\psi d\lambda
+ J \psi \ast d\lambda),
$$
the potential $U$ multiplied by $\j$ from the left is called the
Hopf field $Q = \j U$ of the connection $\D$ on $L$ and the
quantity
$$
{\cal W} = \int_{\Sigma_0} |U|^2 dx \wedge dy
$$
is called the Willmore energy of the connection $\D$.

Although at the beginning this quaternion language had looked very
artificial, at least to us,
it led to an extension of the Weierstrass
representation for surfaces in $\R^4$ \cite{PP}.
Later it was developed into a tool of investigation
based on working out analogies between complex algebraic
geometry and the theory of complex quaternionic line bundles.
It appears that that is effectively applied to a study of special types of
surfaces and B\"acklund transforms between them in the conformal
setting, i.e. not distinguishing between $\R^4$ and $S^4$
\cite{BFLPP,LPP}.
Finally this approach had led to a fabulous
extension of the Pl\"ucker type relations from complex algebraic
geometry onto geometry of complex quaternionic line 
bundles and application of that to obtaining lower bound for the
Willmore functional \cite{FLPP} (see in \S \ref{subsec5.4}). Therewith 
this theory deals in the same
manner with general bundles $L$ not always coming from the surface
theory \cite{FLPP}. The bundles related to surfaces are
distinguished by their degrees: it follows from (\ref{bundle0})
and (\ref{bundle}) that
$$
\deg E_0 = \mathrm{genus}\, (\Sigma_0) - 1 = g-1.
$$

\subsection{Surfaces in $\R^4$}
\label{subsec2.5}

The Grassmannian of oriented two-planes in $\R^4$ is diffeomorphic
to the quadric
$$
y_1^2 + y_2^2 + y_3^2 + y_4^2 = 0, \ \ \ y \in\C P^3.
$$
We take another coordinates in
$y^\prime_1,y^\prime_2,y^\prime_3,y^\prime_4$ in $\C^4$:
$$
y_1 = \frac{i}{2}(y^\prime_1 + y^\prime_2), \ \ \ y_2 =
\frac{1}{2}(y^\prime_1-y^\prime_2),
\ \ \
y_3 = \frac{1}{2}(y^\prime_3 + y^\prime_4), \ \ \ y_4 =
\frac{i}{2}(y^\prime_3-y^\prime_4).
$$
In terms of these coordinates $\widetilde{G}_{4,2}$ is defined by
the equation
$$
y^\prime_1 y^\prime_2 = y^\prime_3 y^\prime_4.
$$
It is clear that there is a diffeomorphism
$$
\C P^1 \times \C P^1 \to \widetilde{G}_{4,2}
$$
given by the Segre mapping
$$
y^\prime_1 = a_2 b_2, \ \ y^\prime_2 = a_1 b_1, \ \ y^\prime_3 = a_2
b_1, \ \ y^\prime_4 = a_1 b_2
$$
where $(a_1:a_2)$ and $(b_1:b_2)$ are the homogeneous coordinates
on different copies of $\C P^1$.

Let us parameterize $x^k_z, k=1,2,3,4$, in terms of these
homogeneous coordinates and put
$$
a_1 = \varphi_1, \ \ a_2 =\bar{\varphi}_2, \ \  b_1 = \psi_1, \ \
b_2 = \bar{\psi}_2.
$$

In difference with the $3$-dimensional situation this
parameterization is not unique even up to the multiplication by
$\pm 1$ and the vector functions $\psi$ and $\varphi$ are defined
up to the gauge transformations
\begin{equation}
\label{gauge1} \left(
\begin{array}{c}
\psi_1 \\ \psi_2
\end{array}
\right) \to \left(
\begin{array}{c}
e^f\psi_1 \\ e^{\bar{f}}\psi_2
\end{array}
\right), \ \ \ \left(
\begin{array}{c}
\varphi_1 \\ \varphi_2
\end{array}
\right) \to \left(
\begin{array}{c}
e^{-f} \varphi_1 \\ e^{-\bar{f}}\varphi_2
\end{array}
\right),
\end{equation}
where $f$ is an arbitrary function. However the mappings
$$
G_\psi = (\psi_1 : \bar{\psi}_2), \ \ \ G_\varphi = (\varphi_1 :
\bar{\varphi}_2)
$$
into $\C P^1$ are correctly defined and split the Gauss map
$$
G  = (G_\psi, G_\varphi): \Sigma \to \widetilde{G}_{4,2} = \C P^1
\times \C P^1.
$$

We have the following formulas for an immersion of the surface:
\begin{equation}
\label{int4} x^k = x^k(0) + \int \left( x^k_z dz + \bar{x}^k_z
d\bar{z}\right), \ \ k=1,2,3,4,
\end{equation}
with
\begin{equation}
\label{int40}
\begin{split}
x^1_z = \frac{i}{2} (\bar{\varphi}_2\bar{\psi}_2 + \varphi_1
\psi_1), \ \ \ \ x^2_z = \frac{1}{2} (\bar{\varphi}_2\bar{\psi}_2 -
\varphi_1 \psi_1),
\\
x^3_z = \frac{1}{2} (\bar{\varphi}_2 \psi_1 + \varphi_1
\bar{\psi}_2), \ \ \ \ x^4_z = \frac{i}{2} (\bar{\varphi}_2 \psi_1 -
\varphi_1 \bar{\psi}_2).
\end{split}
\end{equation}
Of course, as in the three-dimensional case, these formulas define a
surface if and only if the integrands are closed forms or,
equivalently,
$$
\Im x^k_{z\bar{z}} = 0, \ \ \ k=1,2,3,4.
$$
This is rewritten as
\begin{equation}
\label{codazzi4} \left(\bar{\varphi}_2 \psi_1 \right)_{\bar{z}} =
\left(\bar{\varphi}_1 \psi_2 \right)_z, \ \ \ \left(\bar{\varphi}_2
\bar{\psi}_2\right)_{\bar{z}} = - \left(\bar{\varphi}_1
\bar{\psi}_1\right)_z.
\end{equation}
For generic $\varphi$ and $\psi$ these conditions are not written in
terms of Dirac equations.

However there is the following

\begin{theorem}
[\cite{T5}]
\label{theorem4} Let $r:W \to \R^4$ be an immersed surface with a
conformal parameter $z$ and let $G_\psi = (e^{i\theta}\cos \eta:
\sin \eta)$ be one of the components of its Gauss map.

There exists another representative $\psi$ of this mapping $G_\psi =
(\psi_1:\bar{\psi}_2)$ such that it meets the Dirac equation
\begin{equation}
\label{dirac4} \D \psi = 0, \ \ \ \ \ \D = \left(
\begin{array}{cc}
0 & \partial \\
-\bar{\partial} & 0
\end{array}
\right)+ \left(
\begin{array}{cc}
U & 0 \\
0 & \bar{U}
\end{array}
\right)
\end{equation}
with some potential $U$.

A vector function $\psi = (e^{g+i\theta}\cos \eta,e^{\bar{g}}\sin
\eta)$ is defined from the equation
\begin{equation}
\label{maineq}
g_{\bar{z}} = -i\theta_{\bar{z}} \cos^2 \eta,
\end{equation}
whose solution is defined up to addition of an arbitrary
holomorphic function $h$ and the corresponding potential $U$ is
defined by the formula
$$
U =  -e^{\bar{g}-g-i\theta}(i\theta_z \sin \eta \cos \eta +\eta_z)
$$
up to multiplication by $e^{\bar{h}-h}$.

Given the function $\psi$, a function $\varphi$ which represents
another component $G_\varphi$ of the Gauss map meets the equation
\begin{equation}
\label{dirac40}
\D^\vee \varphi = 0, \ \ \ \
\D^\vee = \left(\begin{array}{cc} 0 & \partial \\
-\bar{\partial} & 0 \end{array}\right) +
\left(\begin{array}{cc} \bar{U} & 0 \\
0 & U \end{array}\right).
\end{equation}

Different lifts into $\C^2 \times \C^2$ of the Gauss mapping $G:
\Sigma \to \C P^1 \times \C P^1$ are related by gauge
transformations
\begin{equation}
\label{gauge2}
\left(
\begin{array}{c}
\psi_1 \\ \psi_2
\end{array}
\right) \to \left(
\begin{array}{c}
e^h\psi_1 \\ e^{\bar{h}}\psi_2
\end{array}
\right), \ \ \ \left(
\begin{array}{c}
\varphi_1 \\ \varphi_2
\end{array}
\right) \to \left(
\begin{array}{c}
e^{-h} \varphi_1 \\ e^{-\bar{h}}\varphi_2
\end{array}
\right),
\ \ \
U \to \exp{(\bar{h}-h)}U, \
\end{equation}
where $h$ is an arbitrary holomorphic function on $W$.
\end{theorem}

\begin{corollary}
Any oriented surface in $\R^4$ is defined by the formulas
(\ref{int4}) and (\ref{int40}) where the vector functions $\psi$
and $\varphi$ meet the equations of the Dirac type (\ref{dirac4})
and (\ref{dirac40}):
$$
\D \psi = \D^\vee \varphi = 0.
$$
The induced metric equals
$$
e^{2\alpha} dzd\bar{z} =
(|\psi_1|^2+|\psi_2|^2)(|\varphi_1|^2+|\varphi_2|^2)dz d\bar{z}
$$
and the norm of the mean curvature vector ${\bf H} = \frac{2
x_{z\bar{z}}}{e^{2\alpha}}$ meets the equality
$$
|U| = \frac{|{\bf H}| e^\alpha}{2}.
$$
\end{corollary}

Let us consider the diagonal embedding
$$
\widetilde{G}_{3,2} = \C P^1 \to \widetilde{G}_{4,2} = \C P^1 \times
\C P^1.
$$
If $\varphi$ and $\psi$ generate a surface and lie in the
diagonal: $\varphi = \pm \psi$, then $x^4 = 0$ and we obtain a
Weierstrass representation of the surface in $\R^3$.

The formulas (\ref{int4}) and (\ref{int40}) appeared for inducing
surfaces in \cite{K2}. This corollary demonstrates that they are
general although this has to follow also from \cite{PP} where such a
representation was first indicated in the quaternion language.

We indicate two specific features of the representation of
surfaces in $\R^4$ which were not discussed in the previous
papers:

\bi
\item
given a surface, a representation is not unique and different
representations are related by nontrivial gauge transformations;

\item
a Weierstrass representation of some domain is not always expanded
onto the whole surface and in difference with the
three-dimensional case it needs to solve $\bar{\partial}$-problem
(\ref{maineq}) on the whole surface to obtain a representation of
the surface.
\ei

Indeed, take $\psi$ and $\varphi$ generating surface $\Sigma$, a
domain $W \subset \Sigma$ and a holomorphic function $f$ on $W$
which is not analytically extended outside $W$. Then by
(\ref{gauge1}) we construct from $\psi,varphi$, and $f$ a another
representation of $W$ which is not expanded outside $W$.

{\sc Example. Lagrangian surfaces in $\R^4$.}

We expose the Weierstrass representation of Lagrangian surfaces in
$\R^4$ obtained by Helein and Romon \cite{HR}. A reduction of the
formulas (\ref{int40}) to the formulas from \cite{HR} was
demonstrated by Helein \cite{Helein}).

Let us take the following symplectic form on $\R^4$:
$$
\omega = dx^1 \wedge dx^2 + dx^3 \wedge dx^4.
$$
We recall that an $n$-dimensional submanifold $\Sigma$ of a
$2n$-dimensional symplectic manifold $M^{2n}$ with a symplectic form
$\omega$ is called Lagrangian if the restriction of $\omega$ onto
$\Sigma$ vanishes:
$$
\omega\vert_{\Sigma} = 0.
$$
This means that at any point $x \in \Sigma$ the restriction of
$\omega$ onto the tangent space $T_x \Sigma$ vanishes, i.e., $T_x
\Sigma$ is a Lagrangian $n$-plane in $\R^{2n}$.

The condition that a $2$-plane is Lagrangian in $\R^4$ is written
as
$$
\Im \left(y_1 \bar{y}_2 + y_3 \bar{y}_4\right) = 0
$$
or
$$
|y^\prime_1|^2 - |y^\prime_2|^2 - |y^\prime_3|^2 + |y^\prime_4|^2 =
0.
$$
In terms of $a_1,a_2,b_1$, and $b_2$ it takes the form
$$
|b_1|^2 = |b_2|^2.
$$
Hence the Grassmannian of Lagrangian $2$-planes in $\R^4$ is the
product of manifolds
$$
G^{\mathrm{Lag}}_{4,2} = \C P^1 \times S^1
$$
where $\C P^1$ is parameterized by $(a_1:a_2)$ and $S^1$ is
parameterized by
$$
\beta = \frac{1}{i}\log \frac{b_1}{b_2} \mod 2\pi.
$$
This quantity $\beta$ is called the Lagrangian angle. We conclude
that a surface is Lagrangian if and only if
$$
|\psi_1| = |\psi_2|
$$
in its Weierstrass representation. Let us put
$$
s = \left(\frac{e^{i\beta}}{\sqrt{2}},\frac{1}{\sqrt{2}}\right), \ \
\ (s_1 : s_2) = (\psi_1:\bar{\psi}_2) \in \C P^1
$$
and apply Theorem \ref{theorem4}. We obtain the following formulas:
$$
g = -\frac{i\beta}{2}, \ \ U = -\frac{1}{2}\beta_z, \ \ \psi_1 =
\psi_2 = \frac{1}{\sqrt{2}} e^{i\beta/2}.
$$
For any solution $\varphi$ to the equation $\D^\vee \varphi = 0$,
we obtain a Lagrangian surface defined by $\psi$ and $\varphi$ via
(\ref{int40}. Moreover all Lagrangian surfaces are represented in
this form.

Let
$$
f: \Sigma \to \R^4
$$
be an immersion of an oriented closed surface in $\R^4$. By
Theorem \ref{theorem4}, this surface is locally defined by the
formulas (\ref{int4}) and (\ref{int40}). A globalization is
similar to the case of surfaces in $\R^3$ and on the quaternion
language was described in \cite{PP,BFLPP} however to obtain it one
has to solve a $\bar{\partial}$-problem of the surface \cite{T5}:

\begin{proposition}
Given a Weierstrass representation of an immersion of an oriented
closed surface $\Sigma$ into $\R^4$, the corresponding functions
$\psi$ and $\varphi$ are sections of the $\C^2$-bundles $E$ and
$E^\vee$ over $\Sigma$ which are as follows:

1) $E$ and $E^\vee$ split into sums of pair-wise conjugate line
bundles
$$
E = E_0 \oplus \bar{E}_0, \ \ \ \ E^\vee = E^\vee_0 \oplus
\bar{E}^\vee_0
$$
such that $\psi_1$ and $\bar{\psi}_2$ are sections of $E_0$ and
$\varphi_1$ and $\bar{\varphi}_2$ are sections of $E^\vee_0$;

2) the pairing of sections of $E_0$ and $E^\vee_0$ is a $(1,0)$ form
on $\Sigma$: if
$$
\alpha \in \Gamma(E_0), \ \ \ \beta \in \Gamma(E^\vee_0),
$$
then
$$
\alpha \beta dz
$$
is a correctly defined $1$-form on $\Sigma$;

3) the Dirac equation $\D \psi =0$ implies that $U$ is a section of
the same line bundle $E_U$ as
$$
\frac{\partial \gamma}{\alpha} \in \Gamma(E_U) \ \ \ \mbox{for
$\alpha \in \Gamma(E_0),\ \  \gamma \in \Gamma(\bar{E}_0)$}
$$
and $U\bar{U} dz \wedge d\bar{z}$ is a correctly defined
$(1,1)$-form on $\Sigma$ whose integral over the surface equals
$$
\int_{\Sigma} U\bar{U} dz \wedge d\bar{z} = - \frac{i}{2} {\cal
W}(\Sigma)
$$
where ${\cal W}(\Sigma) = \int_{\Sigma} |{\bf H}|^2 d\mu$ is the
Willmore functional.
\end{proposition}

The gauge transformation (\ref{gauge2}) show that in difference
with the three-dimen\-sio\-nal case $\psi$ are not necessarily
sections of spin bundles.

For tori we derive from Theorem \ref{theorem4} the following result.

\begin{theorem}
[\cite{T5}] \label{torus} Let $\Sigma$ be a torus in $\R^4$ which
is conformally equivalent to $\C/\Lambda$ and $z$ is a conformal
parameter on it.

Then there are vector functions $\psi$ and $\varphi$ and a
function $U$ on $\C$ such that

1) $\psi$ and $\varphi$ give a Weierstrass representation of
$\Sigma$;

2) the potential $U$ of this representation is $\Lambda$-periodic;

3) functions $\psi$, $\varphi$, and $U$ meeting 1) and 2) are
defined up to gauge transformations
\begin{equation}
\label{gauge3} \left(
\begin{array}{c}
\psi_1 \\ \psi_2
\end{array}
\right) \to \left(
\begin{array}{c}
e^{h} \psi_1 \\ e^{\bar{h}} \psi_2
\end{array}
\right),
\\
\left(
\begin{array}{c}
\varphi_1 \\ \varphi_2
\end{array}
\right) \to \left(
\begin{array}{c}
e^{-h} \varphi_1 \\ e^{-\bar{h}}\varphi_2
\end{array}
\right),
\\
U \to e^{\bar{h} -h}U
\end{equation}
where
$$
h(z) = a + bz, \ \ \Im (b\gamma) \in \pi \Z \ \ \mbox{for all
$\gamma \in \Lambda$}.
$$
\end{theorem}

As in the case of surfaces in $\R^3$ in general the vector
functions $\psi$ and $\varphi$ define an immersion of the
universal covering surface $\widetilde{\Sigma}$ of a surface
$\Sigma$ into $\R^4$.

\begin{proposition}
An immersion of $\widetilde{\Sigma}$ converts into an immersion of
$\Sigma$ if and only if
\begin{equation}
\label{period4} \int_{\Sigma} \bar{\psi}_1\bar{\varphi}_1 d\bar{z}
\wedge \omega = \int_{\Sigma} \bar{\psi}_1 \varphi_2 d\bar{z} \wedge
\omega = \int_{\Sigma} \psi_2 \bar{\varphi}_1 d\bar{z} \wedge \omega
= \int_{\Sigma} \psi_2 \varphi_2 d\bar{z} \wedge \omega =0
\end{equation}
for any holomorphic differential $\omega$ on $\Sigma$.
\end{proposition}

For $\psi_1 = \pm \varphi_1, \psi_2 = \pm \varphi_2$ the formula
(\ref{period4}) reduces to (\ref{period3}).

\section{Integrable deformations of surfaces}

\subsection{The modified Novikov--Veselov equation}
\label{subsec3.1}

The hierarchy of modified Novikov--Veselov (mNV) equations was
introduced by Bogdanov \cite{Bogdanov1,Bogdanov2} and each
equation from the hierarchy takes the form of Manakov's
``L,A,B''-triple
$$
\frac{\partial L}{\partial t_n} = [L, A_n] - B_n L,
$$
where $L = \D$ is the Dirac operator
$$
L = \left(
\begin{array}{cc}
 0 & \partial \\ -\bar{\partial} & 0
\end{array}\right) + \left(\begin{array}{cc} U & 0 \\ 0 & U
\end{array}\right)
$$
and $A_n$ and $B_n$ are matrix differential operators such that
the~highest term of $A_n$ takes the form
$$
A_n = \left(\begin{array}{cc} \partial^{2n+1} +
\bar{\partial}^{2n+1} & 0 \\ 0 & \partial^{2n+1} +
\bar{\partial}^{2n+1}
\end{array}\right) + \dots\ .
$$
In difference with ``L,A''-pairs, ``L,A,B''-triple preserves only
the zero energy level of $L$ deforming the corresponding
eigenfunctions. Indeed, we have
$$
\frac{\partial L \psi}{\partial t} = L_t \psi + L \psi_t = L[(A +
\partial_t)\psi] - (A+B)[L\psi].
$$
Therefore if $\psi$ meets the equation \beq \label{aequation}
\frac{\partial \psi}{\partial t} + A \psi = 0 \eeq and $L\psi_0 =
0$ for the initial data $\psi_0=\psi|_{t=t_0}$ of this evolutionary
equation then
$$
L \psi = 0
$$
for all $t \geq t_0$.

For $n=1$ we have the original mNV equation \beq \label{mnv} U_t =
\left(U_{zzz} + 3 U_z V + \frac{3}{2} U V_z \right) +
\left(U_{\bar{z}\bar{z}\bar{z}} + 3 U_{\bar{z}} \bar{V} +
\frac{3}{2} U \bar{V}_{\bar{z}}\right) \eeq
where
\beq
\label{mnv-cons} V_{\bar{z}} = (U^2)_z.
\eeq
We see that if the
initial Cauchy data $U\vert_{t=0}$ is a real-valued function then
the solution is also real-valued. In the case when $U\vert_{t=0}$
depends only on $x$ we have $U = U(x,t)$ and the mNV equation
reduces to the modified Korteweg--de Vries equation \beq
\label{mkdv} U_t = \frac{1}{4} U_{xxx} + 6U_x U^2 \eeq (here $V =
U^2$).

This reduction explains the name since Novikov and Veselov had
introduced in \cite{NV1,NV2} a hierarchy of $(2+1)$-dimensional
soliton equations which take the form of ``L,A,B''-triples for
scalar operators with $L =
\partial \bar{\partial} + U$, the two-dimensional Schr\"odinger
operator, and reduces in $(1+1)$-limit to the Korteweg--de Vries
equation. The original Novikov--Veselov equation takes the form
$$
U_t = U_{zzz}+ U_{\bar{z}\bar{z}\bar{z}} + (V U)_z + (\bar{V}
U)_{\bar{z}}, \qquad V_{\bar{z}} = 3  U_z
$$
and its derivation was later modified by Bogdanov for deriving the
mNV equation.

It from the formulas (\ref{int3}) and (\ref{int30}) of the
Weierstrass representation that just the zero energy level of the
Dirac operator relates to surfaces in $\R^3$. This leads to the
following

\bet [\cite{K1}] \label{t-mnv} Let $U(z,\bar{z},t)$ be a
real-valued solution to the mNV equation (\ref{mnv}). Let $\Sigma$
be a surface constructed via the Weierstrass representation
(\ref{int3}) and (\ref{int30}) from $\psi_0$ such that $\psi_0$
meets Dirac equation $\D \psi_0=0$ with the potential $U =
U(z,\bar{z},0)$. Let $\psi(z,\bar{z},t)$ be a solution to the
equation (\ref{aequation}) with $\psi\vert_{t=0}=\psi_0$.

Then the surfaces $\Sigma(t)$ constructed from $\psi(z,\bar{z},t)$
via the Weierstrass representation give a soliton deformation of
the surface $\Sigma$. \eet

The deformation given by this theorem is called {\it the mNV
deformation of a surface}.

Of course, this theorem holds for all equations of the mNV
hierarchy. The recursion formula for them is still unknown and the
next equations are not written explicitly down until recently
except the case $n=2$ \cite{T1}. Finite gap solutions to the mNV
equations are constructed in \cite{T22} (see also \cite{T23}).

It was established in \cite{T1} that this deformation has a global
meaning for tori and preserves the Willmore functional.

\bet [\cite{T1}] \label{tgeo-mnv} The mNV deformation evolves tori
into tori and preserves their conformal classes and the values of
the Willmore functional. \eet

The proof of this theorem is as follows.
To correctly define this
deformation we need to resolve the constraint (\ref{mnv-cons}) and
for tori that can be done globally as it was shown in \cite{T1}.
We have to take a solution $V$ to (\ref{mnv-cons}) normalized by
the condition that
$$
\int_\Sigma V dz \wedge d\bar{z} = 0.
$$
The form $(U^2)_t dz \wedge d\bar{z}$ is an exact form on a torus
$\Sigma$:
$$
UU_t = \bigg(UU_{zz}-\frac{U_z^2}{2}+\frac{3}{2}\,U^2V\bigg)_z+
\bigg(UU_{{\bar z}{\bar z}}-\frac{U_{\bar z}^2}{2}+\frac{3}{2}\,
U^2{\bar{V}}\bigg)_{\bar z}
$$
and therefore the Willmore functional is preserved:
$$
\frac{d}{dt} \int_\Sigma U^2 dz \wedge d\bar{z} = \int_\Sigma
(U^2)_t dz \wedge d\bar{z} = 0.
$$

The flat structure on a torus admits us to identify differentials
with periodic functions. For instance, formally $U^2 dz d\bar{z}$
is a $(1,1)$-differential and $Vdz^2$ is a quadratic differential.
This is impossible for surfaces of higher genus and therefore that
persists to define globally the mNV deformations of such surfaces.

Some attempt to redefine soliton deformations in completely
geometrical terms was done in \cite{BPP}. Finally it did not
manage to avoid introducing a parameter on a surface however some
interesting geometrical properties of the deformations were
revealed.

After papers \cite{K1,T1} in the framework of affine and Lie
sphere geometry some other soliton deformations of surfaces with
geometrical conservation laws were introduced and studied in
\cite{KonPic,Fer,BF}.

\subsection{The modified Korteweg--de Vries equation}
\label{subsec3.2}

In the case when $U$ depends only on $x = \Re z$ the Dirac
equation $\D \psi = 0$ for functions of the form \beq
\label{revsurf} \psi(z,\bar{z}) = \varphi(x)
\exp\left(\frac{iy}{2}\right) \eeq reduces to the Zakharov--Shabat
problem
$$
L \varphi = 0, \ \ \ L = \left[ \left(
\begin{array}{cc}
0 & 1 \\
-1 & 0
\end{array}
\right) \frac{d}{dx} + \left(
\begin{array}{cc}
q &  - ik \\
- ik & q
\end{array}
\right)\right],\ \ \ \ q =2U,
$$
for $k = \frac{i}{2}$.

Notice that for surfaces of revolution the function $\psi$ takes
the form (\ref{revsurf}) in some conformal coordinate $z = x +iy$
where $y$ is the angle of revolution. However there are many other
surfaces with inner $S^1$-symmetry for which the potential $U$
depends only on $x$. The function $\varphi$ is periodic for tori
of revolution and is fast decaying for spheres of revolution
(\cite{T21}).

The operator $L$ is associated with the modified Korteweg--de
Vries hierarchy of soliton equations which admit the ``L,A''-pair
representation
$$
\frac{d L}{dt} = [L, A_n].
$$
The simplest of them is
$$
q_t = q_{xxx} + \frac{3}{2} q^2 q_x, \ \ \ \ n=1,
$$
$$
q_t = q_{xxxxx} + \frac{5}{2}q^2 q_{xxx} + 10 qq_x q_{xx} +
\frac{5}{2} q_x^3 + \frac{15}{8}q^4 q_x, \ \ \ \ n=2.
$$
The first of them coincides with the reduction (\ref{mkdv}) of the
mNV equation after substituting $q \to 4U$ and rescaling the
temporary parameter $t \to 4t$. In fact, the mKdV hierarchy is the
reduction of the mNV hierarchy for $U=U(x)$.

We see that in the mKdV case we have no constraints of type
(\ref{mnv-cons}) and may easily define mKdV deformations of
surfaces of revolution. Moreover in this case there is a recursion
formula for higher equations:
$$
\frac{\partial q}{\partial t_n} = D^n q_x, \ \ \ \ \ D =
\partial_x^2 + q^2 + q_x \partial^{-1}_x q.
$$

Let us introduce the Kruskal--Miura integrals. Their densities
$R_k$ are defined by the following recursion procedure:
$$
R_1 = \frac{iq_x}{2}-\frac{q^2}{4}, \ \ \ R_{n+1} = -R_{nx}
-\sum_{k=1}^{n-1} R_k R_{n-k}.
$$
It is shown that $R_{2n}$ are full derivatives and only the
integrals
$$
H_k = \int R_{2k-1} dx
$$
dot not vanish identically.

\bet [\cite{T12}] For every $n \geq 1$ the $n$-th mKdV equation
transforms (as the reduction of the mNV deformation) tori of
revolution into tori of revolution preserving their conformal
types and the values of $H_k, k \geq 1$. \eet

A proof of the analogous theorem for spheres of revolution (they
are studied in \cite{T21}) is basically the same as for tori.

We note that the preservation of tori is not a trivial fact. This
$L$ operator also comes into the ``L,A''-pair representation of
the sine-Gordon equation and, thus, this equation also induces a
deformation of surfaces of revolution. However this deformation
closes up tori into cylinders.

We see that
$$
H_1 = - \frac{1}{4} \int q^2 dx = - 4 \int U^2 dx = -\frac{2}{\pi}
\int U^2 dx \wedge dy
$$
and therefore the first Kruskal--Miura integral is proportional to
the Willmore functional. The next integrals are
$$
H_2 = \frac{1}{16} \int (q^2 - 4q_x^2) dx, \ \ \ H_3 =
\frac{1}{32} \int (q^6 - 20 q^2 q_x^2 + 8 q_xx^2) dx.
$$
It is interesting

{\sl what are the geometrical meanings of the functionals $H_k$
and what are extremals of these functionals on compact
surfaces of revolution?}

The mKdV deformations of surfaces of revolution determine
deformations of the revolving curves in the upper half plane.
geometry of such deformations and therewith an interplay between
recursion relations and curve geometry were studied in
\cite{Langer,GarayLanger}.

\subsection{The Davey--Stewartson equation}
\label{subsec3.3}

The mNV equations are themselves reductions for $U=-p=\bar{q}$ of
the Davey--Stewartson equations represented by ``L,A,B''-triples
with
$$
L = \left(
\begin{array}{cc} 0 & \partial \\
-\bar{\partial} & 0
\end{array}
\right) + \left(
\begin{array}{cc} -p & 0 \\
0 & q
\end{array}
\right).
$$
Actually this reduction of the Davey--Stewartson (DS) equations
give more equations which are of the form
$$
U_t = i \left(\partial^{2n} U + \bar{\partial}^{2n} U \right) +
\dots
$$
and
$$
U_t = \partial^{2n+1} U + \bar{\partial}^{2n+1} U + \dots
$$
for $n \geq 1$. The first series does not preserve the reality
condition $U = \bar{U}$ and the second series for $U = \bar{U}$
reduces to the mNV hierarchy.

The first two of these equations are the DS$_2$ equation \beq
\label{ds2} U_t = i(U_{zz}+U_{\bar{z}\bar{z}} + 2(V+\bar{V})U)
\eeq where \beq \label{ds2-cons} V_{\bar{z}} = \partial (|U|^2)
\eeq and the DS$_3$  equation (which is sometimes called the
Davey--Stewartson I equation) \beq \label{ds3} U_t =
U_{zzz}+U_{\bar{z}\bar{z}\bar{z}} + 3(V U_z + \bar{V} U_{\bar{z}})
+ 3(W+W^\prime) U \eeq where \beq \label{ds3-cons} V_{\bar{z}} =
(|U|^2)_z, \ \ \ W_{\bar{z}} = (\bar{U}U_z)_z, \ \ \ W^\prime_z =
(\bar{U}U_{\bar{z}})_{\bar{z}}. \eeq

The Davey--Stewartson equations govern soliton deformations of
surfaces in $\R^4$. As for the case for surfaces in $\R^3$ such
deformations were introduced by Konopelchenko who proved in
\cite{K2} the corresponding analog of Theorem \ref{t-mnv}.

However in this case there are two specific problems:

1) As we already mentioned in \S \ref{subsec2.5} the Weierstrass
representation of a surface in $\R^4$ is not unique. Does the DS
deformations of surfaces geometrically different for different
representations?

2) The constraints for the DS equations are more complicated and
how to resolve the constraints (\ref{ds2-cons}) and
(\ref{ds3-cons}) to obtain global deformations of closed surfaces?

These problems we considered in \cite{T5}.

The answer to the first question demonstrates a big difference
from the mNV deformation:

\bi
\item
{\sl the DS deformations are correctly defined only for surfaces
with fixed potentials $U$ of their Weierstrass representations and
for different choices of the potentials such deformations are
geometrically different.} \ei

It would be interesting to understand the geometrical meanings of
these different deformations of the same surface.

The second question is answered by the following analog of Theorem
\ref{tgeo-mnv}:

\bet 1) Given $V$ uniquely defined by (\ref{ds2-cons}) and the
normalization condition $\int V dz \wedge d\bar{z}=0$, the DS$_2$
equation induces deformation of tori into tori preserving their
conformal classes and the values of the Willmore functional.

2) For \beq \label{ds3-cons-sol} V_{\bar{z}} = (|U|^2)_z, \ \ \
\int V dz \wedge d\bar{z} = 0, \ \ \ W = \partial
\bar{\partial}^{-1} (\bar{u}u_z), \ \ \ W^\prime = \bar{\partial}
\partial^{-1} (\bar{u}u_{\bar{z}}) \eeq the DS$_3$ equation
governs a deformation of tori into tori which preserves their
conformal classes and the Willmore functional. \eet

The surface is deformed via deformations of $\psi$ and $\varphi$
and such deformations involve the operators $A$ from the
``L,A,B''-triple. There much more additional potentials coming in
$A$ and the DS equations as it is explained in \cite{K2}. We do
not explain here the reductions in the formula for $A$ which are
necessary to save closedness of surfaces under deformations. We
only mention that the formula (\ref{ds3-cons}) defines the
periodic potentials $W$ and $W^\prime$ up to constants and the
formula (\ref{ds3-cons-sol}) normalizes these constants. This
normalization is necessary for preserving the Willmore functional.
The resolution of the constraints is exposed in \cite{T5} and we
refer to this paper for all details.

\section{Spectral curves}
\label{sec2}

\subsection{Some facts from functional analysis}
\label{subsec4.1}

Given a domain $\Omega \subset \R^n$, denote by $L_p(\Omega)$ and
$W^k_p$ the Sobolev spaces which are the closures of the space of
finite smooth functions on $\Omega$ with respect to the norms
$$
\|f\|_p = \int_{\Omega} |f(x)|^p dx_1 \dots dx_n
$$
and
$$
\|f\|_{k,p} = \sum_{0 \leq l_1+ \dots +l_n =l \leq k}
\int_{\Omega} \left|\frac{\partial^l f}{\partial^{l_1} x_1 \dots
\partial^{l_n} x_n}\right|^p dx_1 \dots dx_n.
$$
For a torus $T^n = \R^n/\Lambda$ we denote by $L_p(T^n)$ and
$W^k_p(T^n)$ the analogous Sobolev spaces formed by
$\Lambda$-periodic functions. Therewith the integrals in the
definitions of norms are taken over compact fundamental domains of
the translation group $\Lambda$.

\begin{proposition}
\label{proposition1} Given a compact closed domain $\Omega$ in
$\R^n$ or a torus, we have

\begin{itemize}
\item
(Rellich)
there is a natural continuous embedding
$
W^k_p(\Omega) \to L_p(\Omega)
$
which is compact for $k > 0$;

\item
(H\"older)
a multiplication by $u \in L_p$ is a bounded operator
from $L_q$ to $L_r$ with
$
\|uv\|_r \leq \|u\|_p \|v\|_q, \ \ \frac{1}{p} + \frac{1}{q} =
\frac{1}{r};$

\item
(Sobolev)
there is a continuous embedding
$W^1_p(\Omega) \to L_q(\Omega), \ \ q \leq \frac{np}{n-p}$
whose norm is called the Sobolev constant;

\item
(Kondrashov)
for $q < \frac{np}{n-p}$ the Sobolev embedding is compact.
\end{itemize}
\end{proposition}

We shall denote the space of two-component vector functions on a
torus $M = \R^2/\Lambda$ by
$$
L_p =  L_p(M) \times L_p(M), \ \ \ \
W^k_p =  W^k_p(M) \times W^k_p(M)
$$
in difference with the spaces of scalar functions $L_2(M)$ and
$W^1_p(M)$.

Let $H$ be a Hilbert space. An operator $A: H \to H$ is  compact
if for the unit ball $B = \{|x|<1 : x \in H\}$ the closure of its
image $A(B)$ is compact. The spectrum $\Spec A$ of a compact
operator $A$ is bounded and can have a limit point only at zero.

Given a Hilbert space $H$ and an operator $A$ (not necessarily
bounded) denote by $R(\lambda)$ the resolvent of $A$. It is a
operator pencil :
$$
R(\lambda) = (A-\lambda)^{-1}
$$
with singularities at $\Spec A$ and holomorphic in
$\lambda$ outside $\Spec A$.

The Hilbert identity reads
\begin{equation}
\label{hilbert}
R(\mu)R(\lambda) = \frac{1}{\mu-\lambda}(R(\lambda) - R(\mu)),
\end{equation}
or in another notation it is
$$
\frac{1}{A-\mu} \frac{1}{A-\lambda}  =
\frac{1}{\mu - \lambda} \left(\frac{1}{A-\lambda} -
\frac{1}{A-\mu}\right).
$$
Given a resolvent defined in some domain in $\C$, we may extend it
onto $\C$ by using the following consequence of the Hilbert
identity:
$$
R(\mu) = R(\lambda) ((\mu-\lambda) R(\lambda) +1)^{-1}
$$
(notice that $R(\lambda) R(\mu) = R(\mu) R(\lambda)$).

\begin{proposition}
\label{proposition2}
If $R(\lambda)$ is compact for $\lambda = \lambda_0$ and holomorphic in
$\lambda$ near $\lambda_0$ then

1) $R(\mu)$ is compact for any $\mu \in \C \setminus \Spec A$
and the resolvent has poles in points from $\Spec A$;

2) $R(\lambda)$ is holomorphic in $\C \setminus \Spec A$.
\end{proposition}

\subsection{The spectral curve of the Dirac operator with boun\-ded poten\-tials}
\label{subsec4.2}

In this section we explain the scheme of proving the existence of
a spectral curve of the differential operator with periodic
coefficients which we used \cite{T2} for the case of Dirac
operators with bounded potentials. This case covers all Dirac
operators corresponding to immersed tori in $\R^3$.

Let
$$
\D =
\left(
\begin{array}{cc}
0 & \partial \\
-\bar{\partial} & 0
\end{array}
\right)+
\left(
\begin{array}{cc}
U & 0 \\
0 & V
\end{array}
\right) =
\D_0 +
\left(
\begin{array}{cc}
U & 0 \\
0 & V
\end{array}
\right).
$$
Here we denote by $\D_0$ the free Dirac operator:
\begin{equation}
\label{diracfree}
\D_0 =
\left(
\begin{array}{cc}
0 & \partial \\
-\bar{\partial} & 0
\end{array}
\right).
\end{equation}

A Floquet eigenfunction $\psi$ of the operator $\D$ with the
eigenvalue (or the energy) $E$ is a formal solution to the equation
$$
\D \psi = E \psi
$$
which satisfies the following periodicity conditions:
$$
\psi(z+\gamma_j,\bar{z}+\overline{\gamma}_j) = e^{2\pi i
(k,\gamma_j)} \psi(z,\bar{z}) = \mu(\gamma_j)\psi(z,\bar{z}), \ \
j=1,2,
$$
where
$$
(k,\gamma_j) = k_1 \gamma_j^1 + k_2 \gamma_j^2, \ \ \gamma_j = \gamma_j^1 + i \gamma_j^2 \in \C = \R^2, \ \
k = (k_1,k_2).
$$
The quantities $k_1,k_2$ are called the quasimomenta of $\psi$ and
$(\mu_1,\mu_2) = (\mu(\gamma_1)$, $\mu(\gamma_2))$ are the
multipliers of $\psi$.

Let us represent a Floquet eigenfunction $\psi$ as a product
$$
\psi(z,\bar{z}) = e^{2\pi i (k_1 x + k_2 y)} \varphi(z,\bar{z}),
\ \ z=x+iy, \ x,y \in \R,
$$
with a $\Lambda$-periodic function $\varphi(z,\bar{z})$. The
equation $\D \psi = E\psi$  takes the form
$$
\left[\left(
\begin{array}{cc}
0 & \partial \\
-\bar{\partial} & 0
\end{array}
\right) +
\left(
\begin{array}{cc}
U & \pi i (k_1 - ik_2) \\
-\pi i(k_1 + i k_2) & V
\end{array}
\right)\right]
\left(
\begin{array}{c}
\varphi_1 \\ \varphi_2
\end{array}
\right) =
E
\left(
\begin{array}{c}
\varphi_1 \\ \varphi_2
\end{array}
\right).
$$
We have an operator pencil
\begin{equation}
\label{pencil}
\D(k) = \D + T_k
\end{equation}
where
\begin{equation}
\label{tk}
T_k =
\left(
\begin{array}{cc}
0 & \pi i (k_1 - ik_2) \\
-\pi i (k_1 + ik_2) & 0
\end{array}
\right).
\end{equation}
This pencil depends analytically on parameters $k_1,k_2$.

We see that to find a Floquet eigenfunction $\psi$ with the quasimomenta
$k_1,k_2$ and the energy $E$ is the same as to find
a periodic solution $\varphi$ to the equation
$$
\D(k) \varphi = E\varphi.
$$
We consider solutions to this equation from $L_2$.

We take a value of $E_0$ such that the operator $(\D_0 - E_0)$ is
inverted on $L_2$, i.e. there exists the inverse operator
$$
(\D_0-E_0)^{-1}: L_2 \to W^1_2.
$$
We represent $\varphi$ in the form
$$
\varphi = (\D_0-E_0)^{-1} f
$$
and substitute this expression into the equation
$$
(\D(k)-E)\varphi = 0
$$
arriving at the following equation:
$$
(1 + A(k,E)) f =0, \ \ \ f \in L_2,
$$
with
$$
A(k,E) =
\left(
\begin{array}{cc}
U + (E_0 - E) & \pi i(k_1 - ik_2) \\
-\pi i (k_1 + ik_2) & V +(E_0 - E)
\end{array}
\right)
(\D_0 - E_0)^{-1} =
$$
$$
= B(k,E)(\D_0 - E_0)^{-1}.
$$

Finally the problem of existence of Floquet functions
with the quasimomenta $k$ and the energy $E$ reduces
to the solvability of
the equation
$$
(1 + A(k,E))f = 0
$$
in $L_2$. Let us notice that the operator $A(k,E)$ is decomposed in the
following chain of operators:
\begin{equation}
\label{dec}
L_2 \stackrel{(\D_0-E)^{-1}}{\longrightarrow} W^1_2
\stackrel{\mathrm{embedding}}{\longrightarrow} L_2
\stackrel{\mathrm{multiplication}}{\longrightarrow} L_2.
\end{equation}
The first mapping is continuous, the second mapping is compact, and,
assuming that the potentials $U$ and $V$ are bounded, the third mapping
which is the multiplication by $B(k,E)$
is continuous. Therefore, we have

\begin{proposition}
Given bounded potentials $U$ and $V$, the analytic pencil of
operators $A(k,E): L_2 \to L_2$
consists of compact operators.
\end{proposition}

Now we can use the Keldysh theorem \cite{Keld1,Keld2} which is the
Fredholm alternative for analytic operator pencils of the form $[1 +
A(\mu)]$ where $A(\mu)$ is a compact operator for every $\mu$. It
reads that

\begin{itemize}
\item
the resolvent of
a pencil $[1+A(\mu)]:H \to H$
where $A(\mu)$ is an analytic pencil of compact operators
is a meromorphic function of $\mu$. Its singularities
which correspond to solutions of the equation $(1+A(\mu))f = 0$
form an analytic subset $Q$ in the space of parameters $\mu$.
\end{itemize}

In the sequel we consider only Floquet functions with $E=0$.

For the operator $\D$ with potentials $U,V$ we have $\mu=(k,E)
\in\ \C^3$ and we put
\begin{equation}
\label{keld}
Q_0 (U,V) = Q \cap \{E=0\}.
\end{equation}
This set is invariant under translations by vectors from the dual
lattice $\Lambda^{\ast} \subset \R^2 = \C$:
$$
k_1 \to k_1 + \eta_1, \ \ \
k_2 \to k_2 + \eta_2.
$$
We recall  that the dual lattice consists of vectors
$\eta = \eta_1 + i \eta_2$
such that $(\eta,\gamma) = \eta_1\gamma^1 + \eta_2\gamma^2$ for any
$\gamma = \gamma^1 + i \gamma^2 \in \Lambda$.

The spectral curve is defined as
$$
\Gamma = Q_0(U,V)/\Lambda^\ast.
$$

{\sc Remark.}
It is easy to notice that the composition of the operator
$$
(\D(k) - E)^{-1} = (\D_0 - E_0)^{-1}(1+A(k,E))^{-1}: L_2 \to W^1_2
$$
and the canonical embedding $W^1_2 \to L_2$ is the resolvent
$R(k,E)$ of
the operator
$$
\D(k) = \D +
\left(
\begin{array}{cc}
0 & \pi i (k_1-ik_2) \\
-\pi i(k_1+ik_2) & 0
\end{array}
\right).
$$
The intersection of the set of poles of $R(k,E)$ with the plane
$E=0$ is the set $Q_0(U,V)$.

We arrive at the following definitions:

\begin{itemize}
\item
the {\it spectral curve} $\Gamma$ of the operator $\D$ with
potentials $U$ and $V$ is the complex curve $Q(U,V)/\Lambda^\ast$
considered up to biholomorphic equivalence;

\item
on $\Gamma$ there is defined the {\it multiplier mapping}, which
is a local embedding near a generic point:
$$
{\cal M}: \Gamma \to \C^2 \ \ : \ \ {\cal M}(k) = (\mu_1,\mu_2) =
(e^{2\pi i (k,\gamma_1)}, e^{2\pi(k,\gamma_2)}),
$$
where $\gamma_1,\gamma_2$ are generators of $\Lambda \subset \C$
and $(k,\gamma_j) = k_1 \Re \gamma_j + k_2 \Im \gamma_j, j=1,2$.
\footnote{This a mapping depends on a choice of generators
$\gamma_1,\gamma_2$. if the basis $\gamma_1,\gamma_2$ is replaced
by another basis $\widetilde{\gamma}_1 = a \gamma_1 + b \gamma_2,
\widetilde{\gamma}_2 = c \gamma_1 + d \gamma_2$, then ${\cal M} =
(\mu_1,\mu_2)$ is transformed as follows \beq \label{changebasis}
{\cal M} \to \widetilde{\cal M} = (\mu_1^a \, \mu_2^b, \mu_1^c \,
\mu_2^d). \eeq};

\item
to every point of $\Gamma$ there is attached the space of Floquet
functions with given multipliers. The dimension of such spaces, in
general, jumps at singular points of $\Gamma$.
\ei

\begin{proposition}
\label{involutions}
Let $k =(k_1,k_2)$ be the quasimomenta of a
Floquet function of $\D$.

1) If $U = \bar{V}$, $\Gamma$ admits an antiholomorphic involution
$\tau: k \to -\bar{k}$.

2) If $U = \bar{U}$ and $V = \bar{V}$, $\Gamma$ admits an
antiholomorphic involution $k \to \bar{k}$.

3) If $U= \bar{U} = V$, then the composition of involutions from
1) and 2) gives a holomorphic involution $\sigma: k \to -k$.
\end{proposition}

Such conditions are usual for spectral curves (see, for instance,
the case of a potential Schr\"odinger operator in \cite{NV1,NV2})
and for the Dirac operator are explained in
\cite{Schmidt,T22,T23}. The simplest of them is the first one
which is proved by the following evident lemma.

\begin{lemma}
If $U = \bar{V}$, then the transformation
$\varphi \to \varphi^\ast$ given by (\ref{ast})
maps Floquet functions into Floquet functions changing the
quasimomenta as follows: $k \to -\bar{k}$.
\end{lemma}

Let us denote by $\Gamma_\nm$ the normalization of $\Gamma$. The
Riemann surface $\Gamma$  is not algebraic but a complex space for
which the existence of a normalization was proved in \cite{GR}.
Since we are in a one-dimensional situation all singular points
are isolated and the normalization is as follows:

1) if a point $P \in \Gamma$ is reducible, i.e. several branches
of $\Gamma$ intersect at $P$, then these branches are unstacked;

2) for an irreducible  singular point $P$ the normalization
$\Gamma_\nm \to \Gamma$ is a local homeomorphism near $P$ given in
terms of local parameters by some series
$$
k_1 = t^a + \dots, \ \ \ k_2 = t^b + \dots, \ \ a>1, \ b>1.
$$
Here $t$ is a local coordinate near $P$ on $\Gamma_\nm$.

If there are no reducible singular points then the normalization map
$\Gamma_\nm \to \Gamma$ is a homeomorphism.

The genus of the complex curve $\Gamma_\nm$ is called the
geometric genus of $\Gamma$ and is denoted by $p_g(\Gamma)$. It is
said that an operator is {\it finite gap} (on the zero energy)
level if $p_g(\Gamma) < \infty$.

The analog of the arithmetic genus for $\Gamma$ which comes into
theorems of the Riemann--Roch type is always infinite:
$p_a(\Gamma) = \infty$.

We have

\bi
\item
nonsingular points of the normalized spectral curve
$\Gamma_{\mathrm{nm}}$ parameterize (up to multiples) the Floquet
functions $\psi$, $\D \psi = 0$. In difference with $\Gamma$ the
one-to-one parameterization property fails only in finitely many
singular points. \footnote{This follows from the asymptotical
behavior of the spectral curve (see \S \ref{subsec4.3}).}
\end{itemize}

In \S \ref{subsec4.7} we argue that in the case when the genus of
$\Gamma_\nm$ is finite it is better to replace $\Gamma_\nm$ by a
curve $\Gamma_\psi$ whose definition involves the Baker--Akhiezer
function of $\D$.

{\sc Example. The spectral curve for $U=V=0$ (the free operator).} For simplicity, we
assume that $\Lambda = \Z + i \Z$. The Floquet functions are as
follows
$$
\psi^+ = (e^{\lambda_+ z},0), \ \ \psi^- = (0,e^{\lambda_-
\bar{z}})
$$
and are parameterized by a pair of complex line with parameters
$\lambda_+$ and $\lambda_-$. These complex lines form the
normalized spectral curve $\Gamma_{\rm nm}$. Since it is of finite
genus we compactify it by two points at infinities such that
$\psi$ has exponential singularities at these points. The
quasimomenta of of these functions are
$$
k_1 = \frac{\lambda_+}{2\pi i} + n_1, \ \ k_2 = \frac{\lambda_+}{2\pi}
+ n_2, \ \  \ \mbox{for $\psi^+$},
$$
$$
k_1 = \frac{\lambda_-}{2\pi i} + m_1, \ \ k_2 =
-\frac{\lambda_-}{2\pi} + m_2, \ \ \ \mbox{for $\psi^-$},
$$
where $m_j,n_j \in \Z$.
The functions $\psi^+$ and
$\psi^-$ have the same multipliers at the points
$$
\lambda_+^{m,n} = \pi(n+im), \ \ \lambda_-^{m,n} = \pi(n-im), \ \
m,n \in \Z,
$$
which form the resonance pairs. The complex curve $\Gamma$ is
obtained from two complex lines after the pair-wise identification
of points from resonance pairs.

{\sc Remark. Spectral curve and the Kadomtsev--Petviashvili
equation.} We exposed above the scheme which we used for defining
the spectral curves of differential operators with periodic
coefficients in 1985 (this paper was never published although it
is referred in \cite{Krichever}). Very similar scheme as we had
known later was used by Kuchment \cite{Ku} (see also \cite{Ku2}).
However some observation about the Kadomtsev--Petviashvili
equations done in that time is worth to be mentioned. Actually
there are two Kadomtsev--Petviashvili (KP) equations
$$
\partial_x ( u_t + 6 uu_x + u_{xxx}) = -3\varepsilon^2 u_{yy}
$$
with $\varepsilon^2 = \pm 1$. For $\varepsilon = i$ it is called
the KPI equation and for $\varepsilon=1$ it is called the KPII
equation. From the point of view of physics these equations are
drastically different. Both these equations admit similar
``L,A''-pair representations $\dot{L} = [L,A]$ with the $L$
operator
$$
L = \varepsilon \partial_y + \partial^2_x + u.
$$
Here the potential $u$ is double-periodic or, which is the same,
defined on some torus $\R^2/\Lambda$. The free operator equals
$L_0 = \varepsilon \partial_y - \partial^2$ and to prove the
existence of the spectral curve by the scheme used above we need
to take the inverse operator
$$
(L_0 - E_0)^{-1}: L_2 \to W^{2,1}_2
$$
where $W^{2,1}_2$ is the space of functions on the torus such that
$u, u_x, u_{xx}$, and $u_y$ lie in $L_2$.
For simplifying computations, we consider the case when $\Lambda$ is
generated by $(2\pi,0)$ and $(0,2\pi\tau^{-1})$.
Then the Fourier basis in $L_2$ is formed by the functions
$$
e^{i(kx +l\tau y)}, \ \ \ k,l \in \Z.
$$
In this basis the operator $(L_0 - E_0)$ is diagonal and we have
$$
(L_0 - E_0)e^{i(kx +l\tau y)} = (i \varepsilon l \tau - k^2 - E_0)e^{i(kx +l\tau y)}.
$$
Since $\varepsilon = 1$ for KPII, we have for $E_0 > 0$ a bounded
operator
$$
(L_0 - E_0)^{-1}e^{i(kx +l\tau y)} = \frac{1}{i l \tau - k^2 - E_0}e^{i(kx +l\tau y)}.
$$
It is easy to check that if $\varepsilon = i$ then for any $E_0$
either the operator $(L_0-E_0)$ is not inverted or its inverse is
unbounded. That takes the case for any lattice $\Lambda$. One can
deduce from these reasonings that for the operator $L =
i\partial_y + \partial_x^2 + u$ the spectral curve does not exist.
For the heat operator $L = \partial_y + \partial_x^2 + u$ it
exists and is preserved by the KPII equation.

The spectral curve of a two-dimensional periodic differential
operator $L$ on the zero energy level was first introduced in the
paper by Dubrovin, Krichever, and Novikov \cite{DKN} in the case
of Schr\"odinger operator, where it is showed that

1) the periodic operator which is finite gap on the zero energy
level is reconstructed from some algebraic data including this
curve;
\footnote{For the Dirac operator $\D$ see the reconstruction formula
\ref{reconstruction} and its derivation in \cite{T4} and \S \ref{subsec4.7}.}

2) this curve is the first integral of the deformations of $L$
governed by the ``L,A,B''-triples.

\begin{proposition}
[\cite{DKN}] Let $L$ be a two-dimensional periodic differential
operator, $\Gamma$ be its spectral curve, and ${\cal M}$ be the
multiplier mapping.

Let we have the evolution equation
$$
\frac{\partial L}{\partial t} = [L,A]-BL
$$
such that the operator $A$ is also periodic. Then this deformation
of $L$ preserves $\Gamma$ and ${\cal M}$.
\end{proposition}

This result generalizes the conservation law for the spectral
curve of a one-dimensional operator $L$ deformed via the
``L,A''-pair type equation $\frac{\partial L}{\partial t} = [L,A]$
(this was first established for the periodic KdV equation by
Novikov in \cite{Novikov}).

This proposition follows from the deformation equation $\psi_t +
A\psi = 0$ for the Floquet functions which preserves the
multipliers (see \S \ref{subsec3.1} and the equation
(\ref{aequation})). The conservation of the zero level spectrum
was first indicated by Manakov in \cite{Manakov} where the
``L,A,B''-triples were introduced.

\begin{corollary}
The spectral curve $\Gamma$ and the multiplier mapping ${\cal M}$
of the periodic Dirac operator $\D$ are preserved by the modified
Novikov--Veselov and Davey--Stewartson equations.
\end{corollary}

For the mKdV deformations we have two spectral curves: $\Gamma$
defined for the two-dimensional Dirac operator and $\Gamma^\prime$
defined for a one-dimensional operator $L_{\rm mKdV}$ which comes
into the Zakharov--Shabat problem (see \S \ref{subsec3.2}) and the
``L,A''-pair representation for the mKdV equation. These complex
curves are related by the canonical branched two-covering $\Gamma
\to \Gamma_0$ \cite{T22} and both of them are by the mKdV
equation. The complex curve $\Gamma_0$ is uniquely reconstructed
from the Kruskal--Miura integrals $H_k, k=1,\dots$, which are also
first integrals of the mKdV equation.

\subsection{Asymptotic behavior of the spectral curve}
\label{subsec4.3}

The spectral curve of $\D$ is a perturbation of the spectral curve
of the free operator $\D_0$. Although this perturbation could be
rather strong in a bounded domain $|k| \leq C$, outside this
domain it results just in a transformation of double points
corresponding to resonance pairs into handles. Moreover the size
of a handle is decreasing as $|k| \to \infty$ and is estimated in
terms of the perturbation.

Thus we have

1) a compact part $\Gamma_0 = Q_0 \cap \{|k| \leq C\}$ whose
boundary consists in a pair of circles;

2) a complex curve $\Gamma_\infty$ obtained from the planes $k_1 =
ik_2$ and $k_1 = -ik_2$ by removing the domains with $\{|k| \leq
C\}$ and transforming some of double points corresponding to
resonance pairs into handles;

3) $\Gamma_0$ and $\Gamma_\infty$ are glued along their
boundaries;

4) $\Gamma$ has two ends at which ${\cal M}(\Gamma)$ behaves
asymptotically as in the case of the free operator.

This complex curve is the curve obtained from $\Gamma$ by
unstucking double points which correspond to resonant pairs and
survive the perturbation. We denote it again by $\Gamma$.

The operator is finite gap (on the zero energy level) if under the
perturbation $\D_0 \to \D$ only finitely many double points are
transformed into handles.

This picture is typical in soliton theory where the spectral curve
of some operator with potentials is a perturbation of the spectral
curve of the corresponding  free operator and therewith the
perturbation is small for large values of quasimomenta. It was
rigorously established for the two-dimensional Sch\"rodinger
operator by Krichever \cite{Krichever} who used perturbation
theory. In \cite{T3} we proposed to clarify this geometrical
picture for the Dirac operator by using same methods and
formulated the expected statement as Pretheorem.

The theory of spectral curves initiated developing of the analytic
theory of Riemannian surfaces (not only hyperelliptic) of infinite
genus in \cite{FKT0,FKT}.

In \cite{Schmidt} Schmidt proposed another approach to confirm
this asymptotic behavior of the spectral curve. It is based on his
result on the existence of the spectral curves for the Dirac
operators with $L_2$ potentials and the continuiosity of these
curves for weakly converging sequence of potentials.

\begin{theorem}
[\cite{Schmidt}]
\label{l2}
Given $U,V \in L_2(T^2)$, the equation
$$
\D(k) \varphi = (\D + T_k) \varphi = E \varphi
$$
with $k \in \C^2, E\in \C$, has a solution in $L_2$ if and only if
$(k,E) \in Q$ where $Q$ is an analytic subset in $\C^3$. This
subset $Q$ is is formed by poles of the operator pencil
$$
(1 + A_{U,V}(k,E))^{-1}: L_2 \to L_2
$$
where $A_{U,V}(k,E)$ is polynomial in $k,E$. Moreover if
$$
U_n,V_n \stackrel{\mathrm{weakly}}{\longrightarrow}
U_\infty,V_\infty
$$
in $\{\|U\|_{2;\varepsilon} \leq C, \|V\|_{2;\varepsilon} \leq C\}$
\footnote{ Recall that a sequence $\{u_n\}$ in a Hilbert space $H$
weakly converges to $u_\infty$: $u_n
\stackrel{\mathrm{weakly}}{\longrightarrow} u_\infty$ if for any $v
\in H$ we have $\lim_{n\to \infty} \langle u_n, v \rangle = \langle
u_\infty, v \rangle$ where $\langle u,v\rangle$ is the Hilbert
product in $H$.} then
$$
\|A_{U_n,V_n}(k,E) - A_{U_\infty,V_\infty}(k,E)\|_2 \to 0
$$
uniformly near every $k \in \C^2$.
\end{theorem}

We expose the proof of this theorem in Appendix 1. Let us return
to the asymptotic behavior of the spectral curve.

First note the following identity which is checked by
straightforward computations:
\begin{equation}
\label{dress}
\left(
\begin{array}{cc}
e^{-a} & 0 \\
0 & e^{-b}
\end{array}
\right)
\left(
\D_0 +
\left(
\begin{array}{cc}
U & 0 \\ 0 & V
\end{array}
\right) + T_k
\right)
\left(
\begin{array}{cc}
e^{b} & 0 \\
0 & e^{a}
\end{array}
\right) =
\end{equation}
$$
=
\D_0 +
\left(
\begin{array}{cc}
e^{b-a} U & 0 \\ 0 & e^{a-b} V
\end{array}
\right) +
T_k +
\left(
\begin{array}{cc}
0 & a_z \\ -b_{\bar{z}} & 0
\end{array}
\right)
$$
for all smooth functions $a,b: \C \to \C$.

For any $\kappa=(\kappa_1,\kappa_2) \in \Lambda^\ast \subset \C$ define
$\Lambda$-periodic functions
$$
\psi_{\pm \kappa}(z,\bar{z}) = e^{\pm 2\pi i (\kappa_1 x + \kappa_2 y)}
$$
and take the functions $a(z,\bar{z})$ and $b(z,\bar{z})$ in the form
$$
a(z,\bar{z}) = 2\pi i (\alpha_1 x + \alpha_2 y), \ \
b(z,\bar{z}) = 2\pi i ((\alpha_1 - \kappa_1) x + (\alpha_2-\kappa_2) y),
$$
where
$$
\alpha(\kappa) = (\alpha_1,\alpha_2) =
\left(
\frac{\kappa_1 + i\kappa_2}{2},
\frac{-i\kappa_1+\kappa_2}{2}
\right).
$$
The following equalities are clear: $e^{b-a} = \psi_{-\kappa}, \ \
\ a_z = b_{\bar{z}} = 0$. That together with (\ref{dress}) implies

\begin{proposition}
[\cite{Schmidt}]
\label{dress2}
If $\varphi \in L_2$ satisfies the equation
$$
\left[\D_0 +
\left(
\begin{array}{cc}
\psi_{-\kappa} U & 0 \\ 0 & \psi_\kappa V
\end{array}
\right) +
T_k\right] \varphi = 0
$$
then
$\varphi^\prime = \left(\begin{array}{cc} \psi_{-\kappa} & 0 \\ 0 & 1
\end{array}\right)\varphi \in L_2$
meets the equation
$$
(\D + T_{k + \alpha})\varphi^\prime = 0.
$$
Therefore
$$
Q_0(\psi_{-\kappa}U,\psi_\kappa V) = Q_0(U,V) + \alpha(\kappa)
\ \ \ \mbox{for all} \ \ \kappa \in \Lambda^\ast
$$
(here the right-hand side denotes $Q_0(U,V)$ translated by $\alpha$).
\end{proposition}

The functions $\psi_\kappa, \kappa \in \Lambda^\ast$, form a Fourier
basis for $L_2$. The mapping $U \to \widehat{U} = \psi_{\kappa}
U, U = \sum_{\nu \in \lambda^\ast} U_\nu \psi_\nu$, shifts the
Fourier coefficients of $U$: $\widehat{U}_\nu = U_{\nu - \kappa}$.
Therefore, we have
$$
\psi_\kappa U \stackrel{\mathrm{weakly}}{\longrightarrow} 0 \ \ \
\mbox{as $|\kappa| \to \infty$}.
$$
Theorem \ref{l2} (see Appendix 1) and Proposition \ref{dress2}
imply that in a small bounded neighborhood $O(k)$ of $k \in \C^2$
for large $|\kappa|$ the intersection $Q_0(U,V)$ with $O(k) +
\alpha(\kappa)$ is very closed to the intersection of $Q_0(0,0)$
with $O(k)$:
$$
Q_0(U,V) \cap \left[ O(k) + \alpha(\kappa) \right] \approx
Q_0(0,0) \cap O(k) \ \ \ \mbox{as $|\kappa| \to \infty$}.
$$
We conclude that  asymptotically as $|k| \to \infty$
the spectral curve of $\D$ behaves as the spectral curve of the free
operator $\D_0$ on $L_2$.

For $U=V=0$ the spectral curve $\Gamma$ is biholomorphic
equivalent to a pair of two planes (complex lines) defined in
$\C^2$ by the equations
$$
k_2 = ik_1, \ \ \ k_2 = -ik_1,
$$
and glued at infinitely many pairs of points corresponding to the
so-called resonance pairs
$$
\left(k_1 = \frac{\bar{\gamma}_1 n -
\bar{\gamma}_2 m}{\bar{\gamma}_1\gamma_2 -
\gamma_1 \bar{\gamma}_2}, k_2=ik_1 \right)
\leftrightarrow
\left(k_1 = \frac{\gamma_1 n -\gamma_2 m}{\bar{\gamma}_1\gamma_2 -
\gamma_1 \bar{\gamma}_2}, k_2=-ik_1 \right)
$$
where $m,n \in \Z$. Moreover these planes are naturally completed
by a pair of points $\infty_\pm$ which lie at at infinity and are
obtained  the limit $(k_1,\pm ik_1) \to \infty_\pm$ as $k_1 \to
\infty$. Near generic points there is a double covering $\Gamma
\to \C \ \ : \ \ (k_1,k_2) \to k_1$. By Proposition \ref{dress2},
we have

\begin{corollary}
\label{bound} Given a Dirac operator with $L_2$-potentials,
${\cal M}(\Gamma)$ for sufficiently large $|k|$ asymptotically
behaves as
$$
k_2 \approx \pm ik_1.
$$
Therefore it has at most two irreducible components such that
every component contains at least one of these asymptotic ends.
\end{corollary}

The bound for the number of irreducible components is clear, since
other components have to be localized in a bounded domain of
$\C^2$ which is impossible for one-dimensional analytic sets.

Thus we arrive at the definition compatible with one used in the
finite gap integration \cite{DKN,Krichever0}:

\begin{itemize}
\item
if the spectral curve $\Gamma$ of the operator $\D$ is of finite
genus, then this operator is finite gap and we call the completion
of $\Gamma$ by a pair of infinities $\infty_\pm$ the spectral
curve (of a finite gap operator).
\end{itemize}

We finish with the procedure which reconstruct the value of
$$
\int_{\C/\Lambda} UV dx \wedge dy
$$
from $(\Gamma,{\cal M})$ when $\Gamma$ is of finite genus. Near
the asymptotic end where $k_2 \approx ik_1$ we introduce a local
parameter $\lambda_+^{-1}$ such that the multipliers behave as
$$
\mu(\gamma) = \lambda_+ \gamma + \frac{C_0 \bar{\gamma}}{\lambda_+} +
O(\lambda_+^{-2}).
$$
Then
\begin{equation}
\label{willmoreformula} \int_{\C/\Lambda} UV dx \wedge dy = - C_0
\cdot (\mathrm{Area} (\C/\Lambda)
\end{equation}
(see \cite{GS,T3} for the case $U=V$).

The analogous formula for the area of minimal tori in $S^3$ was 
derived by Hitchin in \cite{Hitchin}.

This formula gives us a reason to treat the pair $(\Gamma,{\cal M})$
as a generalization of the Willmore functional. First that was
discussed for tori of revolution in \cite{T12}. In this case the spectral curve is
reconstructed from infinitely many integral quantities known as the
Kruskal--Miura integrals \cite{T12}.

\subsection{Spectral curves of tori}
\label{subsec4.4}

Given a torus $\Sigma$ immersed into the three-dimensional Lie group $G = \R^3$, $SU(2)$ $= S^3$, $\nil$ or
$\sll$ and its the Weierstrass representation, we take
the spectral curve $\Gamma$ of the operator $\D$ coming in this representation.

We call it {\it the spectral curve of the torus} $\Sigma$.

It is defined for all smooth tori and not only for integrable tori
(see \S \ref{subsec4.6}). This definition was originally
introduced for tori in $\R^3$ in \cite{T2} and for tori in $S^3$
in \cite{T21} in its relation to the physical explanation of the
Willmore conjecture. The formula (\ref{willmoreformula}) shows
that the Willmore functional is reconstructed from $\Gamma$ and
the multiplier mapping ${\cal M}$ (at least in the case when
$\Gamma$ is of finite genus).

This definition does not depend on a choice of a conformal
parameter on the torus $\Sigma = \R^2/\Lambda$. The multiplier
mapping ${\cal M}$ depends on a choice of a basis in $\Lambda$ and
the change of a basis results in a simple algebraic transform of
${\cal M}$ (see (\ref{changebasis})).

Let us define the spectral curve for tori in $\R^4$.

In \cite{T5} we explained that the Weierstrass representation for
a surface in $\R^4$ is not unique. The potentials of different
representations of a torus are related by the formula \beq
\label{r4trans} U \to U \exp{(\bar{a} + \overline{bz} - a -bz)}
\eeq where $\Im b\gamma \in \pi \Z$ for all $\gamma \in \Lambda$.
The multiplier mapping ${\cal M}$ depends on the choice of $U$ and
under the transformation (\ref{r4trans}) it is changed as follows:
$$
\mu(\gamma) \to e^{b\gamma} \mu(\gamma), \ \ \ \ \gamma \in \Lambda.
$$
As in the case of tori in $\R^3$ that the integral squared
norm of the potential $U$ is reconstructed from $(\Gamma,{\cal M})$
by the same formula (\ref{willmoreformula}).

The conformal invariance of the Willmore functional led us to the
conjecture which we justified by numerical experiments in
\cite{T12} and was very soon after its formulation was confirmed
in \cite{GS}:

\begin{theorem}
\label{conforminvar}
Given a torus in $\R^3$, its spectral curve
$\Gamma$ and ${\cal M}$ are invariant under conformal
transformations of $\overline{\R}^3$.
\end{theorem}

The proof from \cite{GS} works rigorously for the spectral curves
of finite genus and is as follows. Let us consider the generators
of the conformal group which is $SO(4,1)$ and write down the
deformation equations for a Floquet function $\varphi$ which are
of the form \beq \label{grin} \D \delta \varphi + \delta U \cdot
\varphi. \eeq It is enough to check the invariance only for
inversions and even just for one of them since all they pairwise
conjugated by orthogonal transformations. We take the following
generator for an inversion:
$$
\delta x^1 = -2x^1 x^3, \ \ \ \delta x^2 = -2x^2 x^3, \ \ \
\delta x^3 = (x^1)^2 + (x^2)^2 -(x^3)^2
$$
and compute the corresponding variation of the potential:
$$
\delta U = |\psi_2|^2 - |\psi_1|^2
$$
where $\psi$ generates the torus. In \cite{GS} for this variation
an explicit formula for a solution to (\ref{grin}) is given in
terms of functions meromorphic on the spectral curve. It follows
from this explicit formula that the multipliers are preserved. For
the spectral curve of finite genus these meromorphic functions are
easily defined. For the case of spectral curve of infinite genus
one needs to clarify some analytical details that as we think can
be done and relates on a rigorous and careful treatment of the
asymptotic behavior of the spectral curve.

Another proof of theorem \ref{conforminvar} for isothermic tori
was done in \cite{T3}. It is geometrical and works for spectral
curves of any genus.

\subsection{Examples of the spectral curves}
\label{subsec4.5}

{\sc Products of circles in $\R^4$.}

We consider the tori $\Sigma_{r,R}$ defined by the equations
$$
(x^1)^2 + (x^2)^2 = r^2, \ \ \ (x^3)^2 + (x^4)^2 = R^2.
$$
They are parameterized by the angle variables $x,y$ defined modulo
$2\pi$:
$x^1 = r \cos x, x^2 = r \sin x,
x^3 = R \cos y, x^4 = R \sin y$.
The conformal parameter, the period lattice  and
the induced metric are as follows:
$$
z = x + i \frac{R}{r}y,
\ \ \
\Lambda = \{2\pi m + i2\pi \frac{r}{R} n \ : \ m,n \in \Z\},
\ \ \
ds^2 = r^2 dzd\bar{z}.
$$
By simple computations we obtain the formula for the Gauss map:
$\frac{a_1}{a_2} = - e^{i(y-x)}, \frac{b_1}{b_2} = e^{-i(y+x)}$.
Let us apply Theorem \ref{theorem4} to the mapping
$$
\Sigma_{r,R} \to (b_1:b_2) = \left(\frac{e^{-i(x+y)}}{\sqrt{2}}:
\frac{1}{\sqrt{2}}\right) \in \C P^1.
$$
We have $g = \frac{i(x+y)}{2}$,
$$
U = \frac{1}{4}\left(\frac{r}{R} + i\right)
$$
and the torus $\Sigma_{r,R}$ is defined via the Weierstrass
representation by vector functions
$$
\psi_1 = \psi_2 = \frac{1}{\sqrt{2}}
\exp \left( -\frac{i(x+y)}{2} \right), \ \
\varphi_1 = -\varphi_2 = -\frac{r}{\sqrt{2}}
\exp \left( \frac{i(y-x)}{2} \right).
$$

The values of the Willmore functional on such tori are given by the formula
$$
{\cal W}(\Sigma_{r,R}) = 4 \int_{\Sigma_{r,R}} |U|^2 dx \wedge dy =
\pi^2 \left(\frac{r}{R} + \frac{R}{r}\right)
$$
and attain their minimum at the Clifford torus $\Sigma_{r,r}$
in $\R^4$: ${\cal W}(\Sigma_{r,r}) = 2\pi^2.$

The spectral curve $\Gamma(u)$
of the Dirac operator
$$
\D =
\left(
\begin{array}{cc}
0 & \partial \\
-\bar{\partial} & 0
\end{array}
\right)
+
\left(
\begin{array}{cc}
u & 0 \\
0 & \bar{u}
\end{array}
\right), \ \ \ u = \const,
$$
with the constant potential $U=u$ is the complex sphere with a pair of
marked points (``infinities'') which are $\lambda=0$ and $\lambda=\infty$:
$$
\Gamma(u) = \C P^1.
$$
The normalized Baker--Akhiezer function (or the Floquet function)
equals
$$
\psi(z,\bar{z},\lambda) =
\left(
\begin{array}{c}
\psi_1 \\
\psi_2
\end{array}
\right)
=
\frac{\lambda}{\lambda-u}
\exp\left(\lambda z - \frac{|u|^2}{\lambda} \bar{z} \right)
\left(
\begin{array}{c}
1 \\
-\frac{u}{\lambda}
\end{array}
\right).
$$
The normalization means that the following asymptotics hold:
$$
\psi \approx
\left(
\begin{array}{c}
e^{\lambda_+ z} \\ 0
\end{array}
\right) \ \ \mbox{as $\lambda_+ \to \infty$},
\ \ \
\psi \approx
\left(
\begin{array}{c}
0 \\ e^{\lambda_- \bar{z}}
\end{array}
\right) \ \ \mbox{as $\lambda_- \to 0$}
$$
with the local parameters
$\lambda_+ = \lambda$ near $\lambda=\infty$
and
$\lambda_- = -\frac{|u|^2}{\lambda}$ near $\lambda=\infty$.

For a torus $\Sigma_{r,R}$ we have

\begin{itemize}
\item
the function $\psi$ generating it via (\ref{int40}) equals to
$\psi(z,\bar{z},-u)$, $u = \frac{1}{4}(\frac{r}{R}+i)$, and its
monodromy is as follows $\psi(z+ 2\pi,\bar{z}-2\pi i,-u) = \psi(z
+ i 2\pi \frac{R}{r},\bar{z} -i 2\pi \frac{R}{r},-u) =
-\psi(z,\bar{z},-u)$;

\item
there are exactly four points on the spectral curve $\Gamma(u)$
for which the function $\psi(z,\bar{z},\lambda)$ has the same monodromy as
$\psi(z,\bar{z},-\lambda)$: these are $\lambda =\pm u, \pm \bar{u}$.
Moreover,
$$
\left(
\begin{array}{c}
\psi_1(z,\bar{z},-u) \\ \psi_2(z,\bar{z},-u)
\end{array}
\right)
=
\left(
\begin{array}{c}
- \bar{\psi}_2(z,\bar{z},u) \\ \bar{\psi}_1(z,\bar{z},u)
\end{array}
\right);
$$

\item
the spectral curve $\Gamma(u)$ is smooth.
\end{itemize}

Here $k_1$ and $k_2$ are the quasimomenta of Floquet functions
$\psi(z,\bar{z},\lambda)$.

A periodic potential $U$ is defined up to the gauge transformation
(\ref{gauge3}) which for $b=0$ and $e^{\bar{a}-a} =
-\frac{1+i}{\sqrt{2}}$ transforms the potential $U$ of the
Clifford torus to the potential
$$
\frac{1}{4} (1+i) \to \frac{e^{\bar{a}-a}}{4}(1+i) = -\frac{i}{2\sqrt{2}}
$$
which coincides with the potential of the same torus
considered as a torus in the unit sphere $S^3 \subset \R^4$ \cite{T3}.
This leads to the following questions:

1) {\sl do the spectral curves of a torus in $S^3 \subset \R^4$
defined as for a torus  in $S^3$ and a torus in $\R^4$ always
coincide?}

2) {\sl given a torus in $S^3 \subset \R^4$, does the potential
$U$ of its Weierstrass representation in $\R^4$ is always gauge
equivalent to the potential of its Weierstrass representation in
$S^3$:
$$
U = \frac{(H-i)e^\alpha}{2}
$$
where $H$ is the mean curvature of this torus in $S^3$?}

A positive answer to the second question implies a positive answer
to the second one. We think that the both questions are answered
positively.

\vskip4mm

\noindent {\sc The Clifford torus in $\R^3$.}

The Clifford torus in $\R^3$ is the image of the Clifford torus in
$S^3 \subset \R^4$ under a stereographic projection
$$
(x^1,x^2,x^3,x^4) \to \left(\frac{x^1}{1-x^4},\frac{x^2}{1-x^4},
\frac{x^3}{1-x^4}\right), \ \ \ \sum_k (x^k)^2 = 1.
$$
It is considered up up to conformal transformations of
$\bar{\R}^3$ and hence can be obtained as the following torus of
revolution: given a circle of radius $r=1$ in the $x^1 x^3$ plane
such that the distance between the the circle center and the $x^1$
axis equals to $R=\sqrt{2}$, the Clifford torus is obtained by a
rotation of this circle around the $x^1$ axis.

\begin{theorem}[\cite{T4}]
\label{cliftheorem}
The Baker--Akhiezer function of the Dirac operator $\D$ with the potential
\begin{equation}
\label{pot}
U = \frac{\sin y}{2\sqrt{2}(\sin y  - \sqrt{2})}
\end{equation}
is a vector function $\psi(z,\bar{z},P)$, where $z \in \C$ and $P
\in \Gamma$ such that

\begin{itemize}
\item
the complex curve $\Gamma$ is a sphere $\C P^1 = \bar{\C}$ with
two marked points $\infty_+ = (\lambda=\infty), \infty_- =
(\lambda=0)$ where $\lambda$ is an affine parameter on $\C \subset
\C P^1$ and with two double points obtained by stacking together
the points from the following pairs:
$$
\left(\frac{1+i}{4},\frac{-1+i}{4}\right)
\ \ \ \mbox{and} \ \ \
\left(-\frac{1+i}{4},\frac{1-i}{4}\right);
$$

\item
the function $\psi$ is meromorphic on $\Gamma \setminus \{\infty_\pm\}$ and
has at the marked points (``infinities'') the following asymptotics:
$$
\psi \approx
\left(
\begin{array}{c}
e^{k_+  z} \\ 0
\end{array}
\right) \ \mbox{as $k_+ = \lambda \to \infty$};
\ \
\psi \approx
\left(
\begin{array}{c}
0 \\ e^{k_-  \bar{z}}
\end{array}
\right) \ \mbox{as $k_- = -\frac{|u|^2}{\lambda} \to \infty$}
$$
where $u=\frac{1+i}{4}$ and $k^{-1}_\pm$ are local parameters near
$\infty_\pm$;

\item
$\psi$ has three poles on $\Gamma \setminus\{\infty_\pm\}$ which are
independent on $z$ and have the form
$$
p_1 = \frac{-1+i + \sqrt{-2i-4}}{4\sqrt{2}}, \ \
p_2 = \frac{-1+i - \sqrt{-2i-4}}{4\sqrt{2}}, \ \
p_3 = \frac{1}{\sqrt{8}}.
$$
\end{itemize}

Therewith the geometric genus $p_g(\Gamma)$ and
the arithmetic genus $p_a(\Gamma)$ of $\Gamma$ are as follows:
$$
p_g(\Gamma) = 0, \ \ \ \
p_a(\Gamma) = 2.
$$

The Baker--Akhiezer function satisfies the Dirac equation $ \D
\psi = 0$ with potential $U$ given by (\ref{pot}) at any point
from $\Gamma \setminus \{\infty_+,\infty_-,p_1,p_2,p_3\}$.

The Clifford torus is constructed via the Weierstrass representation
(\ref{int3}) and (\ref{int30}) from the
function
$$
\psi = \psi\left(z,\bar{z},\frac{1-i}{4}\right).
$$
\end{theorem}

It is showed that $\psi$ has the form
$$
\psi_1(z,\bar{z},\lambda) = e^{\lambda z -\frac{|u|^2}{\lambda}\bar{z}}
\left(
q_1 \frac{\lambda}{\lambda-p_1} + q_2 \frac{\lambda}{\lambda-p_2} +
(1-q_1-q_2)\frac{\lambda}{\lambda-p_3}\right),
$$
$$
\psi_2(z,\bar{z},\lambda) = e^{\lambda z -\frac{|u|^2}{\lambda}\bar{z}}
\left(
t_1 \frac{p_1}{p_1-\lambda} + t_2
\frac{p_2}{p_2-\lambda} +
(1-t_1-t_2)\frac{p_3}{p_3-\lambda}\right)
$$
where $u=\frac{1+i}{4}$ and the functions
$q_1,q_2,t_1,t_2$ depend only on $y$ and $2\pi$-periodic with respect to $y$.
They are found from the following conditions
$$
\psi\left(z,\bar{z},\frac{1+i}{4}\right) =
\psi\left(z,\bar{z},\frac{-1+i}{4}\right), \ \ \
\psi\left(z,\bar{z},-\frac{1+i}{4}\right) =
\psi\left(z,\bar{z},\frac{1-i}{4}\right).
$$

\subsection{Spectral curves of integrable tori}
\label{subsec4.6}

It is said that a surface is integrable if the Gauss--Codazzi
equations is the compatibility condition \beq \label{null}
[\partial_x - A(\lambda), \partial_y - B(\lambda)] = 0 \eeq for
the linear problems
$$
\partial_x \varphi = A(\lambda) \varphi, \ \ \
\partial_y \varphi = B(\lambda) \varphi
$$
such that $A$ and $B$ are Laurent series in a spectral parameter.
It is also assumed that $\lambda$ comes nontrivially into this representation. 
For deriving explicit solutions of the {\it zero
curvature equation} (\ref{null}) equation one can use the
machinery of soliton theory and, in particular, of the theory of
integrable harmonic maps which started with papers
\cite{Pohl,ZM,U} and was intensively developed in last thirty
years (the recent statement of this theory is presented in
\cite{Guest,Harm,HeleinBook}). The most complete list of
integrable surfaces in $\R^3$ is given in \cite{Bob2} (see also
\cite{FGFL}).

This theory works well for spheres when it is enough to apply
algebraic geometry of complex rational curves and for tori when
explicit formulas for surfaces are derived in terms of theta
functions of some Riemann surfaces. For surfaces of higher genus
theory of integrable systems does not lead to a substantial
progress. This probably has serious reasons consisting in that
tori are the only closed surfaces admitting flat metrics.

The spectral curves of integrable tori appear as the spectral
curves of operators coming in these auxiliary linear problems.
These complex curves (Riemann surfaces) serve for constructing
explicit formulas for tori in terms of theta functions of these
Riemann surfaces.

It appears that that is not accidentally and these spectral curves
of integrable tori are just special cases of the general spectral
curve defined in \S \ref{subsec4.4} for all tori (not only
integrable).

In \cite{T3} we proved such a coincidence (modulo additional
irreducible components) for constant mean curvature and isothermic
tori in $\R^3$ and for minimal tori in $S^3$. Corollary \ref{bound}
rules out additional components.

{\sc A) Constant mean curvature (CMC) tori in $\R^3$.} By the
Ruh--Vilms theorem, the Gauss map of a surface in $\R^3$ is
harmonic if and only if this surface has a constant mean
curvature \cite{RV}. By the Gauss--Codazzi equations this is
equivalent to the condition that the Hopf
differential $A dz^2$ is holomorphic: 
$$
A_{\bar{z}} = 0.
$$

On a sphere a holomorphic quadratic differential vanishes and therefore, by
the Hopf theorem, CMC spheres in $\R^3$ are exactly round spheres \cite{Hopf}

It was also conjectured by Hopf that all immersed compact CMC surfaces 
in $\R^3$ are just round spheres. Although this conjecture was confirmed for 
embedded surfaces by Alexandrov \cite{Alexandrov} it was disproved for 
immersed surfaces of higher genera.
The existence of CMC tori was established in the early 1980's by 
Wente by means of the Banach space implicit function theorem. 
The first explicit examples were found by Abresch in \cite{ACMC1}
and the analysis of these examples performed in \cite{ACMC2} gave a hint on 
the relation of this problem to integrable systems. 
It was proved later for a CMC torus the
complex curve $\Gamma$ is of finite genus \cite{PS} and this
allowed to apply the Baker--Akhiezer functions to deriving
explicit formulas for such tori in terms of theta functions of
$\Gamma$ (this program was realized by Bobenko in
\cite{Bobenko,Bob1}). The existence of CMC surfaces of very genus greater 
than one was established by Kapouleas also by implicit 
methods \cite{Kap1,Kap2} and the problem of explicit 
description of such surfaces stays open. 
We remark that some other 
interpretation of CMC surfaces in terms of 
an infinite-dimensional integrable system was proposed in \cite{KonTai} and 
is based on the Weierstrass representation.

On a torus a holomorphic quadratic differential is constant (with respect to
a conformal parameter $z$). Given a CMC torus, by dilation of the surface and
a linear transform $z \to az$ of the conformal parameter, it is
achieved that
$$
A dz^2 = \frac{1}{2} dz^2, \ \ \ H =1.
$$
In this event the Gauss--Codazzi equations read
$$
u_{z\bar{z}} + \sinh u = 0
$$
where $u=2\alpha$ and $e^{2\alpha}dzd\bar{z}$ is the metric on the
torus. This equation is the compatibility condition for the
following system
\begin{equation}
 \left[ \frac{\partial}{\partial z} - \frac{1}{2} \left(
\begin{array}{cc}
-u_z & - \lambda \\
- \lambda & u_z
\end{array}
\right)\right]\psi = 0, \ \ \ \left[ \frac{\partial}{\partial
\bar{z}} - \frac{1}{2 \lambda} \left(
\begin{array}{cc}
0 & e^{-u} \\
e^u & 0
\end{array}
\right)\right] \psi = 0. \label{sinh-comm-2}
\end{equation}

Let $\Lambda$ be the period lattice for the torus.
We consider the
linear problem
$$
L \psi =
\partial_z \varphi -
\frac{1}{2} \left(
\begin{array}{cc}
-u_z & 0 \\
0 & u_z
\end{array}
\right) \psi = \frac{1}{2} \left(
\begin{array}{cc}
0 & - \lambda \\
- \lambda & 0
\end{array}
\right) \psi.
$$
Since $L$ is a first order $2\times 2$-matrix operator, for every
$\lambda \in \C$ the system (\ref{sinh-comm-2}) has a
two-dimensional space $V_\lambda$ of solutions and these spaces
are invariant under the translation operators
$$
\widehat{T}_j \varphi (z) = \varphi (z + \gamma_j), \ \ \ j=1,2,
$$
where $\gamma_1$ and $\gamma_2$ are generators of $\Lambda$. The
operators $\widehat{T}_1, \widehat{T}_2$, and $L$ commute and
therefore have common eigenvectors which are glued into a
meromorphic function $\psi (z,\bar{z},P)$ on a two-sheeted
covering
$$
\widehat{\Gamma} \rightarrow \C: P \in \widehat{\Gamma}
\rightarrow \lambda \in \C,
$$
ramified at points where $\widehat{T}_j$ and $L$ are not
diagonalized simultaneously. This the standard procedure for
constructing spectral curves of periodic operators \cite{Novikov}.

To each point $P \in \widehat{\Gamma}$ there corresponds a unique
(up to a constant multiple) Floquet function $\psi(z,\bar{z},P)$
with multipliers $\mu(\gamma_1,P)$ and $\mu(\gamma_2,P)$. The
complex curve $\widehat{\Gamma}$ is compactified by four
``infinities'' $\infty_{\pm}^1, \infty_{\pm}^2$ such that
$\infty_{\pm}^1$ are projected into $\lambda = \infty$ and
$\infty_{\pm}^2$ are projected into $\lambda = 0$ and we may take
$\psi$ meromorphic on $\widehat{\Gamma}$ with the following
essential singularities at the ``infinities'':
$$
\psi(z,\bar{z},P) \approx \exp{\left(\mp\frac{\lambda
z}{2}\right)} \left(
\begin{array}{c} 1 \\ \pm 1
\end{array}\right)
\ \ \mbox{as $P \to \infty_{\pm}^1$},
$$
$$
\psi(z,\bar{z},P) \approx \exp{\left(\mp
\frac{\bar{z}}{2\lambda}\right)} \left(
\begin{array}{c} 1 \\ \pm 1
\end{array}\right)
\ \ \mbox{as $P \to \infty_{\pm}^2$}.
$$
The multipliers tend to $\infty$ as $\lambda \to 0,\infty$.

The complex curve $\Gamma$ admits the involution which preserves
multipliers:
$$
\sigma(\lambda) = - \lambda, \ \ \
\left(
\begin{array}{c}
\varphi_1 \\ \varphi_2
\end{array}
\right)
\to
\left(
\begin{array}{c}
\varphi_1 \\ -\varphi_2
\end{array}
\right), \ \ \
\sigma(\infty^1_{\pm}) = \infty^1_{\mp}, \ \ \
\sigma(\infty^2_{\pm}) = \infty^2_{\mp}.
$$
The complex quotient curve $\widehat{\Gamma}/\sigma$ is called the
spectral curve of a CMC torus.

We have

\begin{proposition}
[\cite{T3}] $\varphi$ meets (\ref{sinh-comm-2}) if and only if
$\psi = (\lambda \varphi_2, e^{\alpha}\varphi_1)^{\top}$ satisfies
the Dirac equation $\D \psi = 0$ with $U = \frac{He^\alpha}{2} =
\frac{e^\alpha}{2}$.
\end{proposition}

Thus we have an analytic mapping of $\Gamma$ onto the  spectral
curve of a general torus defined in  \S \ref{subsec4.4} such that
the mapping preserves the multipliers. This implies these complex
curves coincide up to irreducible components. Together with
Corollary \ref{bound} this implies

\begin{proposition}
The spectral curve of a CMC torus in $\R^3$ coincides with the
(general) spectral curve, of the torus, defined in \S
\ref{subsec4.4}.
\end{proposition}

{\sc B) Minimal tori in $S^3$.} We consider the unit sphere in
$\R^4$ as the Lie group $SU(2)$. For minimal surfaces in $SU(2)$
the derivational equations (\ref{h1}) and (\ref{h2}) are
simplified and we obtain the Hitchin system \cite{Hitchin}
\begin{equation}
\bar{\partial} \Psi - \partial \Psi^\ast + [\Psi^\ast,\Psi] = 0, \
\ \
\bar{\partial} \Psi + \partial \Psi^\ast = 0.
\label{hitchin1}
\end{equation}
The first equation implies that the $SL_2$ connection ${\cal A}=
(\partial + \Psi, \bar{\partial} + \Psi^{\ast})$ on $f^{-1}(TG)$
is flat. In this event the second equation implies that this
connection is extended to an analytic family of flat connections
$$
{\cal A}_{\lambda} = \left(\partial +
\frac{1+\lambda^{-1}}{2}\Psi, \bar{\partial} +
\frac{1+\lambda}{2}\Psi^{\ast} \right)
$$
where ${\cal A} = {\cal A}_1$ and $\lambda \in \C\setminus \{0\}$.
Thus we obtain an ``L,A''-pair with a spectral parameter and
therefore derive that this system is integrable. This trick is
general for integrable harmonic maps.

Let us define the spectral curve.

Let $\Sigma$ be a minimal torus in $SU(2)$ and let
$\{\gamma_1,\gamma_2\}$ be a basis for $\Lambda$.We define
matrices $H(\lambda)$ and $\widetilde{H}(\lambda) \in SL(2,\C)$
which describe the monodromies of ${\cal A}_{\lambda}$ along
closed loops realizing $\gamma_1$ and $\gamma_2$ respectively.
These matrices commute and hence have joint eigenvectors
$\varphi(\lambda,\mu)$ where $\mu$ is a root of the characteristic
equation for $H(\lambda)$:
$$
\mu^2 - \Tr H(\lambda) + 1 = 0.
$$
The eigenvalues
$$
\mu_{1,2} = \frac{1}{2}\left( \Tr H(\lambda) \pm \sqrt{\Tr^2
H(\lambda) -4}\right)
$$
are defined on a Riemann surface $\Gamma$ which is a two-sheeted
covering of $\C P^1$ ramified at the odd zeros of the function
$(\Tr^2 H(\lambda) -4)$ and at $0$ and $\infty$ (multiple zeros
are removed by the normalization). This $\Gamma$ is the spectral
curve of a minimal torus in $SU(2)$ and has finite genus.

Above we expose Hitchin's results which are valid for all harmonic
tori in $S^3$ (this includes both cases of minimal tori in $S^3$ and
harmonic Gauss maps into $S^2 \subset S^3$) \cite{Hitchin}. Now we
have to confine to minimal tori in $S^3$.

Let $\D$ be the Dirac operator associated with this torus and the
spinor $\psi^\prime$ generates the torus via the Weierstrass
representation. Let us
$$
L =
\frac{1}{\sqrt{2}}
\left(
\begin{array}{cc}
\bar{a} & -\bar{b} \\ b & a
\end{array}
\right), \ \ \ a = -i\bar{\psi^\prime}_1 + \psi^\prime_2, \ \ \ b =
- i \psi^\prime_1 + \bar{\psi^\prime}_2.
$$

We have

\begin{proposition}
[\cite{T3}]
The Hitchin eigenfunctions $\varphi$
are transformed by the
mapping
$$
\varphi \to
\psi =
e^{\alpha} \left(
\begin{array}{cc}
0 & i \lambda \\
1 & 0
\end{array}
\right) \cdot L^{-1}
\varphi
$$
into solutions of the Dirac equation $\D\psi = 0$ corresponding to
the torus $\Sigma$ in $S^3$.
\end{proposition}

As in the case of CMC tori in $\R^3$ (see above) this Proposition
together with Corollary \ref{bound} implies

\begin{proposition}
The spectral curve of a minimal torus in $S^3$ coincides with the
(general) spectral curve of the torus (as this curve is defined in
\S \ref{subsec4.4}).
\end{proposition}

\subsection{Singular spectral curves}
\label{subsec4.7}

The perturbation of the free operator could be so strong that
another singularities (not coming from resonance pairs) could
appear in $\Gamma$. If $\Gamma_\nm$ is algebraic then we write
down the corresponding Baker--Akhiezer function
$\psi(z,\bar{z},P)$ such that

1) $\D \psi =0$;

2) $\psi$ is meromorphic on $\Gamma$ and has the following
asymptotics at the infinities:
$$
\psi \approx \left(
\begin{array}{c} e^{\lambda_+ z} \\ 0 \end{array}
\right) \ \ \mbox{as $P \to \infty_+$}, \ \ \ \psi \approx \left(
\begin{array}{c} 0 \\ e^{\lambda_- \bar{z}} \end{array}
\right) \ \ \mbox{as $P \to \infty_-$}
$$
where $\lambda^{-1}_\pm$ are local coordinates near $\infty_\pm$,
$\lambda^{-1}_\pm(\infty_\pm)=0$. We may put
$\lambda_\pm = 2\pi i k_1$.

The function $\psi$ is formed by Floquet functions
$\psi(z,\bar{z},P)$ taken at different points of the spectral
curve such that $\psi$ is meromorphic and has the asymptotics as
above. The function ``draws'' the complex curve
$\Gamma_\psi$ on which it is defined such that no one Floquet function
is counted twice in different points of $\Gamma_\psi$.
There is a chain of mappings
$$
\Gamma_\nm \to \Gamma_\psi \to \Gamma
$$
such that the composition of them is the normalization of $\Gamma$
and the first of them is the normalization of $\Gamma_\psi$. We
have evident inequalities:
$$
p_g(\Gamma) = p_g(\Gamma_\psi) \leq p_a(\Gamma_\psi) < \infty
$$
where $p_a(\Gamma_\psi)$ is the arithmetic genus of $\Gamma_\psi$
which differs from the geometric genus of $\Gamma_\psi$ by the
contribution of singular points.

The function $\psi$ can be pulled back onto a nonsingular curve
$\Gamma$ where it will have exactly $p_a(\Gamma_\psi)+1$ poles
(this follows from the finite gap integration theory).
For the Dirac operator, $p_a(\Gamma_\psi)$ equals ``the number of
poles of its normalized Baker--Akhiezer function'' minus one.

We arrive to the following conclusion:

\begin{itemize}
\item
the Baker--Akhiezer function $\psi$ defines the Riemann surface
$\Gamma_\psi$ in the classical spirit of the Riemann paper as the
surface on which the given function $\psi$ is naturally defined.
This surface is obtained from $\Gamma$ by performing
normalizations of singularities only when the dimension of the
space of Floquet functions at the point is decreased after the
normalization (for instance, this is the case of resonance pairs).

\item
in difference with $\Gamma_\nm$, the complex curve $\Gamma_\psi$ gives a one-to-one
parameterization of all Floquet functions (up to multiples).
\end{itemize}

For minimal tori in $S^3$ this situation is explained in detail in \cite{Hitchin}.

If we would like to construct a torus with finite spectral genus
in terms of theta functions we have to work with the curve
$\Gamma_\psi$ again as we demonstrated that in \S \ref{subsec4.5}
for the Clifford torus.

The following definition of $\Gamma_\psi$ comes from the finite
gap integration theory:

\begin{itemize}
\item
Let $\D$ be a Dirac operator with double-periodic potentials $U$
and $V$ and let $\Gamma_\psi$ be a Riemann surface (probably
singular) of finite arithmetic genus $p_a(\Gamma_\psi) =g$ with
two marked nonsingular points $\infty_\pm$ and local parameters $k_\pm^{-1}$
near these points such that $k^{-1}_\pm(\infty_\pm)=0$.

Let $\psi(z,\bar{z},P)$ be  a Baker--Akhiezer function $\psi$
which is defined on $\C \times \Gamma_\psi
\setminus\{\infty_\pm\}$ such that

1) $\psi$ is meromorphic in $P$ outside $\infty_\pm \in \Gamma$
and has poles at $g+1$ nonsingular points $P_1+ \dots + P_{g+1}$;

2) $\psi$ has the following asymptotics at $\infty_\pm$:
$$
\psi \approx e^{k_+ z} \left[ \left(
\begin{array}{c}
1 \\ 0
\end{array}
\right) + \left(
\begin{array}{c}
\xi^+_1 \\ \xi^+_2
\end{array}
\right) k_+^{-1} + O(k_+^{-2}) \right] \ \ \mbox{as $P \to
\infty_+$},
$$
$$
\psi \approx e^{k_- \bar{z}} \left[ \left(
\begin{array}{c}
0 \\ 1
\end{array}
\right) + \left(
\begin{array}{c}
\xi^-_1 \\ \xi^-_2
\end{array}
\right) k_-^{-1} + O(k_-^{-2}) \right] \ \ \mbox{as $P \to
\infty_-$}
$$
and $\psi$ satisfies the Dirac equation $\D \psi=0$ everywhere on
$\Gamma_\psi$ except the ``infinities'' $\infty_\pm$ and the poles
of $\psi$.

We say that $\Gamma_\psi$ is the {\it spectral curve of a finite
gap operator} $\D$.
\end{itemize}

For a generic divisor $P_1+\dots+P_{g+1}$ such a function is
unique and the potentials are reconstructed from it by the
formulas:
\beq
\label{reconstruction}
 U = -\xi^+_2, \ \ V =\xi^-_1.
\eeq

The attempt to define such a Riemann surface in the case when $p_g(\Gamma) =
\infty$ meets a lot of analytical difficulties.

We refer to \cite{T4} for more detailed exposition of some questions related to
singular spectral curves.

We see in \S \ref{subsec4.5} that for the Clifford torus
$\Sigma_{1,1} \subset \R^4$ the potential is constant and the
spectral curve is a sphere. Moreover
$$
p_g(\Gamma) = p_a(\Gamma_\psi) = 0.
$$
However the potential of its stereographic projection, which is
the Clifford torus in $\R^3$, equals
$$
U = \frac{\sin x}{2\sqrt{2}(\sin x  - \sqrt{2})}
$$
where $x$ is one of the angle variables and, by Theorem
\ref{cliftheorem}, for an operator with this potential we have
$$
p_g(\Gamma) = 0, \ \ \ p_a(\Gamma_\psi) = 2.
$$
Therefore the stereographic projection of the Clifford torus from $S^3$ into
$\R^3$ results in the appearance of singularities in
$\Gamma_\psi$.

This leads to an interesting problem:

\begin{itemize}
\item
{\sl what is the relation between the spectral curve of a torus in
the unit sphere $S^3 \subset \R^4$ and the spectral curve of its
stereographic projection?}
\end{itemize}

We think that the answer to this question is as follows: the
potentials are related by some B\"acklund transformation which
leads to a transformation of the spectral curve. Probably there is
an analogy with such a transformation for the Schr\"odinger
operator exposed in \cite{EK}. We also expect that the answer to
the following question is positive:

\begin{itemize}
\item
{\sl do the images ${\cal M}(\Gamma)$ of the multiplier mappings for
a torus in $S^3$
and for its stereographic projection coincide?}
\end{itemize}

There is another interesting problem:

\begin{itemize}
\item
{\sl characterize the spectral curves of tori in $\R^3$ and $\R^4$.}
\end{itemize}

For tori in $\R^3$ and in $\R^4$ the answers have to be different.
Indeed, it was already mentioned in  \cite{Bobenko} that the
spectral curves for CMC tori in $\R^3$ have to be singular (for
them that results in the appearance of multiple branch points
which are transformed by the normalization into pairs of points
interchanged by the hyperelliptic involution). \footnote{In
\cite{Schmidt} it is showed that for tori in $\R^3$ ${\cal
M}(\Gamma)$ contains a point of multiplicity at least four or a
pair of double points at which the differentials $dk_1$ and $dk_2$
vanish (here $k_1$ and $k_2$ are quasimomenta). We notice that
does not mean that $\Gamma_\psi$ meets the same conditions: for
instance, for the Clifford torus in $\R^3$ the spectral curve
$\Gamma_\psi$ has a pair of double points at which $dk_1$ and
$dk_2$ do not vanish and has no more singular points.} However for
the Clifford torus in $\R^4$ the spectral curve is nonsingular.

\section{The Willmore functional}

\subsection{Willmore surfaces and the Willmore conjecture}
\label{subsec5.1}

The Willmore functional for closed surfaces in $\R^3$ is defined
as \beq \label{willmore} \W(\Sigma) = \int_\Sigma H^2 d\mu \eeq
where $d\mu$ is the induced area form on the surface. It was
introduced by Willmore in the context of variational problems
\cite{Willmore}. Therewith Willmore was first who stated a global
problem of the conformal geometry of surfaces, i.e. the Willmore
conjecture which we discuss later. The Euler--Lagrange equation
for this functional takes the form
$$
\Delta H+2H(H^2-K)=0
$$
where $\Delta$ is the Laplace--Beltrami operator on the surface.
Surfaces meeting this equation are called Willmore surfaces.

We remark that $H = \frac{\varkappa_1+\varkappa_2}{2}$ and, by the
Gauss--Bonnet theorem, for a compact oriented surface $\Sigma$
without boundary we have
$$
\int_\Sigma K d\mu = \int_\Sigma \varkappa_1 \varkappa_2 d\mu =
2\pi \chi(\Sigma)
$$
where $\chi(\Sigma)$  is the Euler characteristic of $\Sigma$. By
adding a topological term to $\W$ we obtain the functional with
the same extremals among closed surfaces and may simplify the
variational problem. For spheres it takes the case when
considering the functional
$$
\widehat{\W} (\Sigma) = \int (H^2-K) d\mu = \W(\Sigma) - 2\pi
\chi(\Sigma)
$$
we conclude that
\footnote{This is similar to the instanton trick
which led to a discovery of selfdual connections.}
$$
\widehat{\W} = \frac{1}{4} \int_\Sigma (\varkappa_1-\varkappa_2)^2
d\mu.
$$
We recall that a point on a surfaces is called an {\it umbilic}
point if $\varkappa_1 = \varkappa_2$ at it. A surface is called
{\it totally umbilic} if all its points are umbilics. By the Hopf
theorem, a totally umbilic surface in $\R^3$ is a domain in a
round sphere or in a plane. For spheres this gives for spheres a
lower estimate for the Willmore functional and a description of
all its minima:

\bi
\item
for spheres
$$
\W(\Sigma) \geq 4\pi $$
and
$\W(\Sigma) = 4\pi$ if and only if $\Sigma$ is a round sphere.
\ei

For surfaces of higher genus this trick does not work.

The functional $\widehat{\W}$ was introduced by to Thomsen
\cite{Th} and Blaschke \cite{Blaschke} who called it the conformal
area for the following reasons:

1) the quantity $(H^2 - K)d\mu$ is invariant with respect to
conformal transformations of the ambient space and therefore,
given a compact oriented surface $\Sigma \subset \R^3$ and a
conformal transformation $G: \overline{\R}^3 \to \overline{\R}^3$
which maps $\Sigma$ into a compact surface, we have
$$
\widehat{\W} (\Sigma) = \widehat{\W}(G(\Sigma));
$$

2) if $\Sigma$ is a minimal surface in $S^3$ and $\pi:S^3 \to
\overline{\R}^3$ is the stereographic projection which maps
$\Sigma$ into $\R^3$, then $\pi(\Sigma)$ is a Willmore surface.

Moreover as it is proved in \cite{Bryant1}

3) outside umbilic points there is defined a quartic differential
$\widehat{A}(dz)^4$ which is holomorphic if the surface is a
Willmore surface;

We expose these results in Appendix 2.

By 2) there are examples of compact closed Willmore surfaces. We
remark that not all compact Willmore surfaces are stereographic
projections of minimal surfaces in $S^3$ (this was first showed
for tori in \cite{Pinkall}).

It follows from 3) that outside umbilics Willmore surfaces admit a
good description which is similar to the description of CMC
surfaces in terms of the holomorphicity of a quadratic Hopf
differential. However there are examples of compact Willmore
surfaces which even have lines consisting of umbilics \cite{BB}.

By 1) the minimum of the Willmore functional in each topological
class of surfaces  is conformally invariant and hence degenerate.
We note that the existence of a minimum which real-analytical
surface was proved for tori by Simon \cite{Simon} and for surfaces
of genus $g \geq 2$ by Bauer and Kuwert \cite{Bauer}. Recently
Schmidt presented the proof of the following result: given a genus
and a conformal class of an oriented surface, the Willmore
functional achieves its minimum on some surface which a priori may
have branch points or be a branched covering of an immersed
surface \cite{Schmidt2}. His technique uses the Weierstrass
representation and some ideas from \cite{Schmidt}. \footnote{See
Appendix 1.}

Bryant started the program of describing all Willmore spheres by
applying the fact that at a holomorphic  quartic differential on a
sphere vanishes and, therefore, Willmore spheres admit description
in terms of algebro-geometric data \cite{Bryant1}. We have

\bi
\item
the image of a minimal surface in $\R^3$ under a Moebius transform
$(x-x_0) \to (x-x_0)/|x-x_0|^2$ is a Willmore surface and any
minimal surface $\Sigma$ with planar ends is mapped by a Moebius
transform. with the center $x_0$ outside the surface, into a
smooth compact Willmore surface $\Sigma^\prime$ such that
$$
\W(\Sigma^\prime) = 4\pi n,
$$
where $n$ is the number of planar ends of $\Sigma$.
\ei

\noindent
Bryant proved that all Willmore spheres are Moebius
inverses of minimal surfaces with planar ends, that the case $n=1$
corresponds to the round spheres, there are no such spheres with
$n=2$ and $3$, and described all Willmore spheres with $n=4$.
Later it was proved in \cite{Peng} that Willmore  spheres exist
for all even $n \geq 6$ and all odd $n \geq 9$. The left cases
$n=5$ and $7$ were finally excluded in \cite{Bryant2}.

The Willmore conjecture states that

\bi
\item
{\sl for tori
$$
\W \geq 2\pi^2
$$
and the Willmore functional attains its minimum on the Clifford torus and
its images under conformal transformations of $\overline{\R}^3$.}
\ei

The Clifford torus was already introduced in \S
\ref{subsec4.5}.

Since the Willmore functional is conformally invariant and the
stereographic projection $\pi:S^3 \to \overline{R}^3$ is
conformal, we do not distinguish the original Willmore conjecture
and its counterpart for tori in $S^3$ for which the Willmore
functional is replaced by \beq \label{wills3} \W_{S^3} = \int (H^2
+1) d\mu, \ \ \ \ \W_{S^3} (\Sigma) = \W(\pi(\Sigma)). \eeq

Willmore introduced his conjecture in \cite{Willmore} where he checked it for round tori of revolution.

It is proved in many special cases:

1) for tube tori, i.e for tori
formed by carrying a circle centered at a closed curve along this curve
such that the circle always lies in the normal plane,
by Shiohama and Takagi \cite{Shiohama} and by Willmore \cite{Willmore2}
(if we admit for the radius of the circle to vary we obtain channel tori
for which the conjecture was established in \cite{HP});

2) for tori of revolution by Langer and Singer \cite{LangerSinger};

3) for tori conformally equivalent to $\R^2/\Gamma(a,b)$ with
$0\le a\le1/2$, $\sqrt{1-a^2}\le b\le 1$
where the lattice $\Gamma(a,b)$ is generated by $(1,0)$ and $(a,b)$
(Li--Yau \cite{LiYau});

4) the previous result of Li and Yau was improved by Montiel and Ros who extended it to the case
$\left(a - \frac{1}{2}\right)^2 + (b-1)^2 \leq \frac{1}{4}$ \cite{MontielRos};

5) for tori in $S^3$ which are symmetric under the antipodal
mapping (Ros \cite{Ros}).

6) since it was also proved by Li and Yau that if a surface has a selfintersection point with multiplicity $n$
then $\W \geq 4\pi n$ \cite{LiYau}, the conjecture is proved for tori with selfintersections.

Some other partial results were obtained in \cite{Ammann,Topping}.

We also mention the paper \cite{Weiner} where the second variation
form of $\W$ for the Clifford torus was computed and it was proved
that this form is non-negative. The second variation formula for
general Willmore surfaces was obtained in \cite{Palmer}.

In general case the conjecture stays open.

We shall discuss some recent approach applied in \cite{Schmidt} in the next paragraph.

By (\ref{wills3}), the following conjecture is a special case, of the Willmore conjecture,
which also stays open:

\bi
\item
{\sl for minimal tori in $S^3$ the volume is bounded from below by
$2\pi^2$ and attains its minimum on the Clifford torus in $S^3$.}
\ei

By the Li--Yau theorem on surfaces with selfintersections this
conjecture follows from the following conjecture by Hsiang and
Lawson:

\bi
\item
{\sl the Clifford torus is the only minimal torus embedded in $S^3$.}
\ei

Since a holomorphic quartic differential on a torus has constant
coefficients, there are two opportunities: it vanishes or it
equals $c(dz)^4, c = {\mathrm const} \neq 0$.

In the first case a torus is obtained as a Moebius image of a
minimal torus with planar ends. For evident reasons it is clear
that there are no such tori with $n=1$ and $2$ ends. The case
$n=3$ was excluded by Kusner and Schmitt who also constructed
examples with $n=4$ \cite{KusnerS}. First examples of minimal
rectangular tori with four planar ends were constructed by Costa
\cite{Costa}. Recently Shamaev constructed such tori for all even
$n \geq 6$ \cite{Shamaev}. However it looks from the construction
that in generic case these tori do not have branch points that was
rigorously proved only for $n=6,8$, and $10$.

In the second case the Codazzi type equations for Willmore tori
without umbilics coincide with the four-particle Toda lattice
\cite{FPPS,Babich}. \footnote{See Appendix 2.} The theta formulas
for such Willmore tori are derived in \cite{Babich} by using Baker--Akhiezer
functions related to this Toda lattice.

Another construction of Willmore tori by methods of integrable
systems was proposed in \cite{Helein0}.

For surfaces of higher genus the candidates for the minima of the
Willmore functional were proposed by Kusner \cite{Kusner}.

There is the conjecture that for tori in $\R^4$ the Willmore
functional $\int |H|^2 d\mu$ attains its minimum on the Clifford
torus in $\R^4$, i.e. the product of two circles of the same radii
(see \cite{Wintgen1,Wintgen2}) . Since this torus is Lagrangian
the last conjecture is weakened by assuming that the Clifford
torus is the minimum for $\W$ in a smaller class of Lagrangian
tori. That is discussed in \cite{Minicozzi} where it is proved
that $\W$ achieves its minimum among Lagrangian tori on some
really-analytical torus.

We do not discuss the generalization of the Willmore functional
for surfaces in arbitrary Riemannian manifolds which is
$$
\int (|H|^2 + \widehat{K}) d\mu
$$
where $\widehat{K}$ is the sectional curvature of the ambient
space along the tangent plane to the surface. The quantity $(|H|^2
- K + \widehat{K})d\mu$ is invariant with respect to conformal
transformations of the ambient space \cite{Chen}.

In \cite{Berdinsky} another generalization of the Willmore
functional for surfaces in three-dimensional Lie groups is
proposed. It is based on the spectral theory of Dirac operators
coming into Weierstrass representations (see also \S
\ref{subsec5.5}).

We also have to mention the Willmore flow which is similar to the
mean curvature flow and decreases the value of $\W$ (see the paper
\cite{KS} and references therein).

We finish this part by a remark on constrained Willmore surfaces
which are, by definition,  critical points of the Willmore
functional restricted onto the space of surfaces with the same
conformal type. It was first observed by Langer that compact
constant mean curvature surfaces in $\R^3$ are constrained
Willmore since for them the Gauss map is harmonic \cite{PS0}. We
refer for the basics of the theory of such surfaces to
\cite{BPPin}.

\subsection{Spectral curves and the Willmore conjecture}
\label{subsec5.2}

As it is showed in \cite{T1} in terms of the potential $U$ of a Weierstrass representation of a torus
in $\R^3$ the Willmore functional
is
$$
\W = 4 \int_M U^2 dx dy.
$$
So it measures the perturbation of the free operator.

We recall that the Willmore conjecture states that this functional
for tori attains its minimum at the Clifford torus for which the
Willmore functional equals $2\pi^2$.

Starting from the observation that the Willmore functional is the
first integral of the mNV flow which deforms tori into tori
preserving the conformal class (see \S \ref{subsec3.1}) we
introduced in 1995 the following conjecture (see \cite{T1}):
\begin{itemize}
\item
{\sl a nonstationary (with respect to the mNV flow) torus cannot
be a local minimum of the Willmore functional.}
\end{itemize}

It was based on the assumption that a minimum of such a
variational problem is nondegenerate and thus has 
to be stable with respect to soliton
deformations which are 
governed by equations from the mNV hierarchy and which preserve the value of 
the Willmore functional. By soliton theory these equations are linearized on 
the Jacobi variety of the normalized spectral curve and generically these 
linear flows span this Jacobi variety which is an Abelian variety of 
complex dimension $p_g(\Gamma)$ or some Prym subvariety of the Jacobi
variety.

Its geometrical analog was formulated in \cite{T2} where we
introduced a notion of the spectral genus of a torus as
$p_g(\Gamma)$:
\begin{itemize}
\item
{\sl given a conformal class of tori in $\R^3$, the minima of the
Willmore functional are attained at tori of the minimal spectral
genus.}
\end{itemize}

In \cite{T2} we proposed the following explanation to the lower
bounds for $\W$: for small perturbations of the zero potential
$U=0$ the Weierstrass representation gives us planes which do not
convert into tori and, since for surfaces in $\R^3$ the Willmore
functional is the squared $L_2$-norm of $U$, the lower bound shows
how large a perturbation of the zero potential has to be to force
the planes to convert into tori.

The strategy to prove the Willmore conjecture after proving the
last conjecture is to calculate the values of the Willmore
functional for tori of the minimal spectral genus (by using the
formula  (\ref{willmoreformula}) or by other means) and to check
the Willmore conjecture.

We already mentioned in this text the paper \cite{Schmidt} by
Schmidt. This paper contains a series of
important results. \footnote{In this paper it is also presented a
proof of the fact that the spectral genus of a constrained Willmore
torus in $\R^3$ is finite. Another proof was presented by
Krichever (unpublished). This fact is nontrivial even for Willmore
tori because the trick which uses soliton theory and was applied
in \cite{Hitchin,PS} to harmonic tori in $S^3$ and constant mean
curvature tori in $\R^3$ (see also \cite{Bob1}) does work not for 
all tori. It is applied only to
tori described by the four-particle Toda lattice without umbilic
points at which $e^\beta = 0$ (see Appendix 2).} For our interests
we expose only the results related to the asymptotic behavior of
the spectral curve. Although until recently we did not go through
all details of \cite{Schmidt} we have to say that

{\it in fact the paper \cite{Schmidt} proposes a proof only for
our conjecture (see above). The value of $p_a(\Gamma_\psi)$ is a
priori unbounded however in \cite{Schmidt} the calculations of
values of the Willmore functional are done only for the minimal
possible values of both $p_g(\Gamma)$ and $p_a(\Gamma_\psi)$.}

Following the spirit of the previous conjectures it is natural to guess
that

\begin{itemize}
\item
{\sl given a conformal class of tori in $\R^3$ and a spectral
genus, the minima of the Willmore functional are attained at tori
with the minimal value of $p_a(\Gamma_\psi)$.}
\end{itemize}

This conjecture also fits in the soliton approach since 
the additional dergees of freedom coming from $p_a(\Gamma)-p_g(\Gamma)$
(or some part of it, is the flows are linearized on the Prym variety)
also correspond to soliton deformations.

By our opinion these conjectures are interesting by their own
means. We remark that proofs of the last two of them together with
calculations of the values of the Willmore functional for tori
with minimal possible values of $p_g(\Gamma)$ and
$p_a(\Gamma_\psi)$ would lead to checking the Willmore conjecture.

We would like also to mention another interesting problem:

\begin{itemize}
\item
{\sl how to generalize this spectral curve theory 
for compact immersed surfaces of
higher genera?}
\end{itemize}

\subsection{On lower bounds for the Willmore functional}
\label{subsec5.3}

In \cite{T21} we established in some special case a lower estimate,
for the Willmore functional, which is quadratic in the dimension of the
kernel of the Dirac operator.

Let represent the sphere as a infinite cylinder $Z$ compactified
by a couple of points such that $z = x+iy$ is a conformal
parameter on $Z$, $y$ is defined modulo $2\pi$, $x \in \R$, and
these two ``infinities'' are achieved as $x \to \pm \infty$.

\begin{lemma}[\cite{T21}]
\label{lemmasphere} Given a sphere in $\R^3$, the asymptotics of
$\psi$ and the potential $U$ are as follows
$$
|\psi_1|^2 + |\psi_2|^2 = C_\pm e^{-|x|} + O(e^{-2|x|}), \ \ \ U =
U_\pm e^{-|x|} + O(e^{-2|x|}) \ \ \ \mbox{as $x \to \pm \infty$},
$$
where $C_\pm$ and $U_\pm$ are constants. If $C_+ = 0$ or $C_- =0$
then there is a branch point at the corresponding marked point $x
= +\infty$ or $x = -\infty$ respectively.

The kernel of $\D$ on the sphere consists of solutions $\psi$ to
the equation $\D \psi =0$ on the cylinder such that $|\psi_1|^2 +
|\psi_2|^2 = O(e^{-|x|})$ as $x \to \pm \infty$.
\end{lemma}

Let us assume that the potential $U($ of the Dirac operator depends on $x$ only.
For instance, such a situation realizes for a sphere of revolution for which
$y$ is an angle of rotation.
However this is not only the case
of spheres of revolution and there are more such surfaces with an
intrinsic $S^1$-symmetry reflected by the potential of the
Weierstrass representation.

\begin{theorem}[\cite{T21}]
\label{willmorenumber}
Let $\D$ be a Dirac operator on $M = S^2$
with a real-valued potential $U=V$ which depends only on the
variable $x$. Then \beq \label{quadratic} \int_M U^2 dx \wedge dy
\geq \pi N^2 \eeq where $N = \dim_\H \Ker \D = \frac{1}{2}
\dim_\C\Ker \D$.
These minimuma are achieved on the potentials
$$
U_N(x) = \frac{N}{2\cosh x}.
$$
\end{theorem}

The proof of this theorem is based on the method of the inverse
scattering problem applied to a one-dimensional Dirac operator.
This quadratic estimate appeared from the trace formulas by
Faddeev and Takhtadzhyan \cite{Faddeev}.

Before giving the proof of Theorem \ref{willmorenumber} we 
expose one of its consequences, i.e.
Theorem \ref{theofriedrich}.

Together with the proof of Theorem \ref{willmorenumber} we 
introduced the following conjecture.

\begin{conjecture}[\cite{T21}]
\label{taimconj} The estimate (\ref{quadratic}) is valid for all
Dirac operators on a two-sphere.
\end{conjecture}

Very soon after the electronic publication of \cite{T21} appeared
Friedrich demonstrated the following corollary of this conjecture:
\footnote{The conjecture was finally proved by Ferus, Leschke, Pedit, and
Pinkall in \cite{FLPP} together with the generalization of
(\ref{quadratic}), the so-called Pl\"ucker formula, for surfaces
of higher genera (we expose that in \S
\ref{subsec5.4}).}

\begin{theorem}
[\cite{F2}] \label{theofriedrich}
Let us assume that Conjecture
\ref{taimconj} holds.
Given an eigenvalue $\lambda$ of the Dirac
operator on a two-dimensional spin-manifold homeomorphic to the
two-sphere $S^2$, the inequality holds \beq \label{estfriedrich}
\lambda^2 \mathrm{Area}(M) \geq \pi m^2(\lambda) \eeq where
$m(\lambda)$ is the multiplicity of $\lambda$.
\end{theorem}

We remark that due to the symmetry (\ref{ast}) of $\Ker D$ the
multiplicity of an eigenvalue is always even. For the case
$m(\lambda)=2$ the inequality (\ref{estfriedrich}) was proved by
B\"ar \cite{Baer}.

{\sc Proof of Theorem \ref{theofriedrich}.} First we recall the
definition of the Dirac operator on a spin-manifold (see
\cite{F3,LM} for detailed expositions).

A spin $n$-manifold $M$ is a Riemannian manifold with a spin
bundle $E$ over $M$ such that at each point $p \in M$ there is
defined a Clifford multiplication
$$
T_p M \times E_p \to E_p
$$
such that
$$
v\cdot w \cdot \varphi + w \cdot v \cdot \varphi = -
2(v,w)\varphi, \ \ \ v,w \in T_p M, \ \psi \in E_p.
$$
We also assume that there is a Riemannian connection $\nabla$
which induces a connection on $E$. Then the Dirac
operator is defined at every point $p \ M$ as
$$
D \varphi = \sum_{k=1}^n e_k \cdot \nabla_{e_k} \varphi
$$
where $e_1,\dots,e_n$ is an orthonormal basis for $T_p M$ and
$\varphi$ is a section of $E$.

We consider as an example a two-dimensional spin manifold $M$ with
a flat metric. The Clifford algebra ${\mathcal Cl}_2$ is
isomorphic to $\H$. Thus we have a $\C^2$-spin bundle over $M$
(here we identify $\H$ with $\C \oplus \C$). For the flat metric
on $M$ the Clifford multiplication is represented by the matrices
$$
e_1 = e_x = \left(\begin{array}{cc} 0 & 1 \\ -1 & 0
\end{array}\right), \ \ \ \ e_2 = e_y = \left(\begin{array}{cc} 0 & -i
\\ -i & 0
\end{array}\right).
$$
It is easy to check that
$$
e_x e_y + e_y e_x = 0, \ \ \ e_x^2 = e_y^2 = - 1.
$$
The Dirac operator $D_0$ is given by the formula
$$
D_0 = e_x\cdot  \partial_x + e_y \cdot \partial_y =
2 \left(\begin{array}{cc} 0 & \partial \\
-\bar{\partial} & o
\end{array}\right) = 2 \D_0
$$
and its square equals to the Laplace operator (up to a sign):
$$
D^2_0 = -\partial^2_x - \partial^2_y.
$$

For a conformally Euclidean metric $e^\sigma dz d\bar{z}$ the
Dirac operator takes the form
$$
D = e^{-3\sigma/4} D_0 e^{\sigma/4}
$$
(see \cite{BFGK}). Therefore the eigenvalue problem
$$
D \varphi = \lambda \varphi
$$
for the Dirac operator associated with such a metric takes the
form
$$
D_0[e^{\sigma/4} \varphi] - \lambda e^{\sigma/2}
[e^{\sigma/4}\varphi] = 0
$$
which we rewrite as
$$
(\D_0 + U)\psi = 0, \ \ \ \ \ U = -\frac{\lambda e^{\sigma/2}}{2},
\ \ \psi = e^{\sigma/4}\varphi.
$$

If Conjecture \ref{taimconj} holds we have the inequality
$$
\int_M U^2 dx \wedge dy = \frac{\lambda^2}{4} \mathrm{Area} (M)
\geq \pi \left( \frac{\dim_\C \Ker (D_0 + U)}{2}\right)^2 = \pi
\frac{m^2(\lambda)}{4}.
$$
This proves Theorem \ref{theofriedrich}.

{\sc Proof of Theorem \ref{willmorenumber}.}
If the potential $U$ depends only on $x$
the linear space of solutions to $\D \psi = 0$ on the sphere
$S^2 = Z \cup \pm \infty = \R_x \times S^1_y \cup \infty$
is spanned by the functions of the form
$\psi(x,y) = \varphi(x)e^{\varkappa y}$ such that
$$
L \varphi :=
\left[\left(
\begin{array}{cc}
0 & \partial_x \\
- \partial_x & 0
\end{array}
\right) + \left(
\begin{array}{cc}
2U  & 0 \\
0 & 2U
\end{array}\right)\right]
\varphi
=
\left(
\begin{array}{cc}
0 & i\varkappa  \\
i\varkappa & 0
\end{array}
\right) \varphi
$$
where $e^{2\pi \varkappa} = - 1$ (this condition defines the spin
bundle over the sphere, see \cite{T12}) and $\varphi$ is
exponentially decaying as $x \to \pm \infty$.  This means that
$\varphi$ is the bounded state of $L$, i.e. $\varkappa$ belongs to
the discrete spectrum which is invariant with respect to the
complex conjugation $\varkappa \to \bar{\varkappa}$. Therefore
$\dim_\C \D=2N$ is twice the number of bounded states meeting the
condition $\Im \varkappa > 0$.

The trace formula (\ref{zs8}) (see Appendix 3) for for $p=q=2U$
takes the form
$$
\int^\infty_{-\infty} U^2(x)dx = -\frac{1}{4\pi}
\int^\infty_{-\infty} \log (1-|b(k)|^2)dk + \sum_{j=1}^N
\Im \varkappa_j.
$$
Given $\dim \Ker \D = N$, the functional $\int_M U^2(x) dx \wedge dy =
2\pi \int^\infty_{-\infty} U^2(x) dx$
achieves its minimum on the potential with the following spectral data:
$$
b(k) \equiv 0, \ \ \varkappa_k = \frac{i(2k-1)}{2}, \ k=1,\dots,N,
$$
and we have
$$
\int_{S^2} U^2(x) dx \wedge dy \geq 2\pi \left(\frac{1}{2} +
\frac{3}{2} + \dots + \frac{N}{2} \right) = \pi N^2.
$$
Actually there is $N$-dimensional family of potentials with such the
spectral data and it is parameterized by $\lambda_1,\dots,\lambda_N$
and moreover this family is invariant with respect to the mKdV
equations. It is easy to show that every such a family contains the
potential $U_N = \frac{N}{2\cosh x}$; $p_N(x) = 2U_N(x) = N/\cosh x$
is the famous $N$-soliton potential of the Dirac operator.

Theorem \ref{willmorenumber} is proved.

We see that the equality in (\ref{quadratic}) is achieved on some
special spheres which are particular cases of the so-called {\it
soliton spheres} introduced in \cite{T21}. By definition, these
are spheres for which potential of the Dirac operator $\D$ is a
soliton (reflectionless) potential $U(x)$. It is also worth to
select a special subclass of soliton spheres distinguished by the
condition that all poles $\varkappa_1,\dots,\varkappa_N, \Im
\varkappa_k > 0$, of the transition coefficient $T(k)$ are of the
form $\frac{(2m+1)i}{2}$, $m \in N$.

Soliton spheres are easily constructed from the spectral data  via
the inverse scattering method (see (\ref{zs9}) in Appendix 3).

We showed in \cite{T21} that

\bi
\item
the lower estimate (\ref{quadratic}) achieves the equality on the
soliton spheres corresponding to the potentials $U_N =
\frac{N}{2\cosh x}$);

\item
generically a soliton sphere is not a surface of revolution. 
\footnote{Indeed, let us denote by $f_1 = \varphi_1(x)e^{\varkappa_1
y},\dots, f_n = \varphi_N e^{\varkappa_N y}$ the distinct
generators of $\Ker \D$. Then any linear combination $f = \alpha_1
f_1 + \dots + \alpha_N f_N$ via the Weierstrass representation
gives rise to a sphere in $\R^3$. If there is a pair of
nonvanishing coefficients $\alpha_j$ and $\alpha_k$ such that $\Im
\varkappa_j \neq \Im \varkappa_k$ then this sphere is not a
surface of revolution.}

\item
the class of soliton spheres is preserved by the mKdV deformations
(note that they are defined by $1+1$-equations) for which the
Kruskal--Miura integrals are integrals of the motion;

\item
soliton spheres corresponding to the potentials $U_N =
\frac{N}{2\cosh x}$ are described in terms of rational functions,
\footnote{From the reconstruction formulas
(\ref{zs9}) it is clear that that holds for all reflectionless
potentials.} i.e.
they can be called rational spheres;

\item
soliton spheres such that each pole $\varkappa_j$ is of the form
$\frac{(2m_j+1)i}{2}$ are critical points of the Willmore
functional restricted onto the class of spheres with
one-dimensional potentials.
\ei

\subsection{The Pl\"ucker formula}
\label{subsec5.4}

Our attempts to prove Conjecture \ref{taimconj} had failed due to
the lack of well-developed inverse scattering method for
two-dimensional operators. However in a fabulous paper \cite{FLPP}
this conjecture was proved together with its generalization for 
surfaces of arbitrary genera
by using methods of algebraic geometry.

As it was mentioned in \cite{PP} the following statement can be
derived from the results of \cite{Aron} (see also \cite{HW}):

\begin{proposition}
Let $E$ be a $C^2$-bundle over a surface $M$ $\psi$ and let $\psi$
be a nontrivial section of $E$ such that $\D \psi =0$. Then the
zeroes of $\psi$ are isolated and for any local complex coordinate
$z$ on $M$ centered at some zero $p$ of the function $\psi$: 
$z(p)=0$, we
have
$$
\psi = z^k \varphi + O(|z|^{k+1})
$$
where $\varphi$ is a local section of $E$ which does not vanish in
a neighborhood of $p$. The integer $k$ is well-defined independent
of choice $z$.
\end{proposition}

This integer number $k$ is called the order of the zero $p$:
$$
\ord_p \psi = k.
$$

Now recall the Gauss equation (see Proposition \ref{prop1} in \S \ref{subsec2.1}):
\beq
\label{gaussquat}
\alpha_{z\bar{z}} + U^2 - |A|^2
e^{-2\alpha} = 0, \ \ \ e^\alpha = |\psi_1|^2+|\psi_2|^2.
\eeq

For simplicity, we assume that $M$ is a sphere and $E$ is a spin
bundle. If $\psi$ vanishes nowhere then it defines a surface in
$\R^3$ and integrating the left-hand side of (\ref{gaussquat})
over $M$ we obtain \beq \label{gaussquat1} \int_M
\alpha_{z\bar{z}} dx \wedge dy + \int_M U^2 dx \wedge dy - \int_M
|A|^2 e^{-2\alpha} dx \wedge dy = 0. \eeq By the Gauss theorem,
the first term equals
$$
-\frac{1}{4} \int_M \left(-4\alpha_{z\bar{z}}e^{-2\alpha}\right)
e^{2\alpha} dx \wedge dy  = -\frac{1}{4} \int_M K d\mu = -\pi,
$$where $K$ is the Gaussian curvature and $d\mu$ is the measure
corresponding to the induced metric. Thus we have \footnote{For
general complex quaternionic line bundles $L =E_0 \oplus E_0$ we
have $\int_M \alpha_{z\bar{z}} dx \wedge dy = \pi \deg E_0 = \pi
d$.}
$$
\int_M U^2 dx \wedge dy = \pi + \int_M |A|^2 e^{-2\alpha} dx
\wedge dy \geq -\int_M \alpha_{z\bar{z}} dx \wedge dy =\pi.
$$

In general for any surface and for any section $\psi$ satisfying
$\D \psi = 0$ (i.e. we do not assume here that $\psi$ does not
vanish anywhere) we have
$$
\int_M U^2 dx \wedge dy = \pi (-\deg E_0 + \sum_p \ord_p \psi) + \int_M |A|^2 e^{-2\alpha} dx
\wedge dy \geq
$$
$$
\geq \pi (-\deg E_0 + \sum_p \ord_p \psi)
$$
(see \cite{PP}). The integrand $|A|^2 e^{-2\alpha}$ has
singularities in the zeros of $\psi$ however the integral
converges and is non-negative.

Returning to the case of spin bundles over spheres ($\deg E_0 = g-1=-1$) and assuming that
$\dim_\H \Ker \D = N$ we take a point
$p$ and choose a function $\psi \in \Ker \D$ such that $\ord_p
\psi = \dim \Ker_\H \D -1 = N-1$.
Now we substitute $\psi$ in (\ref{gaussquat}) and obtain
$$
\int_M U^2 dx \wedge dy = \pi(1 + N-1) + \int_M |A|^2 e^{-2\alpha}
dx \wedge dy \geq \pi N.
$$
However this estimate is too rough since we see from the proof of
Theorem \ref{willmorenumber} that not only the function from $\Ker
\D$ with the maximal order of zeros contributes to lower bounds
for the Willmore functional and it needs to consider the flag of
functions.

In \cite{FLPP} the deep analogy of this problem to the Pl\"ucker
formulas which relate the degrees and the ramification indices of
the curves associated to some algebraic curve in $\C P^n$ was
discovered. This enables to write down this flag and count the
contribution of the whole kernel of $\D$ into the Willmore
functional. Finally this led to the establishing of the estimates
for the Willmore functional which are quadratic in $\dim_\H \Ker
\D$.

To formulate the main result of \cite{FLPP} we introduce some definitions. Let $H$ be a subspace of $\Ker \D$.
For any point $p$ we put
$$
n_0(p) = \min \ord_p \psi  \ \ \ \mbox{for $\psi \in H$}.
$$
Then successively we define
$$
n_k(p) = \min \ord_p \psi \ \ \ \mbox{for $\psi \in H$ such that $\ord_p \psi > n_{k-1}(p)$}.
$$
We have the Weierstrass gap sequence
$$
n_0(p) < n_1(p) < \dots < n_{N-1}(p), \ \ \ N = \dim_\H H,
$$
and a chain of embeddings
$$
H = H_0 \supset H_1 \supset \dots \supset H_{N-1}
$$
with $H_k$ consisting of $\psi$ such that $\ord_p \psi \geq
n_k(p)$. Then we define the order of a linear system $H$ at the
point $p$ as
$$
\ord_p H = \sum_{k=0}^{N-1}(n_k(p) - k) = \sum_{k=0}^{N-1} n_k(p) - \frac{1}{2}N(N-1).
$$
We say that $p$ is a Weierstrass point if $\ord_p H \neq 0$.

Now we can formulate the main result of this theory:

\begin{theorem}
[\cite{FLPP}] \label{theoFLPP} Let $H \subset \Ker \D$ and $\dim_\H
H = N$. Then
\beq
\label{plucker}
\int_M U^2 dx \wedge dy =\pi( N^2(1-g) + \ord H) + {\cal A}(M),
\eeq
where the term ${\cal A}(M)$
is non-negative and reduces to 
the term $\int_M |A|e^{-2\alpha} dx \wedge dy$
in the case of (\ref{gaussquat1}).
\end{theorem}

In fact the main result of
\cite{FLPP} works for Dirac operators $\D$ with complex conjugate potentials $U=\bar{V}$
\footnote{In this case ${\cal W} = \int UV dx \wedge dy = \int |U|^2 dx \wedge dy$
where ${\cal W}$ is the Willmore functional.}
on arbitrary complex quaternionic line bundles of
arbitrary degree $d$ (this is straightforward from the proof) and
explains the term ${\cal A}(M)$ in the terms of dual curves:
$$
\int_M |U|^2 dx \wedge dy = \pi( N((N-1)(1-g)-d) + \ord H) + {\cal A}(M),
$$
where $g$ is the genus of $M$ and $d = \deg L = \deg E_0$.
In Theorem \ref{theoFLPP} we assume that $d=g-1$, i.e. the case which is interesting
for the surface theory.

For $U=0$ we have also ${\cal A}(M)=0$ and the
Pl\"ucker formula reduces to the original Pl\"cker relation for
algebraic curves (see, for instance, \cite{GH}):
$$
\ord H = N((N-1)(g-1)+d).
$$

For $g=0$ we have

\begin{corollary}
[\cite{FLPP}]
Conjecture \ref{taimconj} is valid: $\int U^2 dx \wedge dy \geq \pi N^2$.
\end{corollary}

For $g \geq 1$ we have an effective lower bound only in terms of $\ord H$ because the term quadratic in $N$
vanishes for $g=1$ and is negative for $g > 1$.

For obtaining effective lower bounds for the Willmore functional in \cite{FLPP}
it was proposed to use some special linear systems
$H$. Let $\dim_\H \Ker \D = N$. We take in $\Ker \D$ a $k$-dimensional
linear system $H$ distinguished by
the condition that for all $\psi \in H$ we have $\ord_p \psi \geq N-k$ for some fixed point $p$.
The Weierstrass gap sequence at this point meets the inequality
$$
n_l(p) \geq N-k+l, \ \ \ l=0,\dots,k-1,
$$
and therefore $\ord_p H \geq k(N-k)$.
From (\ref{plucker}) we have
\beq
\label{pluck1}
\int U^2 dx \wedge dy \geq \pi(k^2(1-g) + k(N-k)) = kN - k^2 g.
\eeq

If $g=0$ then the right-hand side attains its maximum at $k=N$
and we have the estimate (\ref{quadratic}).

If $g\geq 1$ then the function $f(x) = xN - x^2 g$ attains its maximum at $x_{\max} = \frac{N}{2g}$.
Therefore the right-hand side in (\ref{pluck1}) attains its maximum either at
$k = \left[\frac{N}{2g}\right]$ or at $\left[\frac{N}{2g}\right]+1$, the integer point which is
closest to $x_{\max}$. From that it is easy to derive the rough lower bound valid for all $g$.
Of course for special cases this bound can be improved as, for instance, in the case $g=1$:

\begin{corollary}
[\cite{FLPP}]
We have
\beq
\label{pluck2}
\int U^2 dx \wedge dy \geq \frac{\pi}{4g} \left( N^2 - g^2 \right)
\eeq
for $g>1$ and
\beq
\label{pluck3}
\int U^2 dx \wedge dy \geq
\begin{cases}
\frac{\pi N^2}{4} & \text{for $N$ even} \\
\frac{\pi (N^2-1)}{4} & \text{for $N$ odd}
\end{cases}
\eeq
for $g=1$.
\end{corollary}

The proof of Theorem \ref{theofriedrich} works straightforwardly for deriving the following corollary.

\begin{corollary}
[\cite{FLPP}]
Given an eigenvalue $\lambda$ of the Dirac
operator on a two-dimensional spin-manifold of genus $g$,
the following inequalities hold:
$$
\lambda^2 \mathrm{Area}(M) \geq
\begin{cases}
\pi m^2(\lambda) & \text{for $g=0$} \\
\frac{\pi}{g}(m^2(\lambda) - g^2) & \text{for $g \geq 1$},
\end{cases}
$$
where $m(\lambda)$ is the multiplicity of $\lambda$.
\end{corollary}

Another important application of (\ref{pluck3}) concerns the lower bounds for
the area of CMC tori in $\R^3$ and minimal tori in $S^3$. One can see from the explicit
construction of the spectral curves (see \S \ref{subsec4.6}) that in both cases
the normalized spectral curves are hyperelliptic curves
$$
\mu^2 = P(\lambda)
$$
such that a pair of branch points correspond to the ``infinities'' $\infty_\pm$.
There are also $2g$ other branch points (here $g$ is the genus of this 
hyperelliptic curve)
at which the multipliers of Floquet functions equal $\pm 1$
(this is by the construction of the spectral curve).
Moreover there are also a pair of points interchanged by the hyperelliptic involution
at which the multipliers are also $\pm 1$ (the tori are constructed via these Floquet functions
as it is shown in \cite{Hitchin} and \cite{Bobenko}). Thus we have the space $F$ with $\dim_\C F =2g+2$
consisting of solutions to $\D \psi =0$ with multipliers $\pm 1$.
Let us take $4$-sheeted covering $\widehat{M}$ of a torus $M$ which doubles both periods.
The pullbacks of the functions from $F$ onto this covering are double-periodic functions,
i.e. they are sections of the same spin bundle over $\widehat{M}$.
The complex dimension of the kernel of $\D$ acting on this spin bundle is at least
$2g+2$ and thus $\dim_\H \Ker \D \geq g+1$.
Applying (\ref{pluck3}) we obtain the lower bounds for 
$\int |U|^2 dx \wedge dy$.
For CMC tori we have $H=1$ and $U=\frac{e^\alpha}{2}$. Thus
$$
\int_{\widehat{M}} U^2 dx \wedge dy = \frac{1}{4}  \mathrm{Area} (\widehat{M}) = \mathrm{Area} (M).
$$
For minimal tori in $S^3$ we have $U = -\frac{ie^\alpha}{2}$ and thus
$\int_{\widehat{M}} |U|^2 dx \wedge dy = \mathrm {Area} (M)$.
We derive

\begin{corollary}
For minimal tori in $S^3$ and CMC tori in $\R^3$ of spectral genus $g$ we have
the following lower bounds for the area:
$$
\mathrm{Area} \geq
\begin{cases}
\frac{\pi (g+1)^2}{4} & \text{for $g$ odd} \\
\frac{\pi ((g+1)^2-1)}{4} & \text{for $g$ even}.
\end{cases}
$$
\end{corollary}

In \cite{FLPP} it is remarked that it follows from \cite{Hitchin}
that for minimal tori in $S^3$ the bound can be improved by replacing $g+1$ by $g+2$.
Moreover it is valid for the energy of all harmonic tori in $S^3$ however we do not discuss
the spectral curves of harmonic tori in this paper.

Recently the genus of the spectral curve was applied by Haskins to
completely another problem: to study special Lagrangian $T^2$-cones in $\C^3$ \cite{Haskins}.
He obtained linear (in the genus) lower bounds for some quantities characterizing the geometric complexity
of such cones and conjectured that these bounds can be improved to quadratic bounds.
We note that the methods of \cite{Haskins} are completely different from methods used in \cite{FLPP}.

The contribution of the term $\ord H$ is easily demonstrated by soliton spheres such that
the poles of the transition coefficient are of the form $(2l+1)i/2$. In this case $\ord \Ker \D$ counts
the gaps in filling these energy levels.

Recently the definition of soliton spheres was generalized in the
spirit of the lower estimates for the Willmore functional: a sphere
is called soliton if for it the ``Pl\"ucker inequality''
$$
\int_M |U|^2 dx \wedge dy \geq \pi( N^2(1-g) + \ord H)
$$
becomes an equality, i.e. ${\cal A}(M)=0$ \cite{Peters}.

As it was showed by Bohle and Peters \cite{BP} 
this class contains many other interesting surfaces.

Before formulating their result we recall that Bryant surfaces are just surfaces of constant mean curvature one in
the hyperbolic three space \cite{Bryant11}. By \cite{BP}, Bryant surface $M$ in the Poincare ball model
$B^3 \subset \R^3$ is a smooth Bryant end if there is a point $p_\infty$ on the asymptotic boundary
$\partial B^3$ such that $M \cup p_\infty$ is a conformally immersed open disc in $\R^3$.
Generally a Bryant surface is called a compact Bryant surface with smooth ends
if it is conformally equivalent to a compact surface
with finitely many punctured points at which the surface have open neighborhoods isometric to smooth Bryant ends.

It is clearly a generalization of minimal surfaces with planar ends.

We have

\bet
[\cite{BP}]
Bryant spheres with smooth ends are soliton spheres.
The possible values of the Willmore functional for such spheres are $4\pi N$ where $N$
is positive natural number which is non-equal to $2,3,5$, or $7$.
\eet

As it was mentioned by Bohle and Peters they were led to this theorem by the observation that the
simplest soliton spheres corresponding to the potentials $U_N = \frac{N}{2\cosh x}$
can be treated as Bryant spheres with smooth ends.
They also announced that all Willmore spheres are soliton spheres
(we remark that by the results of Bryant and Peng the Willmore functional has the same possible
values
for Willmore spheres as for Bryant spheres with smooth ends \cite{Bryant1,Bryant2,Peng}).

\subsection{The Willmore type functionals for surfaces in
three-dimensional Lie groups} \label{subsec5.5}

The formula (\ref{willmoreformula}) shows that it is reasonable to
consider the functional
$$
E(\Sigma) = \int_\Sigma UV dx \wedge dy
$$
for surfaces. For tori it measures the asymptotic flatness of the
spectral curve and for surfaces in $\R^3$ it equals $E =
\frac{1}{4}\W$ (\cite{T1}). In \cite{Berdinsky} this functional was
considered for surfaces in other Lie groups and was called the
energy of a surface. Although the product $UV$ is not always
real-valued for closed surfaces the functional is real-valued and
equals

\bi
\item
for $SU(2)$ \cite{T3}:
$$
E = \frac{1}{4}\int (H^2 +1) d\mu,
$$
i.e. it is a multiple of the Willmore functional;

\item
for $\nil$ \cite{Berdinsky}:
$$
E = \frac{1}{4} \int
\left( H^2 + \frac{\widehat{K}}{4} - \frac{1}{16} \right)
d \mu;
$$

\item
for $\sll$ \cite{Berdinsky}:
$$
E(M) = \frac{1}{4} \int_M \left( H^2 +
\frac{5}{16} \widehat{K} - \frac{1}{4} \right) d \mu;
$$

\item
for surfaces in $\sol$,
since  the potentials have indeterminacies on the zero measure set,
the energy $E$ is correctly defined.
However we do not know until recently its geometric meaning.
\ei

We recall that by $\widehat{K}$ we denote the sectional curvature of the ambient space
along the tangent plane to a surface.

These functionals were not studied and many problems are open:

1) are they bounded from below (some numerical experiments confirm that)?

2) what are their extremals?

3) what are the analogs of the Willmore conjecture for them?

\medskip

\addcontentsline{toc}{section}{Appendix 1. On the existence of the
spectral curve for 
the Dirac operator with $L_2$-potentials}

\subsection*{Appendix 1. On the existence of the spectral curve for
the Dirac operator with $L_2$-potentials}

{\small

In this appendix we expose the proof of Theorem \ref{l2} following
\cite{Schmidt} where, as we think, the exposition is too short.

Moreover the ideas of the proof of this theorem are essential for
proving the main result of \cite{Schmidt2}: the minimum of the
Willmore functional in a given conformal class of surfaces is
constructed as follows. We consider the infimum of the Willmore
functional in this class and take in this class a sequence of
surfaces (or, more precisely their Weierstrass representations)
with the values of the Willmore functional converging to the
infimum. Then there is a weakly converging sequence of potentials
of the corresponding Dirac operators. The Dirac operator with the
limit potential also has a nontrivial kernel (this follows from
the converging of the resolvents) and the desired minimizing
surface is constructed from a function from this kernel by the
Weierstrass representation. Of course it is necessary to control
the smoothness which is possible. However in \cite{Schmidt2} it is
mentioned that one can not say that there are no the absence of
branch points of the limit surface.

The analog of the decomposition (\ref{dec}) is the following
sequence
\begin{equation}
\label{dec2} L_p \stackrel{(\D_0 - E_0)^{-1}}{\longrightarrow}
W^1_p \stackrel{\mathrm{Sobolev's\ embedding}}{\longrightarrow}
L_{\frac{2p}{2-p}}
\stackrel{\mathrm{multiplication}}{\longrightarrow} L_p, \ \ \
p<2.
\end{equation}
All operators coming in the sequence are only continuous and we
can not argue as in \S \ref{subsec4.2}.

Let $M = \C/\Lambda$ be a torus and $z$ be a linear complex
coordinate on $M$ defined modulo $\Lambda$. Denote by
$\rho(z_1,z_2)$ a distance between points $z_1, z_2 \in M$ in the
metric induced by the Euclidean metric on $\C$ via the projection
$\C \to \C/\Lambda$.

The following proposition is derived from the definition of the resolvent
$$
(\D_0 - E) R_0(E) = \delta(z-z^\prime)
$$
by straightforward computations.

\begin{proposition}
The resolvent
$$
R_0(E) = (\D_0 - E)^{-1}: L_2 \to W^1_2  \to L_2
$$
of the free operator $\D_0: L_2 \to L_2$ is an integral matrix
operator
$$
f(z,\bar{z}) \to [R_0(E)f](z,\bar{z}) = \int_M K_0(z,z^\prime,E)
f(z^\prime,\overline{z^\prime}) dx^\prime dy^\prime, \ \ \
z^\prime =x^\prime + i y^\prime,
$$
with the kernel $ K_0(z,\bar{z},z^\prime,\overline{z^\prime},E) =
\left(
\begin{array}{cc}
r_{11} & r_{12} \\
r_{21} & r_{22}
\end{array}
\right), r_{ik} =
r_{ik}(z,\bar{z},z^\prime,\overline{z^\prime},E),$ where
$$
r_{12} = \frac{1}{E} \partial r_{22}, \ \ \ r_{21} = - \frac{1}{E}
\bar{\partial} r_{11},
\ \ \
\frac{1}{E} ( \partial \bar{\partial} + E^2) r_{11} =
\frac{1}{E} ( \partial \bar{\partial} + E^2) r_{22} =
-\delta(z-z^\prime).
$$
\end{proposition}

\begin{corollary}
The integral kernel of $R_0(E)$ equals
$$
K_0(z,z^\prime,E) = \left(
\begin{array}{cc}
-E G & - \partial{G} \\
\bar{\partial} G & - E G
\end{array}
\right),
$$
where $G$ is the (modified) Green on the Laplace operator on the
torus $M$:
$$
(\partial \bar{\partial} + E^2) G(z,z^\prime,E) =
\delta(z-z^\prime).
$$
\end{corollary}

{\sc Example.} Given a torus $M = \C/\{2\pi \Z + 2\pi i \Z\}$, we
have
$$
\delta(z-z^\prime) = \sum_{k,l \in \Z}
e^{i(k(x-x^\prime)+l(y-y^\prime))}, \ \ \ z=x+iy, z^\prime =
x^\prime + i y^\prime,
$$
\begin{equation}
\label{res} G(z,z^\prime,E) = -4 \sum_{k,l \in \Z}
\frac{1}{k^2+l^2 - 4 E^2} e^{i(k(x-x^\prime)+l(y-y^\prime))}.
\end{equation}
For other period lattices $\Lambda$ the analog of series
(\ref{res}) for $G$ looks almost the same and has very similar
analytic properties. We do not write it down and always will refer
to (\ref{res}) when we consider its analytical properties.

The following proposition is clear.

\begin{proposition}
The series (\ref{res}) converges for $E = i\lambda$ with $\lambda
\in \R$ and $\lambda \gg 0$ (i.e. $\lambda$ sufficiently large).
\end{proposition}

For calculating the operator
$$
R_0(k,E) = (\D_0 + T_k -E)^{-1}: L_2 \to W^1_2
\stackrel{\mathrm{embedding}}{\to} L_2
$$
let us use the following identity
$$
(\D_0 + T_k - E) = (1 + T_k(\D_0 -E)^{-1})(\D_0 -E) = (1+T_k
R_0(E))(\D_0-E)
$$
which implies the formula for the resolvent
\begin{equation}
\label{neumann} R_0(k,E) = R_0(E)(1 + T_kR_0(E))^{-1} = R_0(E)
\sum_{l=0}^\infty \left[ -T_k R_0(E)\right]^l
\end{equation}
provided that the series in the right-hand side converges.

{\sc Remark.} Given $p$, $1<p<2$, the symbol
$$
R_0(E) = (\D_0 - E)^{-1} \ \ \ \mbox{or} \ \ \ R_0(k,E) = (\D_0 +
T_k - E)^{-1}
$$
denotes

a) an operator $A: L_p \to W^1_p$;

b) a composition $B: L_p \to L_q$, $q = \frac{2p}{2-p}$, of $A$
and the Sobolev embedding $W^1_p \to L_q$;

c) a composition $C: L_p \to L_p$ of $A$ and the natural embedding
$W^1_p \to L_p$.

\noindent The actions of these operators are the same on the space
of smooth functions which can be considered as embedded into
$W^1_p$ or $L_q$ (all these spaces are the closures of the space
of smooth functions with respect to different norms). Therefore it
is enough to demonstrate or prove all necessary estimates only for
smooth functions and that could be done by using explicit formulas
for the resolvents.

Let us decompose resolvents into sums of integral operators as
follows.

We denote by $\chi_\varepsilon$ the function $\chi_\varepsilon(r)
=
\begin{cases}
0 & \text{for $r > \varepsilon$} \\
1 & \text{for $r \leq \varepsilon$}
\end{cases}
$ defined for $r \geq 0, r \in \R$. Given $\delta >0$, decompose
the resolvent $R_0(k,E)$ into a sum of two integral operators:
$$
R_0(k,E) = R_0^{\leq \varepsilon}(k,E) + R_0^{> \varepsilon}(k,E):
L_p \to L_q
$$
where the ``near'' part $R_0^{\leq \varepsilon}(k,E)$ is defined
by its kernel
$$
K^{\leq \varepsilon}_0(z,\bar{z},z^\prime,\overline{z^\prime},E) =
K_0(z,\bar{z},z^\prime,\overline{z^\prime},E)
\chi_\varepsilon(\rho(z,z^\prime))
$$
and the ``distant'' part $R_0^{> \varepsilon}(k,E)$ has the
following kernel
$$
K^{> \varepsilon}_0(z,\bar{z},z^\prime,\overline{z^\prime},E) =
K_0(z,\bar{z},z^\prime,\overline{z^\prime},E)
(1-\chi_\varepsilon(\rho(z,z^\prime))).
$$

\begin{proposition}
\label{proposition-s1} Given $p$, $1<p<2$,
$\hat{k}=(\hat{k}_1,\hat{k}_2) \in \C^2$ and $\delta, 0<
\delta<1$, there exists a real constant $\lambda_0 >>0$ such that
$$
\|T_k R_0(i\lambda)\|_{L_p \to L_p} < \delta
$$
for all $\lambda > \lambda_0$ and for $k$ sufficiently close to
$\hat{k}$.

Therefore for such $\lambda$ and $k$

1) the series in (\ref{neumann}) converges and defines a bounded
operator from $L_p$ to $W^1_p$, $L_q$ or $L_p$ (this depends on
the meaning of the symbol $R_0(E)$ multiplied from the left with
the series);

2) the norm of the operator
$$
R_0(k,i\lambda): L_p \to W^1_p
$$
is bounded by some constant $r_p$;

3) given $\varepsilon > 0$, we have
$$
\lim_{\lambda \to \infty} \|R^{>\varepsilon}_0(k,i\lambda)\|_{L_p
\to L_q} = 0, \ \ \ q=\frac{2p}{2-p}.
$$
\end{proposition}

This proposition follows from the explicit formula (\ref{res}) for
the kernel of resolvent.

We denote by $r_{\mathrm{inj}}$ the injectivity radius of the
metric on $M$ and introduce the norms $\|\cdot\|_{2;\varepsilon}$
defined for $0< \varepsilon < r_{\mathrm{inj}}$ as follows. Given
$U \in L_2(M)$, we denote by $U\vert_{B(z,\varepsilon)}$ the
restriction of $U$ onto the ball $B(z,\varepsilon) = \{w \in M\ :
\ \rho(z,w) < \varepsilon\}$ and define $\|U\|_{2;\varepsilon}$ as
$$
\|U\|_{2;\varepsilon} = \max_{z \in M}
\left\|U\vert_{B(z,\varepsilon)}\right\|_2.
$$

\begin{proposition}
\label{localnorm} 1) There are the inequalities
$$
\sqrt{\frac{\pi \varepsilon^2}{\vol (M)}} \leq
\|U\|_{2;\varepsilon} \leq \|U\|_2
$$
for all $U \in L_2(M)$.

2)  For all $C > 0$ and $\varepsilon$ the sets
$\{\|U\|_{2;\varepsilon} \leq C\}$ are closed and therefore are
compact in both of the weak and the weak convergence
topologies on $L_2(M)$.
\end{proposition}

{\sc Proof.} It is clear that $\|U\|_{2;\varepsilon} \leq
\|U\|_2$. Moreover we have
$$
\|U\|^2_{2;\varepsilon} \vol(M) \geq \int_M
\int_{B(z,\varepsilon)}
|U(z+z^\prime,\bar{z}+\overline{z^\prime})|^2 dz^\prime dz =
$$
$$
= \int_M \int_{B(0,\varepsilon)}
|U(z+z^\prime,\bar{z}+\overline{z^\prime})|^2 dz^\prime dz =
\int_{B(0,\varepsilon)} \left[\int_M |U(z,\bar{z})|^2 dz \right]
dz^\prime = \pi \varepsilon^2 \|U\|^2_2
$$
where $dz = dx \wedge dy, dz^\prime = dx^\prime \wedge dy^\prime$.
The second statement is known from a course on functional
analysis.

Consider the resolvent
$$
R(k,E) = (\D + T_k  -E)^{-1}: L_p \to L_p.
$$
We again use an identity
$$
(\D + T_k -E) = \left[ 1 + \left(
\begin{array}{cc}
U & 0 \\ 0 & V
\end{array}
\right) (\D_0 + T_k - E)^{-1} \right] (\D_0 + T_k - E)
$$
which implies
$$
R(k,E) =
R_0(k,E) \sum_{l=0}^\infty \left[ - \left(
\begin{array}{cc}
U & 0 \\ 0 & V
\end{array}
\right) R_0(k,E) \right]^l.
$$

\begin{proposition}
\label{proposition-s2} Let  $1<p<2$,
$\hat{k}=(\hat{k}_1,\hat{k}_2) \in \C^2$, $\varepsilon$ be
sufficiently small and $0< \delta<1$. There is $\gamma
> 0$ such that
$$
\left\| \left(
\begin{array}{cc}
U & 0 \\ 0 & V
\end{array}
\right) R_0(k,i\lambda) \right\|_{L_p \to L_p} < \delta
$$
for all $\lambda > \lambda_0$, $k$ sufficiently close to $\hat{k}$
and $U,V$ with $\|U\|_{2;\varepsilon} < \gamma$,
$\|V\|_{2;\varepsilon} < \gamma$.
\end{proposition}

{\sc Proof.} We have an obvious inequality
$$
\|R_0^{\leq \varepsilon}(k,E)\| \leq \|R_0(k,E)\| \ \ \ \mbox{for
all $\varepsilon$.} $$ Let $S_p$ be the Sobolev constant for the
embedding $W^1_p \to L_q$ (see Proposition \ref{proposition1}).
For $\lambda > \lambda_0$ we have
$$
\|R_0(k,i\lambda)\|_{L_p \to W^1_p} \leq r_p
$$
(see Proposition \ref{proposition-s1}). Now consider the
composition of mappings
$$
L_p \stackrel{R_0^{\leq\varepsilon}(k,E)}{\longrightarrow} W^1_p
\stackrel{\mathrm{embedding}}{\to} L_q \stackrel{\times \left(
\begin{array}{cc}
U & 0 \\ 0 & V
\end{array}
\right)}{\longrightarrow} L_p,
$$
where the norm of the first mapping is bounded from above by
$r_p$, the norm of the second mapping is bounded from above by
$S_p$. Let us compute the norm of the third mapping.

Since the integral kernel of $R_0^{\leq \varepsilon}(k,E)$ is
localized in the closed domain $\{\rho(z,z^\prime) \leq
\varepsilon\}$, for any ball $B(x,\alpha)$ we have
$$
\left[ \left(
\begin{array}{cc}
U & 0 \\ 0 & V
\end{array}
\right) R_0^{\leq \varepsilon}(k,E) f
\right]\Big\vert_{B(x,\alpha)} = \left[\left(
\begin{array}{cc}
U & 0 \\ 0 & V
\end{array}
\right) R_0^{\leq \varepsilon}(k,E)\right]
\left(f\vert_{B(x,\alpha+\varepsilon)}\right).
$$
Applying the H\"older inequality to the right-hand side of this
formula we obtain
$$
\left\| \left[ \left(
\begin{array}{cc}
U & 0 \\ 0 & V
\end{array}
\right) R_0^{\leq \varepsilon}(k,E) f
\right]\Big\vert_{B(x,\alpha)} \right\|_p \leq
$$
$$
\leq m \|R^{\leq \varepsilon}_0(k,E)\|_{L_p \to L_q}
\left\|(f\vert_{B(x,\alpha+\varepsilon)}\right\|_p \leq \leq m r_p
S_p \left\|(f\vert_{B(x,\alpha+\varepsilon)}\right\|_p
$$
where $m = \max (\|U\|_{2;\varepsilon},\|V\|_{2;\varepsilon})$.
Now we recall the identity
$$
\int_M \left\|g_{B(x,\alpha)}\right\|^p_p dx = \vol
B(x,\alpha)\|g\|^p_p = \pi \alpha^2 \|g\|^p_p
$$
and applying it to the previous inequality obtain
$$
\left\| \left(
\begin{array}{cc}
U & 0 \\ 0 & V
\end{array}
\right) R_0^{\leq \varepsilon}(k,E) \right\|_p \leq m r_p S_p
\left(1+ \frac{\varepsilon}{\alpha}\right)^{2/p}.
$$
Since we use the Sobolev constant $S_p$ for the torus, we have to
assume that $(\alpha+\varepsilon) < r_{\mathrm{inj}}$. If
\begin{equation}
\label{gamma} m = \max
\left(\|U\|_{2;\varepsilon},\|V\|_{2;\varepsilon}\right) <
\frac{\delta}{r_p S_p \root p \of 4}
\end{equation}
and $\alpha = \varepsilon < r_{\mathrm{inj}}/2$, we have
$$
\left\| \left(
\begin{array}{cc}
U & 0 \\ 0 & V
\end{array}
\right) R_0^{\leq \varepsilon}(k,E) \right\|_p < \delta.
$$
This proves the proposition.

\begin{proposition}
[\cite{Schmidt}] \label{weaktheorem} Let $p$ and $\hat{k} \in
\C^2$ be the same as in Proposition \ref{proposition-s2}, let
$\gamma < (r_p S_p \root p \of 4)^{-1}$ and let $\lambda >> 0$,
i.e. $\lambda$ be sufficiently large. Given sufficiently small
$\varepsilon > 0$, for $U,V$ such that $\|U\|_{2;\varepsilon} \leq
C \leq \gamma, \|V\|_{2;\varepsilon} \leq C \leq \gamma$ the
series
\begin{equation}
\label{res3} R(k,i\lambda) = R_0(k,i\lambda) \sum_{l=0}^\infty
\left[ - \left(
\begin{array}{cc}
U & 0 \\ 0 & V
\end{array}
\right) R_0(k,i\lambda) \right]^l
\end{equation}
converges uniformly near $\hat{k}$ and defines the resolvent of
operator
$$
\D+T_k: L_p \to L_p.
$$
The action of this resolvent on smooth functions is extended to
the resolvent of $\D+T_k$ on the space $L_2$:
$$
(\D + T_k - E)^{-1}: L_2 \to W^1_2
\stackrel{\mathrm{embedding}}{\longrightarrow} L_2.
$$
This is a pencil of compact operators holomorphic in $k$ in near
$\hat{k}$. If $(U_n,V_n)
\stackrel{\mathrm{weakly}}{\longrightarrow} (U_\infty,V_\infty)$
in $\{\|U\|_{2;\varepsilon} \leq C$, $\|V\|_{2;\varepsilon} < C\}$
then the corresponding resolvents converges to the resolvent of
the operator with potentials $(U_\infty,V_\infty)$ in the normed
topology.
\end{proposition}

{\sc Proof.} By Proposition \ref{proposition-s2} and
(\ref{gamma}), for $\lambda > \lambda_0$ we have
$$
\left\| \left(
\begin{array}{cc}
U & 0 \\ 0 & V
\end{array}
\right) R_0^{\leq \varepsilon}(k,i\lambda) \right\|_p  < \sigma =
\gamma r_p S_p \root p \of 4 < 1
$$
near $\hat{k}$. By Proposition \ref{proposition-s1}, for
sufficiently large real $\lambda$, since the norm of the embedding
$L_q \to L_p$ is bounded, we have
$$
\left\| \left(
\begin{array}{cc}
U & 0 \\ 0 & V
\end{array}
\right) R_0^{> \varepsilon}(k,i\lambda) \right\|_p  < 1- \sigma.
$$
This implies that for $\lambda >>0$ we have
$$
\left\| \left(
\begin{array}{cc}
U & 0 \\ 0 & V
\end{array}
\right) R_0(k,i\lambda) \right\|_p \leq \left\| \left(
\begin{array}{cc}
U & 0 \\ 0 & V
\end{array}
\right) \left( R_0^{\leq \varepsilon}(k,i\lambda) +
R_0^{>\varepsilon}(k,i\lambda) \right) \right\|_p  < 1
$$
and the series in (\ref{res3}) uniformly converges near $\hat{k}$
and defines the resolvent of $\D + T_k: L_p \to L_p$.

The action of $R(k,i\lambda)$ on smooth functions is given by
(\ref{res3}) and we extend it to a compact operator on $L_2$ as
follows. Put
$$
B = \sum_{l=0}^\infty \left[ - \left(
\begin{array}{cc}
U & 0 \\ 0 & V
\end{array}
\right) R_0(k,i\lambda) \right]^l
$$
and consider the following composition of operators
$$
L_2 \stackrel{\mathrm{embedding}}{\longrightarrow} L_p
\stackrel{B}{\to} L_p \stackrel{(\D +T_k
-E)^{-1}}{\longrightarrow} W^1_p
\stackrel{\mathrm{embedding}}{\longrightarrow} L_2
$$
where all operators are bounded and the embedding $W^1_p \to L_2$
is compact by the Kondrashov theorem (see Proposition
\ref{proposition1}). This shows that the action of $R(k,i\lambda)$
on smooth functions is extended to a compact operator on $L_2$.
Since the series (\ref{res3}) is holomorphic in $k$ the resolvent
$R(k,i\lambda)$ is also holomorphic in $k$.

Now we are left to prove that the resolvent is continuous in $U$
and $V$. Every entrance of the matrix $\left(
\begin{array}{cc} U & 0 \\ 0 & V
\end{array}\right)$ in any term of (\ref{res3})
is dressed from both sides by the resolvents $R_0(k,i\lambda)$
which are bounded integral operators. Let $l=1$ and let
$K(z,z^\prime,k,i\lambda)$ be the kernel of such an operator. Then
the composition
$$
R_0(k,i\lambda) \left(
\begin{array}{cc} U & 0 \\ 0 & V
\end{array}\right)
R_0(k,i\lambda)
$$
acts on smooth functions as the integral operator with the kernel
$$
F(z,z^{\prime\prime}) = K(z,z^\prime,k,i\lambda) \left(
\begin{array}{cc} U(z^\prime) & 0 \\ 0 & V(z^\prime)
\end{array}\right)
K(z^\prime,z^{\prime\prime},k,i\lambda).
$$
Obviously such an integral operator is continuous with respect to
the weak convergence of potentials $U,V \in L_2(M)$. For other
values of $l$ the proof is analogous. By Proposition
\ref{localnorm}, every term of the series (\ref{res3}) is
continuous with respect to the weak convergence of potentials
$\{\|U\|_{2;\varepsilon} \leq C, \|V\|_{2;\varepsilon} \leq C \}$.
Since the series (\ref{res3}) uniformly converges, the same
continuosity property holds for the sum of the series. This proves
the proposition.

This proposition establishes the existence of the resolvent only
for large values of $\lambda$ where $E = i\lambda$. The resolvent
is extended to a meromorphic function onto the $E$-plane by using
the Hilbert formula (see Proposition \ref{proposition2}).

{\sc Proof of Theorem \ref{l2}.} By Proposition \ref{weaktheorem}
there are $k_0 \in \C^2$ and $E \in \C$ such that the operator
$$
(\D + T_{k_0} - E_0)^{-1}: L_2 \to W^1_2
\stackrel{\mathrm{embedding}}{\longrightarrow} L_2
$$
is correctly defined. Let us substitute the expression $\varphi =
(\D + T_{k_0} - E)^{-1}f$ into the equation
$$
(\D + T_k - E)\varphi = 0
$$
and rewrite this equation in the form
$$
(\D + T_{k_0} - T_{k_0} + T_k -E_0 + E_0 - E) (\D + T_{k_0} -
E_0)^{-1} f = \left[ 1 + A_{U,V}(k,E) \right] f = 0
$$
where
$$
A_{U,V}(k,E) = (T_k - T_{k_0} + E_0 - E)(\D + T_{k_0} - E)^{-1}.
$$
Since the first multiplier in this formula is a bounded operator
for any $k,E$ and the second multiplier is a compact operator,
$A_{U,V}(k,E)$ is a pencil of compact operators which is
polynomial in $k$ and $E$. By applying the Keldysh theorem as in
\S \ref{subsec4.2} we derive the theorem.

Now the spectral curve is defined as usual by the formula
$$
\Gamma = Q_0(U,V)/\Lambda^\ast.
$$

{\sc Remark.} The resolvents of operators on noncompact spaces
does not behave continuously under the weak convergence of
potentials. Indeed, consider the Schr\"odinger operator
$$
L = -\frac{d^2}{dx^2} + U(x)
$$
where $U(x)$ is a soliton potential (so the operator does have
bounded states). The isospectral sequence  of potentials $U_N(x) =
U(x+N)$ weakly converges to the zero potential $U_\infty = 0$ for
which the Schr\"odinger operator has no bounded states. The same
is true for the one-dimensional Dirac operator.

}

\addcontentsline{toc}{section}{Appendix 2. The conformal Gauss map
and the conformal area}

\subsection*{Appendix 2. The conformal Gauss map and the conformal area}

{\small

In this appendix we expose the known results on the Gauss
conformal mapping mostly following \cite{Bryant1,Eschenburg,FPPS}.

We denote by $S_{q,r}$ the round sphere of radius $r$ in $\R^3$
and with the center at $q$ and denote by $\Pi_{p,N}$ the plane, in
$\R^3$, passing through $p$ and with the normal vector $N$. All
such spheres and planes in $\R^3$ are parameterized by a quadric
$Q^4 \subset \R^{4,1}$. Indeed, let
$$
\langle x, y \rangle = x_1 y_1 + \dots + x_4 y_4 - x_5 y_5
$$
be an inner product in $\R^{4,1}$. Put
$$
Q^4 = \{\langle x, x \rangle = 1\} \subset \R^{4,1},
$$
$$
S_{q,r} \rightarrow \frac{1}{r}\left( q, \frac{1}{2}(|q|^2 - r^2
-1),  \frac{1}{2}(|q|^2 - r^2 +1) \right),
$$
$$
\Pi_{q,N} \rightarrow (N, \langle q, N \rangle,\langle q, N
\rangle).
$$

Given a surface $f: \Sigma \to \R^3$, its {\it conformal Gauss
map}
$$
G^c: \Sigma \rightarrow Q^4
$$
corresponds to a point $p \in \Sigma$ a sphere of radius
$\frac{1}{H}$ which touches the surface at $p$ for $H \neq 0$:
$$
G^c(p) = S_{p+N/H,1/H},
$$
and it corresponds to a point the tangent plane at this point for
$H \neq 0$. In terms of the coordinates on $Q^4$ it is written as
$$
G^c(p) = H\cdot X + T
$$
where
$$
X = \left(f, \frac{(f, f) -1}{2}, \frac{(f, f) + 1}{2}\right), \ \
\ \ T = (N, (N, f), (N, f)).
$$
This mapping is a special case of so-called sphere congruences
which is one of the main subjects of conformal geometry (the
recent statement of this theory is presented in \cite{Hertrich}).

We have $\langle X,X \rangle = 0, \langle T,T \rangle  = 1$, and $\langle X,T \rangle = 0$ which implies that
$\langle dX,X \rangle = \langle dT,T\rangle = 0, \langle dT,X\rangle = \langle -dX,T\rangle$.
It is easily checked that
$$
\langle dX, T \rangle = (df,N) = 0, \ \ \
\langle dX, DX \rangle = (df,df) = {\bf I},
$$
$$
\langle dX,dT \rangle = (df,dN) = -{\bf II}, \ \ \
\langle dT, dT \rangle = (dN,dN) = {\bf III},
$$
where the third fundamental form ${\bf III}$ of a surface measures the
lengths of images of curves under the Gauss map and meets the
identity
$$
K \cdot {\bf I} - 2H \cdot {\bf II} + {\bf III}
$$
which relates it to ${\bf I}$ and ${\bf II}$, the first and the second
fundamental forms of a surface.
It implies that
$$
\langle Y_z, Y_z \rangle = \langle Y_{\bar{z}}, Y_{\bar{z}}
\rangle = 0, \ \ \ \langle Y_z, Y_{\bar{z}} \rangle = e^\beta =
\frac{(H^2-K)e^{2\alpha}}{2} = (H^2-K) (f_z,f_{\bar{z}})
$$
where for brevity we denote $G^c$ by $Y$, $z$ is a conformal
parameter on the surface, and ${\bf I} = e^{2\alpha}dz d\bar{z}$
is the induced metric on the surface. We conclude that

\bi
\item
{\sl the conformal Gauss map is regular and conformal outside
umbilic points.}
\ei

It is
clear that $X$ and $Y$ are linearly independent. Outside umbilics
the set of vectors $Y, Y_z, Y_{\bar z}$, and $X$ is uniquely
completed by a vector $Z \in \R^5$ to a basis
$$
\sigma = (Y, Y_z, Y_{\bar z}, X, Z)^T
$$
for $\C^5$, the complexification of $\R^5$,  such that the inner
product in $\R^{4,1}$ takes the form
$$
\left(
\begin{array}{ccccc}
1 & 0 & 0 & 0 & 0\\
0 & 0 & e^{\beta} & 0 & 0\\
0 &  e^{\beta} & 0 & 0 & 0\\
0 & 0 & 0 & 0 & 1\\
0 & 0 & 0 & 1 & 0
\end{array}
\right).
$$

The analogs of the Gauss--Weingarten equations are
$$
\sigma_z = \U \sigma, \ \ \ \sigma_{\bar{z}} = \V \sigma,
$$
$$
\U = \left(
\begin{array}{ccccc}
0 & 1 & 0 & 0 & 0 \\
0 & \beta_z & 0 & C_2 & C_1 \\
- e^{\beta} & 0 & 0 & C_4 & C_3 \\
0 & - e^{-\beta} C_3 & - e^{-\beta} C_1  & C_5 & 0 \\
0 & - e^{-\beta} C_4 & - e^{-\beta} C_2 & 0 & - C_5
\end{array}
\right),
$$
$$
\V = \left(
\begin{array}{ccccc}
0 & 0 & 1 & 0 & 0 \\
- e^{\beta} & 0 & 0 & C_4 & C_3 \\
0 & 0 & \beta_{\bar{z}} & \bar{C}_2 & \bar{C}_1 \\
0 & - e^{-\beta} \bar{C}_1 & - e^{-\beta} \bar{C}_3 & \bar{C}_5 & 0 \\
0 & - e^{-\beta} \bar{C}_2 & - e^{-\beta} \bar{C}_4 & 0 & -
\bar{C}_5
\end{array}
\right)
$$
with
$$
C_1 = \langle Y_{z z}, X \rangle_{4,1}, \ \ \ C_2 = \langle Y_{z
z}, Z \rangle_{4,1}, \ \ \ C_3 = \langle Y_{z \bar{z}}, X
\rangle_{4,1},
$$
$$
C_4 = \langle Y_{z \bar{z}}, Z \rangle_{4,1}, \ \ \ C_5 = \langle
X_z, Z \rangle_{4,1}.
$$

It is checked by straightforward computations that
$$
\Delta Y + 2(H^2 - K) Y = (\Delta H + 2H(H^2 - K)) X
$$
which, in particular,  implies
$$
C_3 = 0, \ \ \ C_4 = \frac{e^{2\alpha}}{4}(\Delta H + 2H (H^2 -
K)).
$$
Here $\Delta = 4e^{-2\alpha}\partial \bar{\partial}$ stands for
the Laplace--Beltrami operator on the surface. Taking this into
account and keeping in mind that $C_4$ is real-valued, we derive
the Codazzi equations for the conformal Gauss map:
\beq
\begin{split}
\beta_{z\bar{z}} + e^{\beta} - (\bar{C}_1 C_2 + C_1\bar{C}_2) e^{-\beta} =0, \\
C_{1\bar{z}} = C_1 \bar{C}_5, \\
C_{2\bar{z}} + C_2 \bar{C}_5 = C_{4 z} - \beta_z C_4 + C_4 C_5, \\
C_{5 \bar{z}} - \bar{C}_{5 z} = e^{-\beta} (C_1\bar{C}_2 - \bar{C}_1 C_2).
\end{split}
\eeq
By straightforward computations,
we obtain
$$
C_1 = A = \langle N, f_{z z} \rangle , \ \ \ e^{\beta} = 2 |A|^2 e^{-2\alpha}.
$$

{\it The conformal area} $V^c$ of $\Sigma$ is the area of its
image in $Q^4$:
$$
V^c(\Sigma) = \int_{\Sigma} (H^2 - K) d\mu
$$
where $d\mu$ is the volume form on $\Sigma$. The Euler-Lagrange
equation for $V^c$ is
$$
\Delta H + 2H(H^2 - K) = 0.
$$
A surface in $\R^3$ is called {\it conformally minimal} (or {\it
Willmore surface}), if it satisfies this equation. We conclude
that

\bi
\item
{\sl conformally minimal surfaces are
exactly surfaces whose $G^c$-images are minimal surfaces in $Q^4$.}
\ei

Given a non-umbilic point $p \in \Sigma$, the tangent space to
$Q^4$ at  $Y(p)$ is spanned by $Y_z, Y_{\bar{z}}, X$, and $Z$. We
see that $Y$ is {\it conformally harmonic}, i.e., $\Delta Y$ is
everywhere orthogonal to tangent planes to $Q^4$, if and only if
the surface is conformally minimal.

It follows from the Gauss--Weingarten equations for $G^c$ and the
Euler--Lagran\-ge equation for $V^c$ that if $\Sigma$ is
conformally minimal, i.e. $C_4 = 0$, then the quartic differential
$$
\omega = \langle Y_{z z}, Y_{z z} \rangle \, \left(d z\right)^4 =
C_1 C_2 \left(dz\right)^4
$$
is holomorphic.

We recall that a holomorphic quartic differential on a $2$-sphere
vanishes: $\omega = 0$, and any such a differential on a torus is
constant: $\omega = {\mathrm{const}} \cdot \left(dz\right)^4$.

A minimal surface in $Q^4$ is called {\it superminimal} if $\omega
=0$.

We put
$$
\varphi = \log \frac{\bar{C}_1}{C_1}.
$$

We notice that $C_1 \equiv 0$ only for a surface consisting of
umbilics and, by the Hopf theorem, this is a domain in a round
sphere in $\R^3$ or in the plane.

If $\omega \equiv 0$ and $C_1 \neq 0$ then $C_2 \equiv 0$ and the
Gauss--Codazzi equations for the conformal Gauss mapping reduce to
$$
\beta_{z\bar{z}} + e^{\beta} = 0, \ \ \ \ \varphi_{z\bar{z}} = 0.
$$
The first of these equations is the Liouville equation and the
second one is the Laplace equation. These equations describe
superminimal surfaces which are not umbilic surfaces.

Let us consider the case when a conformally minimal surface is not
superminimal. Locally by changing a conformal parameter we achieve
that
$$
\frac{1}{2} \langle Y_{z z}, Y_{z z} \rangle_{4,1} = C_1 C_2 =
\frac{1}{2}.
$$
Then  the Gauss--Codazzi equations
take the form
$$
\beta_{z\bar{z}} + e^{\beta} - e^{-\beta} \cosh \varphi = 0, \ \ \
\
\varphi_{z\bar{z}} + e^{-\beta} \sinh \varphi = 0,
$$
which is the four-particle Toda lattice.

}

\addcontentsline{toc}{section}{Appendix 3. The inverse spectral
problem for the Dirac operator on the line}

\subsection*{Appendix 3. The inverse spectral
problem for the Dirac operator on the line and the trace formulas}

{\small

Here being mostly oriented to geometers we expose some facts which
are necessary for proving Theorem \ref{willmorenumber} and
introducing soliton spheres in \S \ref{subsec5.4}. 

The inverse scattering problem for the Dirac operator on the line
was solved in \cite{ZS} similarly to the same problem for the
Schr\"odinger operator $-\partial^2_x + u(x)$ \cite{Faddeev1} (see
also \cite{Marchenko}).

We consider the following spectral problem, i.e. the
Zakharov--Shabat problem,
\begin{equation}
\label{zs1}
L \psi = \left(
\begin{array}{cc}
0 & ik \\
ik & 0
\end{array}
\right) \psi
\end{equation}
where
\begin{equation}
\label{zs2}
L = \left(
\begin{array}{cc}
0 & \partial_x \\
- \partial_x & 0
\end{array}
\right) + \left(
\begin{array}{cc}
p & 0 \\
0 & q
\end{array}
\right).
\end{equation}

We assume that the potentials $p$ and $q$ are fast decaying as $x
\to \pm \infty$. It is clear from the proofs that it is enough to
assume that $p(x)$ and $q(x)$ are exponentially decaying.

For $p=q=0$ for each $k \in \R \setminus\{0\}$ we have a
two-dimensional space of solutions (free waves) spanned by the
columns of the matrix
$$
\Phi_0(x,k) = \left(
\begin{array}{cc}
0 & e^{-ikx} \\
e^{ikx} & 0
\end{array}\right).
$$
For nontrivial $p$ and $q$ for each $k \in \R \setminus\{0\}$ we
have again a two-dimensional space of solutions which are asymptotic
to free waves when $x \to \pm\infty$. These spaces are spanned by
the so-called Jost functions $\varphi^\pm_l, l=1,2$. For defining
these functions we consider the matrices $\Phi^+(x,k)$ and
$\Phi^-(x,k)$ satisfying the integral equations
$$
\Phi^+ (x,k) = \Phi_0(x,k) + \int_x^{+\infty} \Phi_0(x-x^\prime,k)
\left(
\begin{array}{cc}
p & 0 \\
0 & q
\end{array}\right)
\Phi^+(x^\prime,k) \, d x^\prime,
$$
$$
\Phi^- (x,k) = \Phi_0(x,k) + \int^x_{-\infty} \Phi_0(x-x^\prime,k)
\left(
\begin{array}{cc}
p & 0 \\
0 & q
\end{array}\right)
\Phi^-(x^\prime,k) \, d x^\prime.
$$
These equations are of the form $ \Phi^\pm = \Phi_0 + A^\pm
\Phi^\pm$ where $A^\pm$ are operators of the Volterra type and
therefore each of these equations has a unique solution given by
the Neumann series $\Phi^\pm(x,k) = \sum_{l=0}^\infty
\left(A^\pm\right)^l \Phi_0(x,k)$. The columns of $\Phi^\pm$ are
the Jost functions $\varphi^\pm_l, l=1,2$. We see that by the
construction the Jost functions behave asymptotically as the free
waves:
$$
\varphi^\pm_1 \approx \left( \begin{array}{c} 0 \\ e^{ikx}
\end{array}\right), \ \ \
\varphi^\pm_2 \approx \left( \begin{array}{c} e^{-ikx} \\ 0
\end{array}\right) \ \ \mbox{as $x \to \pm \infty$}.
$$

By straightforward computations it is obtained that

\bi
\item
given a pair of solutions $\theta = \left(\begin{array}{c}
\theta_1 \\ \theta_2 \end{array}\right)$ and $\tau =
\left(\begin{array}{c} \tau_1 \\ \tau_2 \end{array}\right)$ to
(\ref{zs1}), the Wronskian $W = \theta_1 \tau_2 - \theta_2 \tau_1$
is constant. In particular, we have
\begin{equation}
\label{zs21}
\det \Phi^\pm(x,k) = -1.
\end{equation}
\ei

In the sequel we assume that the potentials $p$ and $q$ are complex conjugate:
$$
p =\bar{q}.
$$

It is also checked by  straightforward computations that
the transformation
\begin{equation}
\label{zs3}
\psi = \left(\begin{array}{c} \psi_1 \\
\psi_2 \end{array}\right) \to \psi^\ast = \left(\begin{array}{c}
-\bar{\psi}_2
\\ \bar{\psi}_1
\end{array}\right)
\end{equation}
maps solutions to (\ref{zs1}) to solutions of the same equation.
In particular, it follows from the asymptotics of the Jost
functions that they are transformed as follows
\begin{equation}
\label{zs4}
\varphi^\pm_1 \stackrel{\ast}{\longrightarrow} -
\varphi^\pm_2, \ \ \ \varphi^\pm_2
\stackrel{\ast}{\longrightarrow} \varphi^\pm_1.
\end{equation}

Since the Jost functions $\varphi^+_l,l=1,2$, and
$\varphi^-_l,l=1,2$, give bases for the same space, they are related
by a linear transform
$$
\left(
\begin{array}{c}
\varphi^-_1 \\
\varphi^-_2
\end{array}
\right) = S(k)
\left(
\begin{array}{c}
\varphi^+_1 \\
\varphi^+_2
\end{array}
\right).
$$
It follows from (\ref{zs21}) that $\det S(k) = 1$ and we derive
from (\ref{zs4}) that

\bi
\item
{\it the scattering matrix} $S(k)$ is unitary: $S(k) \in SU(2)$,
i.e.
$$
S(k) =
\left(
\begin{array}{cc}
\overline{a(k)} & -\overline{b(k)} \\
b(k) & a(k)
\end{array}
\right), \ \ \ |a(k)|^2 + |b(k)|^2 = 1.
$$
\ei

The following quantities
$$
T(k) = \frac{1}{a(k)}, \ \ \ R(k) = \frac{b(k)}{a(k)}
$$
are called {\it the transmission coefficient} and {\it the
reflection coefficient} respectively. The operator $L$ is called
reflectionless if its reflection coefficient vanishes: $R(k) \equiv
0$.

The vector functions $\varphi_1^- e^{-ikx}$ and
$\varphi_2^+e^{ikx}$ are analytically continued onto the lower
half-plane $\Im k <0$, and the vector functions
$\varphi_2^-e^{ikx}$ and $\varphi_1^+e^{-ikx}$ are analytically
continued onto the upper half-plane $\Im k >0$.

Without loss of generality it is enough to prove that for
$\varphi_1^- e^{-ikx}$. This function satisfies the equation of
the Volterra type
$$
f(x,k) = \left(
\begin{array}{c} 0 \\ 1 \end{array}\right) -
\int_{-\infty}^x \left(\begin{array}{cc} 0 & -e^{-2ik(x-x^\prime)}
\\ 1 & 0 \end{array}\right)
\left(\begin{array}{cc} p & 0 \\ 0 & q \end{array}\right)
f(x^\prime,k)dx^\prime
$$
and since the integral kernel decays exponentially for $\Im k <0$
the Neumann series for its solution converges in this half-plane.

This implies that

\bi
\item
$T(k)$ is analytically continued onto the upper-half plane $\Im k
\geq 0$. \ei

It is shown that
\bi
\item
$a(k)$ vanishes nowhere on $\R \setminus \{0\}$;

\item
the poles of $T(k)$ correspond to {\it bounded states}, i.e., to
solutions to (\ref{zs1}) which decay exponentially as $x \rightarrow
\pm \infty$. These solutions are $\varphi^+_1(x,\varkappa)$ and
$\varphi^-_2(x,\varkappa)$ where $a(\varkappa) = 0$ and, therefore,
\beq \label{zs5} \varphi^-_2(x,\varkappa) = \mu(\varkappa)
\varphi^+_1(x,\varkappa), \ \ \ \mu(\varkappa) \in \C, \eeq and the
multiplicity of each eigenvalue $\varkappa$ equals to one;

\item
$T(k)$ has only simple poles in $\mbox{Im}\, k > 0$ and for
exponentially decaying potentials there are finitely many such
poles;

\item
since the set of solutions to (\ref{zs1}) is invariant under
(\ref{zs3}), the discrete spectrum of $L$ is preserved by the
complex conjugation $\varkappa \rightarrow \bar{\varkappa}$ and is
formed by the poles of $T(k)$ and their complex conjugates. \ei

{\it The spectral data} of $L$ consist of

1) the reflection coefficient $R(k), k \neq 0$;

2) the poles of $T(k)$ in the upper-half plane $\Im \varkappa >0$:
$\varkappa_1,\dots,\varkappa_N$;

3) the quantities  $\lambda_j=i\gamma_j \mu_j, j=1,\dots,N$, where
$\gamma_j = \gamma(\varkappa_j)$ is the residue of $T(k)$ at
$\varkappa_j$ and $\mu_j = \mu(\varkappa_j)$ (see (\ref{zs5})).

If the potential $p=\bar{q}$ is real-valued then
$$
\varphi^{\pm}_j(x,-k) = \overline{\varphi^{\pm}_j(x,k)} \ \ \
\mbox{for $k \in \R \setminus \{0\}$}
$$
and this implies that
$$
a(k) = \overline{a(-k)}, \ \ \ R(k) = \overline{R(-k)}, \ \ \ T(k)
= \overline{T(-\bar{k})},
$$
the poles of $T(k)$ are symmetric with respect to the imaginary
axis, and
$$
\lambda_j = \bar{\lambda}_k \ \ \ \mbox{for $\varkappa_j =
-\bar{\varkappa}_k$}.
$$

Now by applying the Fourier transform (with respect to $k$) to
both sides of the equality
\beq
\label{zs6}
T(k) \varphi^-_2 = R(k)\varphi^+_1
+ \varphi_2^+
\eeq
after some substitutions we write the equations (\ref{zs6}) for the
components of the vector functions in the form of the Gelfand--Levitan--Marchenko equations
$$
B_2(x,y) + \int_x^{+\infty} B_1(x,x')\Omega(x'+y) \, d x' = 0,
$$
$$
\Omega(x+y) - B_1(x,y) + \int_x^{+\infty} B_2(x,x') \Omega(x'+y)\, d x' = 0
$$
for $B_1$ and $B_2$ with
$$
\Omega(z) = \frac{1}{2\pi}\int_{-\infty}^{+\infty} R(k)e^{-i k z}
\, d k - \sum_{j=1}^N \lambda_j e^{i\varkappa_j z}
$$
where $y>x$ and there are the following limits
$$
\lim_{y \to \infty} B_k(x,y) = 0, \ \ \lim_{y \to x+} B_k(x,y) =
B_k(x,x), \ \ \ k=1,2,
$$
These equations are the Volterra type and are resolved uniquely.
The reconstruction formulas for the potentials are as follows:
\beq
\label{zs7}
p(x) = - 2 B_1(x,x), \ \ \
p(x)q(x) = p(x)\overline{p(x)} =
2 \frac{dB_2(x,x)}{dx}.
\eeq
See the detailed derivation of this formulas from \cite{ZS},
for instance, in \cite{Ablowitz}.

In the sequel, for simplicity,
we assume that the potential $p(x)$ is real-valued.

In \cite{Faddeev} a series of formulas expressing the Kruskal integrals in terms of the spectral data,
i.e. the so-called trace formulas,
is derived. We mention only the formula for the first nontrivial integral:
\beq
\label{zs8}
\int^\infty_{-\infty} p^2(x)dx = -\frac{1}{\pi} \int^\infty_{-\infty} \log (1-|b(k)|^2)dk + 4 \sum_{j=1}^N
\Im \varkappa_j.
\eeq

For reflectionless operators the reconstruction procedure reduces to
algebraic equations (see details, for instance, in
\cite{ZS,Ablowitz,T21}). The spectral data consist of the poles
$\varkappa_k$ and the corresponding quantities $\lambda_j$,
$j=1,\dots,N$. Put
$$
\Psi(x) = (-\lambda_1 e^{i\varkappa_1 x}, \dots, -\lambda_N e^{i\varkappa_N x}),
$$
$$
M_{j k}(x) =
\frac{\lambda_k}{i(\varkappa_j +\varkappa_k)}e^{i(\varkappa_j+\varkappa_k)x}, \ \ \
j,k=1,\dots,N.
$$
We have
\beq
\label{zs9}
\begin{split}
p(x) =
2 \frac{d}{d x}
\arctan \frac{\Im \det (1 + i M(x))}{\Re \det (1 + i M(x))},
\\
\varphi^+_1(x,k) =
\left(
\begin{array}{c}
\langle \Psi(x) \cdot (1+M^2(x))^{-1} | W(x,k) \rangle  \\
e^{i k x}  -
\langle \Psi(x) \cdot (1+M^2(x))^{-1} M(x) | W(x,k) \rangle
\end{array}
\right)
\end{split}
\eeq
where
$\langle u | v \rangle = u_1 v_1 + \dots + u_N v_N$ and
$$
W(x,k) =
\left(
\frac{i}{\varkappa_1 + k}e^{i ( \varkappa_1 + k)x}, \dots ,
\frac{i}{\varkappa_N + k}e^{i ( \varkappa_N + k)x}
\right).
$$

}

\medskip

{\bf Acknowledgement.}
This survey was started and a large part of it was written 
during author's stay at Max Plank Institute of Mathematics (MPIM) in 
October 2003 -- January 2004 and the final proofreadings were done during his 
stay at the same institute in December 2005.

\end{document}